\documentclass[a4paper,11pt,twoside]{article}
\usepackage{pst-fill,pst-grad}
\usepackage{textcomp}
\usepackage[english]{babel}
\usepackage[utf8x]{inputenc}
\usepackage{graphicx}
\usepackage{multicol}
\usepackage{float}
\usepackage{fancyhdr}
\usepackage[matrix,arrow,curve]{xy}
\usepackage{pstricks} 
\usepackage{amsmath,amsfonts,verbatim,afterpage,theorem,euscript,mathrsfs,amssymb}
\usepackage{array}
\usepackage{dsfont}
\usepackage{hyperref}

\numberwithin{equation}{section}

\def \vu{\vec{u}}

\def \vw{\vec{\omega}}
\def \va{\vec{a}}

\def \vU{\vec{U}}
\def \vB{\vec{B}}

\def \vUU{\vec{\mathcal{U}}}

\def \vn{\vec{\nabla}}
\def \vb{\vec{b}}
\def \vf{\vec{f}}
\def \vg{\vec{g}}
\def \vF{\vec{F}}
\def \vG{\vec{G}}

\def \vphi{\vec{\varphi}}

\def \vpphi{\vec{\phi}}

\def \grad{\vec{\nabla}}
\def \rot{\vec{\nabla}\wedge}

\def \p{\mathfrak{p}}
\def \q{\mathfrak{q}}
\def \R{\mathbb{R}^{3}}
\def \M{\mathcal{M}_{t,x}}

\def \Ao{\mathbf{A}_{\rho}}

\def \Poq{\mathbf{P}_{\rho}}

\def \TTo{\mathbf{\Theta}_{\rho}}

\newtheorem{defi}{Definition}[section]
\newtheorem{propo}{Proposition}[section]
\newtheorem{lem}{Lemma}[section]
\newtheorem{theo}{Theorem}[section]
\newtheorem{coro}{Corollary}[section]
\newtheorem{rema}{Remark}[section]
\title{\bf Interior $\epsilon$-regularity theory for the solutions of the magneto-micropolar equations with a perturbation term}
\usepackage{authblk}
\author{Diego Chamorro\footnote{\emph{diego.chamorro@univ-evry.fr}} }
\author{David Llerena}
\affil{\footnotesize LaMME, Univ. Evry, CNRS, Universit\'e Paris-Saclay, 91037, Evry, France.}

\begin{document}
\sloppy
\maketitle
\begin{scriptsize}
\abstract{We develop here a particular version of the partial regularity theory for the Magneto-Micropolar equations (MMP) where a perturbation term is added. These equations are used in some special cases, such as in the study of the evolution of liquid cristals or polymers, where the classical Navier–Stokes equations are not an accurate enough model. The incompressible Magneto-Micropolar system is composed of three coupled equations: the first one is based in the Navier-Stokes system, the second one considers mainly the magnetic field while the last equation introduces the microrotation field representing the angular velocity of the rotation of the fluid particles. External forces are considered and a specific perturbation term is added as it is quite useful in some applications.}\\

\end{scriptsize}

%\tableofcontents
%%%%%%%%%%%%%%%%%%%%%%%%%%%%%%%%%%%%%%
\section{Introduction}
Micropolar equations were introduced in 1966 by Eringen \cite{Eringen} and were first studied mathematically in 1997 by Galdi \& Rionero \cite{Galdi}. 
Some very recent results concerning the regularity of the solution to this system were obtained in \cite{GuWang, LoMelo} (see also the references there in). In this article we will consider a slightly more general framework by introducing a magnetic field, some external forces and a perturbation term. We will develop here the $\epsilon$-regularity criterion which was not, to the best of our knowledge, treated before for this type of problem. The incompressible 3D-Magneto-Micropolar (MMP) system studied in this article is of the following form
\begin{equation*}
\begin{cases}
\partial_t \vU = \Delta \vU  -(\vU \cdot \grad)\vU  +(\vB\cdot \grad)\vB-\grad p +\frac{1}{2}\rot \vw+ \vF-div(\vU\otimes \va +\va \otimes \vU),\\[3mm]
\partial_t \vB = \Delta \vB -(\vU \cdot \grad)\vB +(\vB\cdot \grad)\vU 
+ \vG,\\[3mm]
\partial_t \vw = \Delta \vw +\grad div(\vw)-\vw -(\vU \cdot \grad)\vw
+\frac{1}{2}\rot	\vU ,\\[3mm]
div(\vU)=div(\vF)=div(\vB)=div(\vG)=div(\va)=0,\\[3mm]
\vU (0,x)=\vU_0(x), \;\; \vB(0,x)=\vB_0(x),\;\;\vw(0,x)= \vw_0(x),\;\; x\in \R
\mbox{ and  }div(\vU _0)=div(\vB_0)=0.
\end{cases}
 \end{equation*}
Here $\vU$ denotes the fluid velocity field, $\vB$ is the magnetic field, $\vw$ is the field of microrotation representing the angular velocity of the rotation of the fluid particles and $p$ is the scalar pressure. The quantities $\vF$ and $\vG$ represent external forces (assumed divergence free) and they are given as well as the initial data $\vU_0, \vB_0$ and $\vw_0$.\\

The perturbation $\va$ which appears in the first equation above in the term $div(\vU\otimes \va +\va \otimes \vU)$ is a given divergence free vector field and the presence of this particular type of perturbation is mainly inspired by quantitative studies for the rate of possible blow-up for the Navier-Stokes equations (see in particular the article \cite{BP2}), see also the book \cite[Section 12.6]{PGLR1} for other interesting applications of this type of perturbation. As pointed out in the Remark \ref{RemarkHipotesis} below, the assumptions over $\va$ will have some impact in the general set of hypotheses needed in order to perform our computations.\\ 

Now, in order to simplify the computations we introduce the Elsasser formulation, which was initially used for the Magnetohydrodynamics equations (MHD) see \cite{elsasser1950}: indeed, by a suitable change of variables we will obtain a more symmetric problem and if we define $\vu=\vU+\vB$, $\vb=\vU-\vB$, $\vf=\vF+\vG$ and $\vg=\vF-\vG$, then for all $ x\in \R$ we can write
\begin{equation}
\begin{cases}\label{EquationMMP}
\partial_t \vu  = \Delta \vu  -(\vb \cdot \grad)\vu-\grad p +\frac{1}{2}\rot\vw+ \vf-\frac{1}{2}div((\vu+\vb) \otimes \va+\va \otimes (\vu+\vb)),\\[3mm]
\partial_t \vb = \Delta \vb -(\vu \cdot \grad)\vb -\grad p +\frac{1}{2}\rot\vw +\vg-\frac{1}{2}div((\vu+\vb) \otimes \va +\va \otimes (\vu+\vb)),\\[3mm]
\partial_t \vw = \Delta \vw +\grad div(\vw)-\vw-\frac{1}{2}((\vu+\vb)\cdot \grad)\vw+\frac{1}{4}\rot(\vu+\vb),\\[3mm]
div(\vu)=div(\vb)=div(\vf)=div(\vg)=div(\va)=0,\\[3mm]
\vu (0,x)=\vu_0(x), \;\; \vb(0,x)=\vb_0(x),\;\;\vw(0,x)=\vw_0(x)
\mbox{ and } div(\vu _0)=div(\vb_0)=0.
\end{cases}
\end{equation}
It is worth to remark here that as long as we want to perform a generic study for the functions $\vu$ and $\vb$, this previous system presents a simpler framework and thus, for the rest of the article we will focus ourselves in this formulation. We remark also that since $div(\vu)=div(\vb)=0$, then we can deduce from (\ref{EquationMMP}) that the pressure $p$ satisfies the equation
\begin{equation}\label{FormulePressionIntro}
2\Delta p=-div((\vb\cdot \vn)\vu)-div((\vu\cdot \vn)\vb)-div(div((\vu+\vb) \otimes \va+\va \otimes (\vu+\vb))),
\end{equation}
and we see from this expression that the pressure $p$ is only determined by the couple $(\vu, \vb)$ (recall that $\va$ is given) and we will see how to exploit this relationship later on. \\

We are interested here in studying some properties of (local) weak solutions of the system (\ref{EquationMMP}) and in order to fix the notation we consider now $\Omega$ a bounded subset of $]0, +\infty[\times \mathbb{R}^3$ of the form
\begin{equation}\label{DefConjuntoOmega}
\Omega=]a,b[\times B(x_{0},r), \quad \mbox{with} \quad 0<a<b<+\infty,\; x_{0}\in \R \mbox{ and } 0<r<+\infty. 
\end{equation}
and we will say that $(\vu, \vb, \vw)\in L^\infty_tL^2_x\cap L^2_t\dot{H}^1_x(\Omega)$ satisfies the MMP equations (\ref{EquationMMP}) in the weak sense if for all $\vphi, \vpphi, \vec{\psi} \in \mathcal{D}(\Omega)$ such that $div(\vphi)=div(\vpphi)=0$, we have
$$
\begin{cases}
\langle\partial_{t}\vu-\Delta \vu +(\vb\cdot\vn)\vu-\frac{1}{2}\rot\vw-\vf+\frac{1}{2}div((\vu+\vb) \otimes \va+\va \otimes (\vu+\vb))|\vphi \rangle_{\mathcal{D}'\times \mathcal{D}}=0,\\[3mm]
\langle\partial_{t}\vb-\Delta \vb +(\vu\cdot\vn)\vb-\frac{1}{2}\rot\vw -\vg+\frac{1}{2}div((\vu+\vb) \otimes \va +\va \otimes (\vu+\vb))|\vpphi \rangle_{\mathcal{D}'\times \mathcal{D}}=0,\\[3mm]
\langle\partial_t \vw - \Delta \vw -\grad div(\vw)+\vw+\frac{1}{2}((\vu+\vb)\cdot \grad)\vw-\frac{1}{4}\rot(\vu+\vb)|\vec{\psi}\rangle_{\mathcal{D}'\times \mathcal{D}}=0.
\end{cases}
$$
Note that if $(\vu, \vb, \vw)$ are solutions of the previous system, then due to the expression (\ref{FormulePressionIntro}) there exists a pressure $p$ such that (\ref{EquationMMP}) is fulfilled in $\mathcal{D}'$. \\

Based in the classical Navier-Stokes problem, we can study at least two main regularity theories for the MMP equations: the \emph{local} regularity theory (also known as the Serrin criterion, see \cite{OLeary}, \cite{Serrin1}) and the $\epsilon$-regularity criterion (also known as the \emph{partial} regularity theory, based in the seminal work of Caffarelli, Kohn and Nirenberg \cite{CKN}, see also \cite{MHDKukavica} and \cite{Kukavica, Kukavica0, Kukavica1}).\\

As said previously, in this article we want to develop a particular version of the $\epsilon$-regularity criterion for the system (\ref{EquationMMP}) and we need to impose some assumptions over the functions $\vu, \vb$ and $\vw$ as well as some hypothesis over pressure $p$ and from now on we will always assume that we have the following controls
\begin{align}\label{HypothesesLocal1}
\begin{aligned}
&\vu, \vb, \vw \in L_t^{\infty}L_x^2\cap L_t^2 \dot H_x^1(\Omega),\\
&p\in L_{t,x}^{\frac{3}{2}}(\Omega)\cap L^{\frac52}_tL^1_x(\Omega),\qquad \va\in L_{t,x}^6(\Omega),\qquad \vw\in L_{t,x}^\infty(\Omega),\\
&\vf, \vg \in L_{t,x}^{\frac{10}{7}}(\Omega),
\end{aligned}
\end{align}
where $\Omega$ is a subset of $\mathbb{R}\times \mathbb{R}^3$ of the form (\ref{DefConjuntoOmega}).
\begin{rema}\label{RemarkHipotesis}
The conditions over $\vu, \vb, \vw$ and $\vf, \vg$ are rather classical in the setting of equations arising from fluid dynamics.  Note that the connection between the perturbation term $\va$ and the pressure $p$ is explicitly given in the relationship (\ref{FormulePressionIntro}) above, thus if we assume the local integrability condition $\va\in L_{t,x}^6(\Omega)$ (which appears naturally in some recent results, see \cite{BP2}), then following our computations we need to impose the condition $p\in L_{t,x}^{\frac{3}{2}}(\Omega)\cap L^{\frac52}_tL^1_x(\Omega)$ for the pressure. Observe that conditions of the form $ L^{q}_tL^1_x(\Omega)$ for the pressure were also considered in the setting of the Navier-Stokes equations, see \cite{Vasseur}. Finally note that the $(L_{t,x}^\infty)_{loc}$ information is usually asked in regularity theory, but in this work we only assume it for the variable $\vw$ (not for $\vu$ nor for $\vb$) and this will be crucial to study the term $\grad div(\vw)$ which appears in the micropolar equation (\ref{EquationMMP}). See also Remark \ref{Remarque_HypothesesVW} below, where alternative and more general assumptions are given for the variable $\vw$.
\end{rema}
\begin{rema}
We do not claim here any optimality on the space $L_{t,x}^6(\Omega)$ related to the perturbation term and we believe that it is perhaps possible to consider a slightly more general perturbation term by asking $\va\in L_{t,x}^m(\Omega)$ for $m\geq 5$, however, as far as we can see, this will introduce some quite difficult technical problems and will probably induce some extra hypotheses over the pressure. On the other hand, if we assume some additional information (say $\va\in L_{t}^2\dot{H}^1_x(\Omega)$), then  we can relax the hypotheses on the pressure and work only with $p\in L_{t,x}^{\frac{3}{2}}(\Omega)$.
\end{rema}

Once this local framework is clear, we can now introduce a special class of weak solutions:
\begin{defi}[Suitable solutions]\label{Def_SuitableSolutions} 
Let $(\vu, \vb, \vw, p)$ be a weak solution over $\Omega$ for the perturbed magneto-micropolar equations (\ref{EquationMMP}) which satisfies the local hypotheses (\ref{HypothesesLocal1}) stated above. 
We say that $(\vu, \vb, \vw, p)$ is a suitable solution if the distribution $\mu$ given by the expression
\begin{eqnarray}
\mu&=&-\partial_t(|\vu|^2+|\vb|^2+|\vw|^2)+\Delta(|\vu|^2+|\vb|^2+|\vw|^2)-2(|\grad\otimes\vu|^2+|\grad\otimes\vb|^2+|\grad\otimes\vw|^2)\notag\\
&&-div\left((|\vb|^2+2p)\vu+(|\vu|^2+2p)\vb+\frac{1}{2}(|\vu|^2+|\vb|^2)\vw\right)+2\vn div(\vw)\cdot \vw-2|\vw|^2\notag\\
&&+(\rot \vw)\cdot(\vu+\vb)+\frac{1}{2}(\rot(\vu+\vb))\cdot\vw+2(\vf+\vg)\cdot(\vu+\vb)\label{Def_Mu}\\
&&+div((\vu+\vb) \otimes \va+\va \otimes (\vu+\vb))\cdot (\vu+\vb),\notag
\end{eqnarray}
is a non-negative locally finite measure on $\Omega$.
\end{defi}
\begin{rema}
It is worth noting here that the local hypotheses stated in (\ref{HypothesesLocal1}) guarantee that each one of the terms in the previous expression is meaningful. 
\end{rema}
%%%%%%%%%%%%%%%%%%%%%%%%%%%%%%%%%%%%%%
The main purpose of this article is to prove the following theorem which gives a gain of regularity in space and time variables for suitable solutions.
%%%%%%%%%%%%%%%%%%%%%%%%%%%%%%%%%%%
\begin{theo}\label{HolderRegularity_theorem}
Let $\Omega$ be a subset of the form (\ref{DefConjuntoOmega}). 
Let $(\vu, \vb, \vw, p)$ be a weak solution on $\Omega$ of the magneto-micropolar equations (\ref{EquationMMP}). Assume that
\begin{enumerate}
\item[1)] $(\vu, \vb, \vw, p, \vf, \vg, \va)$ satisfies the conditions (\ref{HypothesesLocal1}),\\
\item[2)] $(\vu, \vb, \vw, p)$ is suitable in the sense of Definition \ref{Def_SuitableSolutions},\\
\item[3)] we have the following local information on $\vf$ and $\vg$: 
$\mathds{1}_{\Omega}\vf \in \M^{\frac{10}{7}, \tau_a}$ and 
$\mathds{m1}_{\Omega}\vg\in \M^{\frac{10}{7}, \tau_b}$ for some indexes
$\tau_a,\tau_b>\frac{5}{2-\alpha}$ with $0<\alpha<\frac{1}{12}$.
\end{enumerate}
There exists a positive constant $\epsilon^*$ which depends only  on $\tau_a$ and $\tau_b$ such that, if for some $(t_0,x_0)\in \Omega$, we have
\begin{equation*}
\limsup_{r\to 0}\frac{1}{r}\iint _{]t_0-r^2, t_0+r^2[\times B(x_0,r)}|\grad \otimes \vu|^2+|\grad \otimes \vb|^2+|\grad \otimes \vw|^2dx ds<\epsilon^*,
\end{equation*}
then $(\vu,\vb,\vw)$ is H\"older regular (in the time and space variables) of exponent $\alpha$ in a neighborhood of $(t_0,x_0)$ for some small $\alpha$ in the interval $0<\alpha<\frac{1}{12}$.
\end{theo}
%%%%%%%%%%%%%%%%%%%%%%%%%%%%%%%%%%%%%%
Some remarks are in order here. 
\begin{itemize}
\item Following standard procedures it is possible to construct Leray-type weak solutions for the problem (\ref{EquationMMP}). However we are only interested here to study the local behavior (for some points of the subset $\Omega$) of the solutions of such equations. 
\item The hypothesis over the pressure $p$ (\emph{i.e.} $p\in L_{t,x}^{\frac{3}{2}}(\Omega)$) is useful to give a sense to the quantities $div(p\vu)$ and $div(p\vb)$ that are present in the definition of the measure $\mu$ given in (\ref{Def_Mu}). It is worth noting that in the setting of the classical Navier-Stokes equation this hypothesis can be removed and a generic pressure $p\in \mathcal{D}'$ can be considered. See \cite{CML} for the details.
\item Some additional hypothesis over the external forces $\vf$ and $\vg$ are stated in Morrey spaces. We will see in the computations below that this functional framework is particularly well suited to the study of the regularity for this type of equations.\\
\end{itemize}

The plan of the article is the following: in Section \ref{Secc_MorreyIntro}  we recall some notation and useful facts about our framework. In Section \ref{Secc_ProofMaintheorem} we establish a first gain of regularity under some particular hypotheses stated in terms of Morrey spaces. The rest of the article (Sections \ref{Secc_LocalEnergy}, \ref{Secc_Inductive} and \ref{Secc_MasEstimates}) is devoted to the proof of these hypotheses.
%%%%%%%%%%%%%%%%%%%%%%%%%%%%%%%%%%%%%%
\section{Notation and functional spaces}\label{Secc_MorreyIntro}
For $1\leq p, q\leq +\infty$ we characterize the Lebesgue space $L^p([0,+\infty[, L^q(\R))$ as the set of measurable functions $\vf:[0,+\infty[\times \R\longrightarrow \R$ such that $\displaystyle{\|\vf\|_{L^p_tL^q_x}=\left(\int_{0}^{+\infty}\|\vf(t,\cdot)\|_{L^q}dt\right)^{\frac1p}<+\infty}$ with the usual modifications when $p=+\infty$ or $q=+\infty$. We also define the space $L^p([0,+\infty[, \dot{H}^s(\R))$ with $1\leq p\leq +\infty$ and $s>0$ as the set of distributions such that $\displaystyle{\|\vf\|_{L^p_t\dot{H}^s_x}=\left(\int_{0}^{+\infty}\|\vf(t,\cdot)\|_{\dot{H}^s}dt\right)^{\frac1p}<+\infty}$
where $\dot{H}^s(\R)$ is the usual homogeneous Sobolev space. See the books \cite{PGLR0} and \cite{PGLR1} for details about these functional spaces.\\

We recall now the notions of parabolic  H\"older and Morrey spaces and for this we need first to consider the homogeneous space $(\mathbb{R}\times \R, d, \mu)$ where $d$ is the parabolic (quasi)distance given by 
$d\big((t,x), (s,y)\big)=|t-s|^{\frac{1}{2}}+|x-y|$ and where $\mu$ is the usual Lebesgue measure $d\mu=dtdx$. Associated to this distance, we  define homogeneous (parabolic) H\"older spaces $\dot{\mathcal{C}}^\alpha(\mathbb{R}\times \R, \R)$ where $0<\alpha<1$ by the usual condition:
\begin{equation}\label{Holderparabolic}
\|\vphi\|_{\dot{\mathcal{C}}^\alpha}=\underset{(t,x)\neq (s,y)}{\sup}\frac{|\vphi(t,x)-\vphi(s,y)|}{\left(|t-s|^{\frac{1}{2}}+|x-y|\right)^\alpha}<+\infty,
\end{equation}
and this formula studies H\"older regularity in both time and space variables. Now, for $1< p\leq q<+\infty$, parabolic Morrey spaces $\mathcal{M}_{t,x}^{p,q}$ are defined as the set of measurable functions $\vphi:\mathbb{R}\times\R\longrightarrow \R$ that belong to the space $(L^p_{t,x})_{loc}$ such that $\|\vphi\|_{M_{t,x}^{p,q}}<+\infty$ where
\begin{equation}\label{DefMorreyparabolico}
\|\vphi\|_{\mathcal{M}_{t,x}^{p,q}}=\underset{x_{0}\in \R, t_{0}\in \mathbb{R}, r>0}{\sup}\left(\frac{1}{r^{5(1-\frac{p}{q})}}\int_{|t-t_{0}|<r^{2}}\int_{B(x_{0},r)}|\vphi(t,x)|^{p}dxdt\right)^{\frac{1}{p}}.
\end{equation}
These spaces are generalization of usual Lebesgue spaces, note in particular that we have $\mathcal{M}_{t,x}^{p,p}=L_{t,x}^p$. See \cite{Adams} for more details on these spaces. We refer the readers to the book \cite{PGLR1} for a general theory concerning the Morrey spaces and Hölder continuity applied to the analysis of PDEs from fluid mechanics. Here are some useful fact concerning Morrey spaces:
\begin{lem}[H\"older inequalities]\label{lemma_Product}
\begin{itemize}
\item[]
\item[1)]If $\vf, \vg:\mathbb{R} \times \R\longrightarrow \R$ are two functions such that $\vf\in \M^{p,q} (\mathbb{R} \times \R)$ and $\vg\in L^{\infty}_{t,x} (\mathbb{R} \times \R)$, then for all $1\leq p\leq q<+\infty$ we have
$\|\vf\cdot\vg\|_{\M^{p,q}}\leq  C\|\vf\|_{\M^{p, q}} \|\vg\|_{L^{\infty}_{t,x}}$.
\item[2)]If $\vf, \vg:\mathbb{R} \times \R\longrightarrow \R$ are two functions that belong to the space $\M^{p,q} (\mathbb{R} \times \R)$ then we have the inequality $\|\vf\cdot\vg\|_{\M^{\frac{p}{2}, \frac{q}{2}}}\leq  C\|\vf\|_{\M^{p, q}} \|\vg\|_{\M^{p, q}}$.
\item[3)] More generally, let $1\leq p_0 \leq q_0 <+\infty$, $1\leq p_1\leq q_1<+\infty$ and $1\leq p_2\leq q_2<+\infty$. If $\tfrac{1}{p_1}+\tfrac{1}{p_2}\leq \frac{1}{p_0}$ and $\tfrac{1}{q_1}+\tfrac{1}{q_2}=\tfrac{1}{q_0}$, then for two measurable functions $\vf, \vg:\mathbb{R} \times \R\longrightarrow \R$ such that $\vf\in \mathcal{M}^{p_1, q_1}_{t,x}$ and $\vg\in \mathcal{M}^{p_2, q_2}_{t,x}$, we have the following H\"older inequality in Morrey spaces $\|\vf\cdot \vg\|_{\mathcal{M}^{p_0, q_0}_{t,x}}\leq \|\vf\|_{\mathcal{M}^{p_1, q_1}_{t,x}}\|\vg\|_{\mathcal{M}^{p_2, q_2}_{t,x}}$.
\end{itemize}
\end{lem} 
\begin{lem}[Localization]\label{lemma_locindi}
Let $\Omega$ be a bounded set of $\mathbb{R} \times \R$. If we have $1\leq p_0 \leq p_1$, $1\leq p_0\leq q_0 \leq q_1<+\infty$ and if the function $\vf:\mathbb{R} \times \R\longrightarrow \R$  belongs to the space $\M^{p_1,q_1} (\mathbb{R} \times \R)$ then we have the following localization property $\|\mathds{1}_{\Omega}\vf\|_{\M^{p_0, q_0}} \leq C\|\mathds{1}_{\Omega}\vf\|_{\M^{p_1,q_1}}\leq C\|\vf\|_{\M^{p_1,q_1}}$.
\end{lem} 
%%%%%%%%%%%%%%%%%%%%%%%%%%%%%%%%%%%%%%
\section{A parabolic gain of regularity: the first step}\label{Secc_ProofMaintheorem}

The proof of Theorem \ref{HolderRegularity_theorem} is essentially based in the following regularity result for parabolic equations which is stated here in the framework of (parabolic) Morrey spaces:
%%%%%%%%%%%%%%%%%%%%%%%%%%%%%%%%%%%%%
\begin{propo}\label{parabolicHolder}
For $\vec{\Phi}, \vec{\Psi}: [0,+\infty[\times\R \longrightarrow \R$ two vector fields, we consider the following equation
\begin{equation}\label{parabolicHolder1}
\begin{cases}
\partial_t \vec{v}=\Delta \vec {v}+\vec {\Phi}+\sigma(D) \vec{\Psi} ,\\[2mm]
\vec v(0,x)=0,
\end{cases}
\end{equation}
where $\sigma$ is a smooth function on $\R\setminus \{0\}$, homogeneous of exponent\footnote{\emph{i.e.} $\sigma(\lambda \xi)=\lambda \sigma(\xi)$ for all $\lambda>0$.} 1 and 
$\sigma(D)$ is the Fourier multiplier operator of symbol $\sigma$ (acting component-wise).\\

Assume that we have $\vec \Phi \in \mathcal{M}_{t,x}^{\p_0,\q_0}$ and  $\vec{\Psi}\in \mathcal{M}_{t,x}^{\p_0,\q_1}$ with $1< \p_0\le \q_0<\q_1$ where we have $\frac{1}{\q_0}=\frac{2-\alpha}{5}$,
$\frac{1}{\q_1}=\frac{1-\alpha}{5}$, and $0<\alpha<1$.
Then the function $\vec {v} $ equal to $0$ for $t\le 0$ 
and to
\begin{equation*}
\vec v(t,x)=\int_0^t e^{(t-s)\Delta }(\vec \Phi (s,\cdot)+\sigma(D) \vec{\Psi}(s,\cdot))ds \qquad \mbox{for } t>0,
\end{equation*}
is a solution of equation (\ref{parabolicHolder1}) that is Hölderian of exponent $\alpha$ in the sense of (\ref{Holderparabolic}). 
\end{propo}
See \cite[Proposition 13.4]{PGLR1} for a proof of this result, see also \cite{Ladyzhenskaya1}.\\

We will apply this proposition to our system (\ref{EquationMMP}) but, as we only assume the controls (\ref{HypothesesLocal1}) over a subset $\Omega$ of the form (\ref{DefConjuntoOmega}), we need to localize our framework and for this we first fix the point $(t_0,x_0)\in \Omega$ considered in the hypotheses of Theorem \ref{HolderRegularity_theorem} and then for a small enough radius $0<\mathfrak{r}<1$, we consider the parabolic ball
\begin{equation}\label{Def_BoulesQ}
Q_\mathfrak{r}(t_0,x_0)=]t_0-\mathfrak{r}^2,t_0+\mathfrak{r}^2[\times B(x_0,\mathfrak{r}),
\end{equation}
such that $Q_{5\mathfrak{r}}(t_0,x_0)\subset \Omega$ (these parabolic balls will be denoted by $Q_\mathfrak{r}$ for simplicity). Note here that since by  (\ref{DefConjuntoOmega}) we have $\Omega=]a,b[\times B(x_{0},r)$ with  $0<a<b<+\infty$ and $x_{0}\in \R$, then the condition $Q_{5\mathfrak{r}}(t_0,x_0)\subset \Omega$ guarantees the fact that $t_0-\mathfrak{r}^2>0$ and thus the time interval $]t_0-\mathfrak{r}^2,t_0+\mathfrak{r}^2[$ does not contain the origin: this condition is important in order to obtain a system of the form (\ref{parabolicHolder1}) for which the initial data is such that $\vec v(0,x)=0$. Now, we construct an auxiliary non-negative function $\phi: \mathbb{R}\times \R\longrightarrow \mathbb{R}$ such that $\phi\in \mathcal{C}_0^{\infty}(\mathbb{R}\times \R)$, $supp(\phi)\subset ]-\tfrac{1}{16},\tfrac{1}{16}[\times B(0,\tfrac{1}{4})$ and such that $\phi \equiv 1\;\; \text{on}\; \; ]-\tfrac{1}{64},\tfrac{1}{64}[\times B(0,\tfrac{1}{8})$ and for a fixed $\mathbf{R}>0$ such that 
\begin{equation}\label{Def_R}
4\mathbf{R}<\mathfrak{r}, 
\end{equation}
we define the localizing function $\eta$ by $\eta(t,x)=\phi\left(\frac{t-t_0}{\mathbf{R}^2},\frac{x-x_0}{\mathbf{R}}\right)$ (remark that we have $supp \; \eta\subset Q_{\mathbf{R}}$) and we define the vector $\vUU=\eta(\vu+\vb+\vw)$. As we can observe, we have the identity $\vec {\mathcal{U}}=\vu+\vb+\vw$
over a small neighborhood of the point $(t_0,x_0)$ and the support of the
variable $\vec{\mathcal{U}}$ is contained in the parabolic ball
$Q_{\mathbf{R}}(t_0,x_0)\subset Q_{\mathfrak{r}}(t_0,x_0)\subset \Omega$. Moreover, this localization forces the property $\vUU(0,\cdot)=0$, we can thus consider the following problem:
\begin{equation}\label{EquationPourUtilisationMorrey}
\begin{cases}
\partial_t \vUU=\Delta \vUU+\vec{\mathcal{B}}+\vn \beta -div(\mathbb{B}),\\[3mm]
\vUU(0,x)=0,
\end{cases}
\end{equation}
where the vector $\vec{\mathcal{B}}$ is given by
\begin{eqnarray}
\vec{\mathcal{B}}&=&(\partial_t\eta- \Delta \eta)(\vu+\vb+\vw)-2\sum_{i=1}^3 (\partial_i \eta)(\partial_i (\vu+\vb+\vw))-\eta\bigg((\vb \cdot \grad)\vu+(\vu \cdot \grad)\vb\bigg)\label{BMorrey}\\
&+&(\vn\eta) \bigg(\frac{div }{(-\Delta)}div\left(\vb\otimes \vu+\vu\otimes \vb+(\vu+\vb)\otimes \va+\va\otimes(\vu+\vb)\right)\bigg)+\eta(\rot\vw)+\eta(\vf+\vg)\notag\\
&+&\bigg((\vu+\vb) \otimes \va+\va \otimes (\vu+\vb)\bigg)\cdot \vn\eta- (\vn \eta) div(\vw)-\eta \left(\vw+\frac{1}{2}((\vu+\vb)\cdot \grad)\vw+\frac{1}{4}\rot(\vu+\vb)\right),\notag
\end{eqnarray}
the scalar function $\beta$ is given by\\
\begin{equation}\label{BetaMorrey}
\beta=\eta div(\vw)-\eta \frac{div }{(-\Delta)}div\left(\vb\otimes \vu+\vu\otimes \vb+(\vu+\vb)\otimes \va+\va\otimes(\vu+\vb)\right),
\end{equation}
and the tensor $\mathbb{B}$ is given by
\begin{equation}\label{TensorBMorrey}
\mathbb{B}=\eta( \va\otimes (\vu+\vb)+(\vu+\vb)\otimes\va).
\end{equation}
Indeed, in order to verify that we have the equation (\ref{EquationPourUtilisationMorrey}) with the terms (\ref{BMorrey}), (\ref{BetaMorrey}) and (\ref{TensorBMorrey}) above, we compute $\partial_t\vUU$ and we have 
\begin{eqnarray*}
\partial_t\vUU&=&(\partial_t\eta)(\vu+\vb+\vw)+\eta \Delta (\vu+\vb+\vw)-\eta\bigg((\vb \cdot \grad)\vu+(\vu \cdot \grad)\vb\bigg)-2\eta\grad p+\eta(\rot\vw) \\
&+&\eta(\vf+\vg)- \eta\left(div((\vu+\vb) \otimes \va+\va \otimes (\vu+\vb))\right)+ \eta \left(\grad div(\vw)-\vw-\frac{1}{2}((\vu+\vb)\cdot \grad)\vw+\frac{1}{4}\rot(\vu+\vb)\right).
\end{eqnarray*}
We use now the identity 
\begin{eqnarray*}
\eta \Delta (\vu+\vb+\vw)&=&\Delta (\eta (\vu+\vb+\vw)) - \Delta \eta(\vu+\vb+\vw)-2\sum_{i=1}^3 (\partial_i \eta)(\partial_i (\vu+\vb+\vw))\\
&=&\Delta\vUU- \Delta \eta(\vu+\vb+\vw)-2\sum_{i=1}^3 (\partial_i \eta)(\partial_i (\vu+\vb+\vw)),
\end{eqnarray*}
to obtain the expression
\begin{eqnarray}
\partial_t\vUU&=&\Delta\vUU+(\partial_t\eta- \Delta \eta)(\vu+\vb+\vw)-2\sum_{i=1}^3 (\partial_i \eta)(\partial_i (\vu+\vb+\vw))\label{SystemeLocalise}\\
&&-\eta\bigg((\vb \cdot \grad)\vu+(\vu \cdot \grad)\vb\bigg)-2\eta\grad p+\eta(\rot\vw) +\eta(\vf+\vg)- \eta\left(div((\vu+\vb) \otimes \va+\va \otimes (\vu+\vb))\right)\notag\\
&&+ \eta \left(\grad div(\vw)-\vw-\frac{1}{2}((\vu+\vb)\cdot \grad)\vw+\frac{1}{4}\rot(\vu+\vb)\right),\notag
\end{eqnarray}
which is the first step to obtain an equation of the form (\ref{EquationPourUtilisationMorrey}). We need now to organize the expression above in a suitable manner and for this we need to rewrite three particular terms, indeed, since we have the identities $\eta (\grad p)=\grad(\eta p)-(\grad \eta) p$, $\eta div(\va\otimes \vu)= div(\eta (\va\otimes \vu) )- (\va\otimes \vu)\cdot \grad \eta$ and $\eta \vn div(\vw)=\vn(\eta div(\vw))-(\vn \eta) div(\vw)$, we obtain 
\begin{eqnarray*}
\partial_t\vUU&=&\Delta\vUU+(\partial_t\eta- \Delta \eta)(\vu+\vb+\vw)-2\sum_{i=1}^3 (\partial_i \eta)(\partial_i (\vu+\vb+\vw))-\eta\bigg((\vb \cdot \grad)\vu+(\vu \cdot \grad)\vb\bigg)+2(\grad \eta) p\\
&+&\eta(\rot\vw) +\eta(\vf+\vg)+\bigg((\vu+\vb) \otimes \va+\va \otimes (\vu+\vb)\bigg)\cdot \vn\eta- (\vn \eta) div(\vw)\\
&-&\eta \left(\vw+\frac{1}{2}((\vu+\vb)\cdot \grad)\vw+\frac{1}{4}\rot(\vu+\vb)\right)\\
&-& 2\grad(\eta p)+\vn(\eta div(\vw))- div\bigg(\eta\big((\vu+\vb) \otimes \va+\va \otimes (\vu+\vb)\big)\bigg).
\end{eqnarray*}
We recall now that, from the expression (\ref{FormulePressionIntro}) and using the fact that $div(\vu)=div(\vb)=0$, we have the following identity for the pressure $p=\frac{div }{2(-\Delta)}div\left(\vb\otimes \vu+\vu\otimes \vb+(\vu+\vb)\otimes \va+\va\otimes(\vu+\vb)\right)$ so we can finally write
\begin{eqnarray*}
\partial_t\vUU&=&\Delta\vUU+(\partial_t\eta- \Delta \eta)(\vu+\vb+\vw)-2\sum_{i=1}^3 (\partial_i \eta)(\partial_i (\vu+\vb+\vw))-\eta\bigg((\vb \cdot \grad)\vu+(\vu \cdot \grad)\vb\bigg)\\
&+&(\vn\eta) \bigg(\frac{div }{(-\Delta)}div\left(\vb\otimes \vu+\vu\otimes \vb+(\vu+\vb)\otimes \va+\va\otimes(\vu+\vb)\right)\bigg)+\eta(\rot\vw)+\eta(\vf+\vg)\\
&+&\bigg((\vu+\vb) \otimes \va+\va \otimes (\vu+\vb)\bigg)\cdot \vn\eta- (\vn \eta) div(\vw)-\eta \left(\vw+\frac{1}{2}((\vu+\vb)\cdot \grad)\vw+\frac{1}{4}\rot(\vu+\vb)\right)\\
&+&\vn\bigg(\eta div(\vw)-\eta \frac{div }{(-\Delta)}div\left(\vb\otimes \vu+\vu\otimes \vb+(\vu+\vb)\otimes \va+\va\otimes(\vu+\vb)\right)\bigg)\\
&-& div\bigg(\eta\big((\vu+\vb) \otimes \va+\va \otimes (\vu+\vb)\big)\bigg)=\Delta \vUU+\vec{\mathcal{B}}+\vn \beta -div(\mathbb{B}),
\end{eqnarray*}
which is (\ref{EquationPourUtilisationMorrey}) as announced with the terms $\vec{\mathcal{B}}$, $\beta$ and $\mathbb{B}$ given in (\ref{BMorrey}), (\ref{BetaMorrey}) and (\ref{TensorBMorrey}), respectively.\\

%%%%%%%%%%%%%%%%%%%%%%%%%%%%%%%%%%%%%%
Once we have deduce the equation (\ref{EquationPourUtilisationMorrey}), in order to obtain the conclusion of the Theorem \ref{HolderRegularity_theorem} it is enough by Lemma \ref{parabolicHolder} to verify that we have
$$\vec{\mathcal{B}}\in \M^{\p_0, \q_0}  \qquad\mbox{and} \qquad \beta, \mathbb{B}\in \M^{\p_0, \q_1},$$
where $1< \p_0\le \q_0<\q_1$ with $\frac{1}{\q_0}=\frac{2-\alpha}{5}$, $\frac{1}{\q_1}=\frac{1-\alpha}{5}$ and $0<\alpha<\frac{1}{12}$. In the next proposition we will prove that under some extra hypothesis over the quantities $\vu, \vb, \vw$ (that will be proven in the next sections) the terms $\vec{\mathcal{B}}$, $\beta$ and $\mathbb{B}$ belong to the suitable Morrey spaces mentioned above.
%%%%%%%%%%%%%%%%%%%%%%%%%%%%%%%%%%%%%%
\begin{propo}\label{HolderRegularityproposition}
Let $R_1$, $R_2$ be positive numbers such that 
\begin{equation}\label{DefinitionRayons}
0<\mathbf{R}<R_2<R_1<4\mathbf{R},
\end{equation}
where $\mathbf{R}$ is fixed by the condition (\ref{Def_R}) above. Let $(\vu, \vb, \vw,  p)$ be a suitable solution for the equations MMP (\ref{EquationMMP}) over $\Omega$. Assume that we have the following points:
\begin{enumerate}
\item[1)] $\mathds{1}_{Q_{R_1}}\vu$, $\mathds{1}_{Q_{R_1}}\vb$, $\mathds{1}_{Q_{R_1}}\vw\in \M^{3,\tau_0}$ for $\frac{11}{2}>\tau_0>\frac{5}{1-\alpha}$ (recall that $0<\alpha<\frac{1}{12}$),\\
\item[2)] $\mathds{1}_{Q_{R_1}}\grad \otimes\vu$, $\mathds{1}_{Q_{R_1}}\grad \otimes\vb$, $\mathds{1}_{Q_{R_1}}\grad \otimes\vw \in \M^{2,\tau_1}$ with $\frac{1}{\tau_1}=\frac{1}{\tau_0}+\frac{1}{5}$,\\
\item[3)] $\mathds{1}_{Q_{R_2}}div(\vw)\in \M^{\frac{6}{5},\frac{12}{5}}$,
\item[4)] $\mathds{1}_{Q_{R_2}}\vu$, $\mathds{1}_{Q_{R_2}}\vb$, $\mathds{1}_{Q_{R_2}}\vw \in \M^{3,\delta}$ with $\delta\gg 1$ is such that $\frac{1}{\delta}+\frac{1}{\tau_0}\le \frac{1-\alpha}{5}$,\\
\item[5)] For all $1\le i,j \le 3$ we have 
$$\mathds{1}_{Q_{R_2}}\frac{\partial_i \partial_j}{(-\Delta)}(u_i b_j)\in \M^{\p,\q},\quad \mathds{1}_{Q_{R_2}}\frac{\partial_i \partial_j}{(-\Delta)}(u_i a_j)\in \M^{\p,\q}, \quad\mbox{ and }\quad  \mathds{1}_{Q_{R_2}}\frac{\partial_i \partial_j}{(-\Delta)}(b_i a_j)\in \M^{\p,\q},$$
with $\p_0\leq \p<+\infty$ and $\q_0<\q_1\leq \q<+\infty$.\\
\item[6)] $\mathds{1}_{Q_{R_1}}\vf \in \M^{\frac{10}{7}, \tau_a}, \mathds{1}_{Q_{R_1}}\vg \in \M^{\frac{10}{7}, \tau_b}$ for $\tau_a, \tau_b>\frac{5}{2-\alpha}$.
\end{enumerate}
If moreover $\va \in L^6_{t,x}(\Omega)$, then we have that the term $\vec{\mathcal{B}}$ defined in (\ref{BMorrey}) belongs to the Morrey space $\M^{\p_0,\q_0}$ with $1< \p_0\le \frac{6}{5}$ and $\frac{5}{2}<\q_0<3$ where $\frac{1}{\q_0}=\frac{2-\alpha}{5}$ with $0<\alpha<\frac{1}{12}$ and the terms $\beta, \mathbb{B}$ defined in (\ref{BetaMorrey}) and (\ref{TensorBMorrey}), respectively, belong to the Morrey space $\M^{\p_0, \q_1}$ with $\frac{1}{\q_1}=\frac{1-\alpha}{5}$.
\end{propo}
\begin{rema}\label{RemarqueIndicesTau1}
Note that, since $\q_0=\frac{5}{2-\alpha}$, $\tau_0>\frac{5}{1-\alpha}$ and $\frac{1}{\tau_1}=\frac{1}{\tau_0}+\frac{1}{5}$ then we easily obtain $\q_0<\tau_1<\tau_0$. Moreover, since $\q_1=\frac{5}{1-\alpha}$ we have $\q_1<\tau_0<\frac{11}{2}$. Remark also that since $0<\alpha<\frac{1}{12}$ we can set $\frac{11}{2}>\tau_0$ and $\frac{11}{2}>\tau_1$.
\end{rema}
Note that the conclusion of this proposition is exactly the input needed to apply Proposition \ref{parabolicHolder} from which we will obtain the wished gain of regularity. \\

%%%%%%%%%%%%%%%%%%%%%%%%%%%%%%%%%%%%%%
\noindent \textbf{Proof of the Proposition \ref{HolderRegularityproposition}}
In order to prove this proposition, and for the time being, let us take for granted the assumptions \emph{1) - 6)} above and let us prove that the quantities $\vec{\mathcal{B}}$, $\beta$ and $\mathbb{B}$ belong to the announced Morrey spaces. 
\begin{itemize}
\item {\bf For $\vec{\mathcal{B}}$.} We write, for $1<\p_0\le \frac{6}{5}$ and $\frac{5}{2}<\q_0<3$ where $\frac{1}{\q_0}=\frac{2-\alpha}{5}$ with $0<\alpha<\frac{1}{12}$:
\begin{eqnarray}
\|\vec{\mathcal{B}}\|_{\M^{\p_0,\q_0}}&\leq &\underbrace{\|(\partial_t\eta- \Delta \eta)(\vu+\vb+\vw)\|_{\M^{\p_0,\q_0}}}_{(1)}+2\sum_{i=1}^3\underbrace{\|(\partial_i \eta)(\partial_i (\vu+\vb+\vw))\|_{\M^{\p_0,\q_0}}}_{(2)}\notag\\
&+&\underbrace{\left\|\eta\bigg((\vb \cdot \grad)\vu+(\vu \cdot \grad)\vb\bigg)\right\|_{\M^{\p_0,\q_0}}}_{(3)}\label{EstimationMorreyTermeMathcalB}\\
&+&\underbrace{\left\|(\vn\eta) \bigg(\frac{div }{(-\Delta)}div\left(\vb\otimes \vu+\vu\otimes \vb+(\vu+\vb)\otimes \va+\va\otimes(\vu+\vb)\right)\bigg)\right\|_{\M^{\p_0,\q_0}}}_{(4)}\notag\\
&+&\underbrace{\|\eta(\rot\vw)\|_{\M^{\p_0,\q_0}}}_{(5)}+\underbrace{\|\eta(\vf+\vg)\|_{\M^{\p_0,\q_0}}}_{(6)}+\underbrace{\left\|\bigg((\vu+\vb) \otimes \va+\va \otimes (\vu+\vb)\bigg)\cdot \vn\eta\right\|_{\M^{\p_0,\q_0}}}_{(7)}\notag\\
&+&\underbrace{\|(\vn \eta) div(\vw)\|_{\M^{\p_0,\q_0}}}_{(8)}+\underbrace{\left\|\eta \left(\vw+\frac{1}{2}((\vu+\vb)\cdot \grad)\vw+\frac{1}{4}\rot(\vu+\vb)\right)\right\|_{\M^{\p_0,\q_0}}}_{(9)}.\notag
\end{eqnarray}
For the first term of (\ref{EstimationMorreyTermeMathcalB}), since we have $\mathds{1}_{Q_{R_1}}\vu$, $\mathds{1}_{Q_{R_1}}\vb$, $\mathds{1}_{Q_{R_1}}\vw\in \M^{3,\tau_0}$ for $\tau_0>5$ and since we have the support property $supp\;  (\partial_t \eta - \Delta \eta)\subset Q_{\mathbf{R}}$, it follows by Lemma \ref{lemma_locindi} (as we have $1\leq\p_0\leq \frac65<3$ and $\q_0<3<\tau_0$) that
\begin{equation*}
\| (\partial_t \eta - \Delta \eta)(\vu+\vb+\vw)\|_{\M^{\p_0,\q_0}}\le C\| \mathds{1}_{Q_{\mathbf{R}}}(\vu+\vb+\vw)\|_{\M^{\p_0,\q_0}}\le\| \mathds{1}_{Q_{R_1}}(\vu+\vb+\vw)\|_{\M^{3,\tau_0}}<+\infty,
\end{equation*}
where we used the information available in the point \emph{1)} of the Proposition \ref{HolderRegularityproposition}.\\

For the second term of (\ref{EstimationMorreyTermeMathcalB}), using Hölder's inequality in Morrey spaces (see the third point of Lemma \ref{lemma_Product}) we have $\| (\partial_i \eta)(\partial_i(\vu+\vb+\vw))\|_{\M^{\p_0,\q_0}}\le \| \mathds{1}_{Q_{\mathbf{R}}}\partial_i \eta\|_{\M^{\p_1,{q}_1}}\|\mathds{1}_{Q_{\mathbf{R}}}\partial_i(\vu+\vb+\vw)\|_{\M^{2,{q}_2}}$
where $\frac{1}{p_1}+\frac{1}{2}\le \frac{1}{\p_0}$ and $\frac{1}{{q}_1}+\frac{1}{q_2}=\frac{1}{\q_0}$. Since $\frac{5}{2}<\q_0=\frac{5}{2-\alpha}<3$, we have that $q_2$ can be chosen such that $q_2<\tau_1=\frac{5\tau_0}{5+\tau_0}$ and thus using Lemma \ref{lemma_locindi} (recall that $\mathbf{R}<R_1$)  and the point \emph{2)} of Proposition \ref{HolderRegularityproposition}, we obtain 
\begin{equation*}
\| (\partial_i \eta)(\partial_i(\vu+\vb+\vw)\|_{\M^{\p_0,\q_0}}\le C\|\mathds{1}_{Q_{R_1}}\grad \otimes(\vu+\vb+\vw)\|_{\M^{2,\tau_1}}<+\infty.
\end{equation*}
For the term (3) of (\ref{EstimationMorreyTermeMathcalB}), since $1< \mathfrak{p}_0\leq \tfrac{6}{5}$ and $\tfrac{5}{2}<\mathfrak{q}_0<3$, by Lemma \ref{lemma_locindi} (recall that $\mathbf{R}<R_2<R_1$), by the H\"older inequalities in Morrey spaces and using the information of points \emph{2)-4)}, we have:
\begin{eqnarray}
\left\|\eta \left( (\vb\cdot\vn)\vu+(\vu\cdot\vn)\vb\right)\right\|_{\mathcal{M}_{t,x}^{\p_0, \mathfrak{q}_0}} \leq  C\left\|\mathds{1}_{Q_{\mathbf{R}}} \left( (\vb\cdot\vn)\vu+(\vu\cdot\vn)\vb\right)\right\|_{\mathcal{M}_{t,x}^{\frac{6}{5}, \mathfrak{q}_0}}\qquad \qquad\label{EstimationUVBMorreyLoc1}\\
\leq C\left( \|\mathds{1}_{Q_{R_2}}\vb\|_{\mathcal{M}_{t,x}^{3, \delta}}\|\mathds{1}_{Q_{R_1}}\vn\otimes\vu\|_{\mathcal{M}_{t,x}^{2, \tau_{1}}}+\|\mathds{1}_{Q_{R_2}}\vu\|_{\mathcal{M}_{t,x}^{3, \delta}}\|\mathds{1}_{Q_{R_1}}\vn\otimes\vb\|_{\mathcal{M}_{t,x}^{2, \tau_{1}}}\right)<+\infty\notag,
\end{eqnarray}
where we have $\frac{1}{\delta}+\frac{1}{\tau_{1}}\leq \frac{1}{\mathfrak{q}_0}$, but since $\frac{1}{\mathfrak{q}_0}=\frac{2-\alpha}{5}$ and $\frac{1}{\tau_{1}}=\frac{1}{\tau_{0}}+\frac{1}{5}$, the previous conditions is equivalent to $\frac{1}{\delta}+\frac{1}{\tau_{0}}\leq \frac{1-\alpha}{5}$, which is exactly the condition stated in the point \emph{4)} of the Proposition \ref{HolderRegularityproposition}.\\

For the term (4) of (\ref{EstimationMorreyTermeMathcalB}), due to the symmetry of the information available in the point \emph{5)} of the Proposition \ref{HolderRegularityproposition}, it is enough to study the following term for $1\leq i,j\leq 3$ and due to the support properties of the function $\eta$, we obtain
$$\left\|(\vn\eta) \frac{\partial_i\partial_j }{(-\Delta)}a_ib_j\right\|_{\M^{\p_0,\q_0}}\leq C\left\|\mathds{1}_{Q_{R_2}} \frac{\partial_i\partial_j }{(-\Delta)}a_ib_j\right\|_{\M^{\p_0,\q_0}}\leq C\left\|\mathds{1}_{Q_{R_2}} \frac{\partial_i\partial_j }{(-\Delta)}a_ib_j\right\|_{\M^{\p,\q}}<+\infty,$$
where we applied Lemma \ref{lemma_locindi} with $\p_0\leq \p$ and $\q_0\leq\q$.\\

For the terms (5), (8) and (6) of (\ref{EstimationMorreyTermeMathcalB}) can be treated in the same manner, indeed, by the assumption \emph{2)} of Proposition \ref{HolderRegularityproposition} we have
$$\|\eta(\rot\vw)\|_{\M^{\p_0,\q_0}}\leq C\|\mathds{1}_{Q_{R_1}}\vn\otimes\vw\|_{\M^{2,\tau_1}}<+\infty, \qquad \|(\vn \eta) div(\vw)\|_{\M^{\p_0,\q_0}}\leq C\|\mathds{1}_{Q_{R_1}}\vn\otimes\vw\|_{\M^{2,\tau_1}}<+\infty,$$
where we used Lemma \ref{lemma_locindi} with $\p_0\leq\frac65<2$ and $\q_0<\tau_1$ (see Remark \ref{RemarqueIndicesTau1}). By essentially the same arguments, using the point \emph{6)} of Proposition \ref{HolderRegularityproposition} (and since $\p_0\leq\frac65<\frac{10}{7}$ and $\tau_a, \tau_b>\frac{5}{2-\alpha}=\q_0$) we have $\|\eta(\vf+\vg)\|_{\M^{\p_0,\q_0}}\leq C(\|\mathds{1}_{Q_{R_1}}\vf\|_{\M^{\frac{10}{7}, \tau_a}}+\|\mathds{1}_{Q_{R_1}}\vg\|_{\M^{\frac{10}{7}, \tau_b}})<+\infty$.\\

For the term (7) of (\ref{EstimationMorreyTermeMathcalB}), as we have the same information over $\vu$ and $\vb$ we only need to study (for $1\leq i,j,k\leq 3$) the terms of the form $\|u_ia_j \partial_k\eta\|_{\M^{\p_0,\q_0}}$ and we have
\begin{eqnarray*}
\left\|u_ia_j \partial_k\eta\right\|_{\M^{\p_0,\q_0}}&\leq &C\left\|\mathds{1}_{Q_{\mathbf{R}}}u_ia_j \right\|_{\M^{\p_0,\q_0}}\leq C\|\mathds{1}_{Q_{R_2}}u_i\|_{\M^{3,\delta}}\|\mathds{1}_{Q_{R_1}}a_j\|_{\M^{2,\tau_1}}\\
&\leq & C\|\mathds{1}_{Q_{R_2}}u_i\|_{\M^{3,\delta}}\|\mathds{1}_{Q_{R_1}}a_j\|_{\M^{6,6}}<+\infty,
\end{eqnarray*}
where we used the H\"older inequality in Morrey spaces, Lemma \ref{lemma_locindi} (with $2<6$ and $\tau_1<6$ by Remark \ref{RemarqueIndicesTau1}), the point \emph{4)} of Proposition \ref{HolderRegularityproposition} and the fact that $\M^{6,6}=L^6_{t,x}$.\\

For the term (9) of (\ref{EstimationMorreyTermeMathcalB}) we easily deduce 
$\|\eta \vw\|_{\M^{\p_0,\q_0}}\leq C\|\mathds{1}_{Q_{R_1}}\vw\|_{\M^{3,\tau_0}}<+\infty$ (by Lemma \ref{lemma_locindi}  since $\p_0<3$ and $\q_0<\tau_0$). Due to the symmetry of the information available for the terms $\vu, \vb$ and $\vw$ and  following the same ideas displayed in (\ref{EstimationUVBMorreyLoc1}), we have $\|\eta ((\vu+\vb)\cdot \grad)\vw\|_{\M^{\p_0,\q_0}}<+\infty$. Finally, since by the point \emph{3)} of Proposition \ref{HolderRegularityproposition} we have $\mathds{1}_{Q_{R_1}}\rot\vu$ and $\mathds{1}_{Q_{R_1}}\rot\vb\in \M^{2,\tau_1}$, and since $\p_0<2$ and $\q_0<\tau_1$, by Lemma \ref{lemma_locindi} we obtain
$\|\eta\rot(\vu+\vb)\|_{\M^{\p_0,\q_0}}\leq C(\|\mathds{1}_{Q_{R_1}}\rot\vu\|_{\M^{2,\tau_1}}+\|\mathds{1}_{Q_{R_1}}\rot\vb\|_{\M^{2,\tau_1}})<+\infty$. We thus have:
$$\left\|\eta \left(\vw+\frac{1}{2}((\vu+\vb)\cdot \grad)\vw+\frac{1}{4}\rot(\vu+\vb)\right)\right\|_{\M^{\p_0,\q_0}}<+\infty.$$

\item {\bf For $\beta$.} By the expression (\ref{BetaMorrey}), we have, for $1<\p_0\le \frac{6}{5}$ and  $\q_1=\frac{5}{1-\alpha}$ with $0<\alpha<\frac{1}{12}$, 
\begin{equation}\label{EstimationMorreyTermeBeta}
\|\beta\|_{\M^{\p_0, \q_1}}\leq \|\eta div(\vw)\|_{\M^{\p_0, \q_1}}+\left\|\eta \frac{div div}{(-\Delta)}\left(\vb\otimes \vu+\vu\otimes \vb+(\vu+\vb)\otimes \va+\va\otimes(\vu+\vb)\right)\right\|_{\M^{\p_0, \q_1}}.
\end{equation}
Since by the point \emph{3)} of Proposition \ref{HolderRegularityproposition} we have $\mathds{1}_{Q_{R_2}}div(\vw)\in \M^{\frac{6}{5},\frac{12}{5}}$ and since $\p_0\leq \frac65$ and $\q_1<\frac{12}{5}$ (see Remark \ref{RemarqueIndicesTau1}), then, by Lemma \ref{lemma_locindi} we have for the first term fo the right-hand side above:
$$\|\eta div(\vw)\|_{\M^{\p_0, \q_1}}\leq C\|\mathds{1}_{Q_{R_2}} div(\vw)\|_{\M^{\p_0, \q_1}}\leq C\|\mathds{1}_{Q_{R_2}} div(\vw)\|_{\M^{\frac{6}{5},\frac{12}{5}}}<+\infty.$$
For the second term of the right-hand side of (\ref{EstimationMorreyTermeBeta}), we use the point \emph{5)} of Proposition \ref{HolderRegularityproposition} and due to the symmetry of the information available, it is enough to study, for $1\leq i,j\leq 3$ the term $\|\eta \frac{\partial_i \partial_j}{(-\Delta)}(u_i b_j)\|_{\M^{\p_0, \q_1}}$, and we write
$$\left\|\eta \frac{\partial_i \partial_j}{(-\Delta)}(u_i b_j) \right\|_{\M^{\p_0, \q_1}}\leq C\left\|\mathds{1}_{Q_{R_2}} \frac{\partial_i \partial_j}{(-\Delta)}(u_i b_j) \right\|_{\M^{\p_0, \q_1}}\leq C\left\|\mathds{1}_{Q_{R_2}} \frac{\partial_i \partial_j}{(-\Delta)}(u_i b_j) \right\|_{\M^{\p, \q}}<+\infty,$$
where we applied Lemma \ref{lemma_locindi} with $\p_0<\p$ and $\q_1<\q$.\\

\item {\bf For $\mathbb{B}$.} By (\ref{TensorBMorrey}) we need to study the quantity $\|\mathbb{B}\|_{\M^{\p_0, \q_1}}=\|\eta( \va\otimes (\vu+\vb)+(\vu+\vb)\otimes\va)\|_{\M^{\p_0, \q_1}}$, for the sake of simplicity we only study $\|\eta \va\otimes \vu\|_{\M^{\p_0, \q_1}}$ as the other terms can be treated in the same manner. We thus have
$$\|\eta \va\otimes \vu\|_{\M^{\p_0, \q_1}}\leq C\|\mathds{1}_{Q_{R_2}} \va\otimes \vu\|_{\M^{\frac65, \q_1}}\leq C\|\mathds{1}_{Q_{R_2}} \va\|_{\M^{2, \delta'}}\|\mathds{1}_{Q_{R_2}} \vu\|_{\M^{3,\delta}},$$
where we used the H\"older inequalities in Morrey spaces with $\frac{1}{\q_1}=\frac{1}{\delta}+\frac{1}{\delta'}$. Since by the point \emph{4)} of Proposition \ref{HolderRegularityproposition} the index $\delta\gg1$ can be chosen big enough such that $\delta'<6$, thus we have by Lemma \ref{lemma_locindi}:
$$\|\eta \va\otimes \vu\|_{\M^{\p_0, \q_1}}\leq C\|\mathds{1}_{Q_{R_2}} \va\|_{\M^{6, 6}}\|\mathds{1}_{Q_{R_2}} \vu\|_{\M^{3,\delta}}<+\infty,$$
since $\M^{6, 6}=L^6_{t,x}$.
\end{itemize}
We have proven that $\vec{\mathcal{B}}\in \M^{\p_0, \q_0}$ and $\beta, \mathbb{B}\in \M^{\p_0, \q_1}$ where $1< \p_0\le \q_0<\q_1$ with $\frac{1}{\q_0}=\frac{2-\alpha}{5}$, $\frac{1}{\q_1}=\frac{1-\alpha}{5}$ and $0<\alpha<\frac{1}{12}$, and thus the proof of Proposition \ref{HolderRegularityproposition} is finished.\hfill$\blacksquare$
%%%%%%%%%%%%%%%%%%%%%%%%%%%%%%%%%%%%%%%%%%%%%%%%%%%%%%%%%%%%%%%%%%%%%%%%%%%%
\section{Local Energy Inequality and Useful estimates}\label{Secc_LocalEnergy}
In order to obtain some of the assumptions stated in Proposition \ref{HolderRegularityproposition}, we will exploit the information given by the local energy estimate that can be deduced from the structure of the equation (\ref{EquationMMP}). We know from the work of Scheffer \cite{Scheffer1, Scheffer} that the use of a special test function is particularly helpful to obtain good estimates. We will use the following function:
\begin{lem}\label{testfonc}
Let $0<\rho\le 1$ and $0<r<\frac{\rho}{2}$. Let $\phi \in \mathcal{C}_0^{\infty}
(\mathbb{R}\times \R)$  be defined by the formula 
\begin{equation}\label{Def_TestFuncionSchaefer1}
\phi(s,y)=r^2\omega\left(\frac{s-t}{\rho^2},\frac{y-x}{\rho}\right)\theta\left(\frac{s-t}{r^2}\right)\mathfrak{g}_{(4r^2+t-s)}(x-y), \quad 0<r<\frac{\rho}{2}\le 1,
\end{equation}
where $\omega \in \mathcal{C}_0^{\infty}(\mathbb{R}\times \R)$ is non-negative 
function supported on the parabolic ball $Q_1(0,0)$ and is equal to 1 on $Q_{\frac{1}{2}}(0,0)$ (see formula (\ref{Def_BoulesQ})), $\theta$ is a non-negative smooth function such that $\theta=1$ on $]-\infty,1[$ and $\theta=0$ on $]2,+\infty[$ and $\mathfrak{g}_t(\cdot)$ is the usual heat kernel. Then, we have the following points
\begin{itemize}
\item[1)]the function $\phi$ is a bounded non-negative function, and its support is contained in the parabolic ball $Q_\rho$, and for all $(s,y)\in Q_r(t,x)$ we have the lower bound $\phi\ge \frac{C}{r}$,
\item[2)] for all $(s,y)\in Q_\rho (t,x)$ with $0<s<t+r^2$ we have $\phi(s,y)\le \frac{C}{r}$,
\item[3)] for all $(s,y)\in Q_\rho(t,x) $ with $0<s<t+r^2$ we have $\grad \phi(s,y)\le \frac{C}{r^2}$,
\item [4)] moreover, for all $(s,y)\in Q_\rho(t,x) $ with $0<s<t+r^2$ we have $|(\partial_s+\Delta)\phi(s,y)|\le C\frac{r^2}{\rho^5}$.
\end{itemize}
\end{lem}
See the book \cite[Section 13.9]{PGLR1} for a proof of this lemma. See also the Appendix B of \cite{MHDKukavica}.\\

Now, with the help of this function we have the local energy inequality:
%%%%%%%%%%%%%%%%%%%%%%%%%%%%%%%%%%%%%%
\begin{propo}\label{Propo_InegaliteEnergieLocale}
Let $(\vu, \vb, \vw, p)$ be a weak solution of the MMP equation (\ref{EquationMMP}) over a subset $\Omega$ of the form (\ref{DefConjuntoOmega}) and assume that $\phi$ is the function given in (\ref{Def_TestFuncionSchaefer1}). Then the local energy inequality for the MMP equation is given by
\begin{eqnarray}
&& \int_{\R}[( |\vu|^2+|\vb|^2+|\vw|^2)\phi](\tau,x)dx+2\int_{s<\tau}\int _{\R}[(|\grad \otimes \vu|^2+|\grad \otimes \vb|^2+|\grad \otimes \vw|^2 )\phi] (s,x)dxds\notag\\
&&+2\int_{s<\tau}\int_{\R}[(div(\vw))^2\phi](s,x)dxds+2\int_{s<\tau}\int_{\R}[|\vw|^2\phi](s,x)dxds\notag\\
& \leq&\int_{s< \tau}\int_{\R}[(\partial_t\phi +\Delta \phi)( |\vu|^2+|\vb|^2+|\vw|^2)](s,x)dxds+\int_{s< \tau}\int_{\R}[(|\vu|^2 +2p)\vb\cdot \vn\phi](s,x)dxds\label{EstimationEnergieLocale}\\
&&+\int_{s< \tau}\int_{\R}[(|\vb|^2+2p)\vu\cdot \vn\phi](s,x)dxds+\int_{s< \tau}\int_{\R}(\rot\vw) \cdot[\phi(\vu+\vb)](s,x) dxds\notag\\
&&+2\int_{s< \tau}\int_{\R}[\vf \cdot (\phi \vu)+\vg \cdot (\phi \vb)](s,x)dxds+2\int_{s<\tau}\int_{\R}\left|[div(\vw)(\vn\phi\cdot \vw)](s,x)\right|dxds\notag\\
&&+\int_{s< \tau}\int_{\R}\left|\bigg[[(\va\cdot \vn)(\vu+\vb)]\cdot (\phi (\vu+\vb))\bigg](s,x)\right|dxds+\int_{s< \tau}\int_{\R}\bigg[\big((\vu+\vb)\cdot \vn\big)\big(\phi(\vu+\vb)\big)\cdot\va\bigg] (s,x)dxds\notag\\
&&+\frac{1}{2}\int_{s< \tau}\int_{\R}[|\vw|^2(\vu+\vb)\cdot \grad\phi](s,x) dxds
+\frac{1}{2}\int_{s< \tau}\int_{\R}[\rot(\vu+\vb)]\cdot(\phi \vw)(s,x)dxds.\notag
\end{eqnarray}
\end{propo}
%%%%%%%%%%%%%%%%%%%%%%%%%%%%%%%%%%%%%%
\textbf{Proof.} In order to deduce the local energy inequality announced, we multiply the three first equations of the system (\ref{EquationMMP}) by
$\phi\vu$, $\phi\vb$ and $\phi\vw$ respectively and we integrate in the space variable to obtain
\begin{eqnarray*}
\int_{\R}\partial_t \vu\cdot(\phi \vu)dx&=&\int_{\R}\left( \Delta \vu  -(\vb \cdot \grad)\vu-\grad p +\frac{1}{2}\rot\vw+ \vf-\frac{1}{2}div((\vu+\vb) \otimes \va+\va \otimes (\vu+\vb))\right)\cdot(\phi \vu)dx,\\
\int_{\R}\partial_t \vb\cdot(\phi \vb)dx&=&\int_{\R}\left(\Delta \vb -(\vu \cdot \grad)\vb -\grad p +\frac{1}{2}\rot\vw +\vg-\frac{1}{2}div((\vu+\vb) \otimes \va +\va \otimes (\vu+\vb))\right)\cdot(\phi \vb)dx,\\
\int_{\R}\partial_t \vw\cdot(\phi \vw)dx&=&\int_{\R}\left( \Delta \vw +\grad div(\vw)-\vw-\frac{1}{2}((\vu+\vb)\cdot \grad)\vw+\frac{1}{4}\rot(\vu+\vb)\right)\cdot(\phi \vw)dx.
\end{eqnarray*}
Recalling that we have the generic identity $\partial_t \vec c\cdot (\phi \vec c)=\frac{1}{2}\partial_t(|\vec c|^2 \phi)-\frac{1}{2}|\vec c|^2 \partial _t \phi$ as well as the formulas $\displaystyle{\int _{\R}\Delta \vec c \cdot (\phi \vec c)dx=\frac{1}{2}\int _{\R}|\vec c|^2\Delta \phi dx-\int _{\R}|\grad \otimes \vec c|^2\phi dx}$ and $\displaystyle{\int _{\R} [(\vec c\cdot \grad)\vec d]\cdot (\phi \vec d)dx=-\frac{1}{2} \int _{\R} |\vec d|^2\vec c \cdot \grad\phi dx}$  which are valid for any (smooth) divergence free vector field $\vec c$, we obtain after some integration by parts and after an integration in the time variable:
\begin{eqnarray*}
&& \int_{\R}[( |\vu|^2+|\vb|^2+|\vw|^2)\phi](\tau,x)dx+2\int_{s<\tau}\int _{\R}[(|\grad \otimes \vu|^2+|\grad \otimes \vb|^2+|\grad \otimes \vw|^2 )\phi] (s,x)dxds\\
&&+2\int_{s<\tau}\int_{\R}[(div(\vw))^2\phi](s,x)dxds+2\int_{s<\tau}\int_{\R}[|\vw|^2\phi](s,x)dxds\\
& \leq&\int_{s< \tau}\int_{\R}[(\partial_t\phi +\Delta \phi)( |\vu|^2+|\vb|^2+|\vw|^2)](s,x)dxds+\int_{s< \tau}\int_{\R}[(|\vu|^2 +2p)\vb\cdot \vn\phi](s,x)dxds\\
&&+\int_{s< \tau}\int_{\R}[(|\vb|^2+2p)\vu\cdot \vn\phi](s,x)dxds+\int_{s< \tau}\int_{\R}(\rot\vw) \cdot[\phi(\vu+\vb)](s,x) dxds\\
&&+2\int_{s< \tau}\int_{\R}[\vf \cdot (\phi \vu)+\vg \cdot (\phi \vb)](s,x)dxds-2\int_{s<\tau}\int_{\R}[div(\vw)(\vn\phi\cdot \vw)](s,x) dxds\\
&&-\int_{s< \tau}\int_{\R}div\left((\vu+\vb) \otimes \va+\va \otimes (\vu+\vb)\right)\cdot[\phi (\vu+\vb)](s,x)dxds\\
&&+\frac{1}{2}\int_{s< \tau}\int_{\R}[|\vw|^2(\vu+\vb)\cdot \grad\phi](s,x) dxds
+\frac{1}{2}\int_{s< \tau}\int_{\R}[\rot(\vu+\vb)]\cdot(\phi \vw)(s,x)dxds,
\end{eqnarray*}
since $\vu, \vb$ and $\va$ are divergence free vector fields, we easily see that the quantity $\displaystyle{\int_{\R}div\left((\vu+\vb) \otimes \va+\va \otimes (\vu+\vb)\right)\cdot[\phi (\vu+\vb)]dx}$ can be rewritten as 
\begin{eqnarray*}
&&\int_{\R}[(\va\cdot \vn)\vu]\cdot (\phi \vu)+[(\va\cdot \vn)\vb]\cdot (\phi \vu)+[(\va\cdot \vn)\vu]\cdot (\phi \vb)+[(\va\cdot \vn)\vb]\cdot (\phi \vb)\\[2mm]
&&+[(\vu\cdot \vn)\va]\cdot (\phi \vu)+[(\vb\cdot \vn)\va]\cdot (\phi \vu)+[(\vu\cdot \vn)\va]\cdot (\phi \vb)+[(\vb\cdot \vn)\va]\cdot (\phi \vb)\; dx,
\end{eqnarray*}
for the last line above we will use the identity $\displaystyle{\int_{\R}[(\vec c\cdot \grad)\vec d]\cdot (\phi \vec e )dx=-\int_{\R}[(\vec c\cdot \grad)(\phi\vec e)]\cdot  \vec d\;dx}$ which is valid for divergence free vector fields, and using the bilinear structure of the terms, we have
\begin{eqnarray*}
\int_{\R}div\left((\vu+\vb) \otimes \va+\va \otimes (\vu+\vb)\right)\cdot[\phi (\vu+\vb)]dx&=&\int_{\R}\bigg[[(\va\cdot \vn)(\vu+\vb)]\cdot (\phi (\vu+\vb))\bigg]dx\\
&&-\int_{\R}\bigg[\big((\vu+\vb)\cdot \vn\big)\big(\phi(\vu+\vb)\big)\cdot\va\bigg]dx,
\end{eqnarray*}
and we finally obtain
\begin{eqnarray*}
&& \int_{\R}[( |\vu|^2+|\vb|^2+|\vw|^2)\phi](\tau,x)dx+2\int_{s<\tau}\int _{\R}[(|\grad \otimes \vu|^2+|\grad \otimes \vb|^2+|\grad \otimes \vw|^2 )\phi] (s,x)dxds\\
&&+2\int_{s<\tau}\int_{\R}[(div(\vw))^2\phi](s,x)dxds+2\int_{s<\tau}\int_{\R}[|\vw|^2\phi](s,x)dxds\\
& \leq&\int_{s< \tau}\int_{\R}[(\partial_t\phi +\Delta \phi)( |\vu|^2+|\vb|^2+|\vw|^2)](s,x)dxds+\int_{s< \tau}\int_{\R}[(|\vu|^2 +2p)\vb\cdot \vn\phi](s,x)dxds\\
&&+\int_{s< \tau}\int_{\R}[(|\vb|^2+2p)\vu\cdot \vn\phi](s,x)dxds+\int_{s< \tau}\int_{\R}(\rot\vw) \cdot[\phi(\vu+\vb)](s,x) dxds\\
&&+2\int_{s< \tau}\int_{\R}[\vf \cdot (\phi \vu)+\vg \cdot (\phi \vb)](s,x)dxds+2\int_{s<\tau}\int_{\R}\left|[div(\vw)(\vn\phi\cdot \vw)](s,x)\right|dxds\\
&&+\int_{s< \tau}\int_{\R}\left|\bigg[[(\va\cdot \vn)(\vu+\vb)]\cdot (\phi (\vu+\vb))\bigg](s,x)\right|dxds+\int_{s< \tau}\int_{\R}\bigg[\big((\vu+\vb)\cdot \vn\big)\big(\phi(\vu+\vb)\big)\cdot\va\bigg] (s,x)dxds\\
&&+\frac{1}{2}\int_{s< \tau}\int_{\R}[|\vw|^2(\vu+\vb)\cdot \grad\phi](s,x) dxds
+\frac{1}{2}\int_{s< \tau}\int_{\R}[\rot(\vu+\vb)]\cdot(\phi \vw)(s,x)dxds,
\end{eqnarray*}
and this ends the proof of Proposition \ref{Propo_InegaliteEnergieLocale}.\hfill$\blacksquare$\\

%%%%%%%%%%%%%%%%%%%%%%%%%%%%%%%%%%%%%%
Once we have obtained this inequality, we will make use of the properties of the test function $\phi$ given in Lemma \ref{testfonc} in order to obtain suitable controls that will be used in the next section. Indeed, by introducing some scaled quantities it would be possible to exploit the previous inequality (\ref{EstimationEnergieLocale}) to deduce by an inductive argument some stability of this scaled quantities in terms of Morrey spaces.\\

In this sense we have the following definition. 
%%%%%%%%%%%%%%%%%%%%%%%%%%%%%%%%%%%%%%
\begin{defi}[Scaled Quantities]\label{Def_Invariants}
Let $ \gamma>0$. For all $(t,x)\in \mathbb{R}\times \R$, we consider the 
following scaled functions: 
\begin{eqnarray*}
&&\vu_ \gamma(t,x)= \gamma\vu( \gamma^2 t, \gamma x),\quad\vb_ \gamma(t,x)= \gamma\vb( \gamma^2 t, \gamma x),\quad\vw(t,x)
= \gamma \vw( \gamma^2 t, \gamma x)\\[2mm]
&&p_ \gamma(t,x)= \gamma^2 p( \gamma^2 t, \gamma x),\quad \vf_ \gamma(t,x) =  \gamma^3 \vf( \gamma^2 t,  \gamma x)\quad\mbox{and}\quad \vg_ \gamma(t,x) =  \gamma^3 \vg( \gamma^2 t,  \gamma x).
\end{eqnarray*}
Now we define the following invariant quantities with respect of the previous scaling:
\begin{align*}
\mathcal{A}_r(t,x)&=\sup_{t-r^2<s<t+r^2}\frac{1}{r}\int_{B(x,r)}|\vu(s,y)|^2dy 
&&\alpha_r(t,x)=\frac{1}{r} \iint _{Q_r(t,x)}|\grad \otimes \vu(s,y)|^2dyds\\
\mathcal{B}_r(t,x)&=\sup_{t-r^2<s<t+r^2}\frac{1}{r}\int_{B(x,r)}|\vb(s,y)|^2dy 
&&\beta_r(t,x)=\frac{1}{r}\iint _{Q_r(t,x)}|\grad \otimes \vb(s,y)|^2dyds\\
\mathcal{C}_r(t,x)&=\sup_{t-r^2<s<t+r^2}\frac{1}{r}\int_{B(x,r)}|\vw(s,y)|^2dy 
&& \gamma_r(t,x)=\frac{1}{r}\iint _{Q_r(t,x)}|\grad \otimes \vw(s,y)|^2dyds\\
\lambda_r(t,x)&=\frac{1}{r^2}\iint _{Q_r(t,x)}|\vu(s,y)|^3dyds 
&&\zeta_r(t,x)=\frac{1}{r^2}\iint _{Q_r(t,x)}|\vb(s,y)|^3dyds\\
\sigma_r(t,x)&=\frac{1}{r^2}\iint _{Q_r(t,x)}|\vw(s,y)|^3dyds\\
\mathcal{W}_r(t,x)&=\frac{1}{r}\iint _{Q_r(t,x)}|div( \vw)(s,y)|^2dyds&&
\mathcal{H}_r(t,x)=\frac{1}{r^3}\iint _{Q_r(t,x)}|\vw(s,y)|^2dyds\\
\mathcal{F}_r(t,x)&=\frac{1}{r^{\frac{5}{7}}}\iint _{Q_r(t,x)}|\vf(s,y)|^{\frac{10}{7}}dyds
&&\mathcal{G}_r(t,x)=\frac{1}{r^{\frac{5}{7}}}\iint _{Q_r(t,x)}|\vg(s,y)|^{\frac{10}{7}}dyds\\
\mathcal{P}_r(t,x)&=\frac{1}{r^2}\iint _{Q_r(t,x)}|p (s,y)|^{\frac{3}{2}}dyds.\\
\end{align*}
\end{defi}
\begin{rema} \label{Remarque_DefinitionScaled}
From the definition above we easily deduce the identities
$(r\mathcal{A}_r)^{\frac12}=\|\vu\|_{L^\infty_tL^2_x(Q_r)}$, $(r\mathcal{\alpha}_r)^{\frac12}=\|\vu\|_{L^2_t\dot{H}^1_x(Q_r)}$ and $r^{\frac{4}{3}}\mathcal{P}_r^{\frac{2}{3}}=\|p\|_{L_{t,x}^{\frac{3}{2}}(Q_r)}$ and similar identities for the variables $\vb$ and $\vw$.
\end{rema}
As announced, we will use these quantities to deduce two main estimates, which are stated in Proposition \ref{Propo_FirstEstimate} and Proposition \ref{Propo_SecondEstimate}. In the next lemma we prove some useful relationships between some of the previous terms given above. 
%%%%%%%%%%%%%%%%%%%%%%%%%%%%%%%%%%%%%%
\begin{lem}\label{Lem_norml3}
For any small $0<r<1$ such that $Q_r\subset \Omega$ and under the general hypotheses stated in (\ref{HypothesesLocal1}), there exists an absolut constant $C$, such that
\begin{equation*}
\lambda_r^{\frac{1}{3}}\le C(\mathcal{A}_r+\alpha_r)^{\frac{1}{2}},\quad\quad \zeta_r^{\frac{1}{3}}\le C(\mathcal{B}_r+\beta_r)^{\frac{1}{2}}\quad\text{and}\quad \sigma_r^{\frac{1}{3}}\le C(\mathcal{C}_r+ \gamma_r)^{\frac{1}{2}}.
\end{equation*}
\end{lem}
\textbf{Proof.} We only detail the proof of the first estimate as the two others follow the same arguments. Thus, by the expression of $\lambda_r$ given in the Definition \ref{Def_Invariants} and Hölder's inequality, we have the estimate $\lambda_r^{\frac 13} = \frac{1}{r^{\frac 23}} \|\vu\|_{L_{t,x}^{3} (Q_r)}\leq C \frac{1}{r^{\frac 12}} \|\vu\|_{L_{t,x}^{\frac{10}{3}} (Q_r)}$. Now, 
using an interpolation inequality we have the control $\|\vu(t,\cdot)\|_{L^{\frac{10}{3}}(B_r)}\leq \|\vu(t,\cdot)\|_{L^{2}(B_r)}^{\frac{2}{5}}\|\vu(t,\cdot)\|_{L^{6}(B_r)}^{\frac{3}{5}}$ and applying the Hölder inequality with respect to the time variable, we obtain $\|\vu\|_{L_{t,x}^{\frac{10}{3}} (Q_r)} \leq \|\vu\|_{L_t^{\infty}L_x^{2} (Q_r)}^{\frac25}\|\vu\|_{L_t^{2}L_x^{6} (Q_r)}^{\frac{3}{5}}$. For the $L_t^2L_x^6$ norm of $\vu$, we use the classical Gagliardo-Nirenberg inequality (see \cite{Brezis}) to obtain $\|\vu\|_{L_t^{2}L_x^{6} (Q_r)} \leq C\left(\|\vn \otimes\vu\|_{L_t^{2}L_x^{2} (Q_r)} +\|\vu\|_{L_t^{\infty}L_x^{2} (Q_r)}\right)$ and using Young's inequalities we have
\begin{eqnarray*}
\|\vu\|_{L_{t,x}^{\frac{10}{3}} (Q_r)} &\leq& C \|\vu\|_{L_t^{\infty}L_x^{2} (Q_r)}^{\frac 25}\left(\|\vn \otimes \vu\|_{L_t^{2}L_x^{2} (Q_r)}^{\frac 35}+\|\vu\|_{L_t^{\infty}L_x^{2} (Q_r)}^{\frac35} \right) \leq C\left(\|\vu\|_{L_t^{\infty}L_x^{2} (Q_r)}+\|\vn\otimes \vu\|_{L_t^{2}L_x^{2} (Q_r)}\right).
\end{eqnarray*}
Noting that $\|\vu\|_{L_t^{\infty}L_x^{2} (Q_r)}=r^{\frac12}\mathcal{A}_r^{\frac12}$ and $\|\vn\otimes \vu\|_{L^2_tL_x^{2} (Q_r)}=r^{\frac12}\alpha_r^{\frac12}$, we finally obtain $\lambda_r^{\frac 13} \leq  C  (\mathcal{A}_r + \alpha_r)^{\frac12}$ and Lemma \ref{Lem_norml3} is proven. \hfill$\blacksquare$
%%%%%%%%%%%%%%%%%%%%%%%%%%%%%%%%%%%%%%
\subsection{A first estimate}
We give now the first general inequality that  bounds all the terms given in the Definition \ref{Def_Invariants}.
\begin{propo}[First Estimate]\label{Propo_FirstEstimate}
Under the hypotheses of Theorem \ref{HolderRegularity_theorem}, for $0<r<\frac{\rho}{2}<1$, we have 
\begin{align}
\mathcal{A}_r+\mathcal{B}_r+\mathcal{C}_r&+\alpha_r+\beta_r+ \gamma_r +\mathcal{W}_r+ r^2\mathcal{H}_r \le C\frac{r^2}{\rho ^2}( \mathcal{A}_\rho+\mathcal{B}_\rho+ \mathcal{C}_\rho)+C\frac{\rho^2}{r^2} \alpha_\rho^{\frac{1}{2}}(\mathcal{A}_\rho+\mathcal{B}_\rho+\beta_\rho)\notag\\
&+ C\frac{\rho ^2}{r^2}\mathcal{P}_\rho^{\frac{2}{3}}( \mathcal{B}_\rho+\beta_\rho)^{\frac{1}{2}}+C\frac{\rho^2}{r^2} \beta_\rho^{\frac{1}{2}}(\mathcal{B}_\rho+\mathcal{A}_\rho+\alpha_\rho)+ C\frac{\rho ^2}{r^2}\mathcal{P}_\rho^{\frac{2}{3}}(\mathcal{A}_\rho+\alpha_\rho)^{\frac{1}{2}}+ C\frac{\rho^{2}}{r}\gamma_{\rho}^{\frac{1}{2}}(\mathcal{A}_{\rho}^{\frac{1}{2}}+\mathcal{B}_{\rho}^{\frac{1}{2}})\notag\\
&+C\frac{\rho}{r}\left(\mathcal{F}_{\rho}^{\frac{7}{10}}(\mathcal{A}_\rho+\alpha_\rho)^{\frac{1}{2}}+\mathcal{G}_{\rho}^{\frac{7}{10}}(\mathcal{B}_\rho+\beta_\rho)^{\frac{1}{2}}\right)+C\frac{\rho^3}{r^2}\mathcal{W}_{\rho}^{\frac{1}{2}}\|\vw\|_{L^\infty_{t,x}(\Omega)}\label{Formula_FirstEstimate}\\
&+C\bigg((\mathcal{A}_\rho+\alpha_\rho)^{\frac{1}{2}}+(\mathcal{B}_\rho+\beta_\rho)^{\frac{1}{2}}\bigg) \left(\frac{\rho^{2}}{r^2} +\frac{\rho}{r}\right) (\alpha_\rho^{\frac{1}{2}}+\beta_\rho^{\frac{1}{2}})\times  \rho^\frac16\|\va\|_{L^6_{t,x}(\Omega)}\notag\\
&+C\frac{\rho^2}{r^2}\mathcal{C}_\rho^{\frac{1}{2}} \gamma_{\rho}^{\frac{1}{2}}
\left((\mathcal{A}_\rho+\alpha_\rho)^{\frac{1}{2}}+(\mathcal{B}_\rho+\beta_\rho)^{\frac{1}{2}}\right)+C\frac{\rho^2}{r}(\alpha_{\rho}+\beta_{\rho})^{\frac{1}{2}}\mathcal{C}_{\rho}^{\frac{1}{2}}.\notag
\end{align}
\end{propo}
\begin{rema}\label{Remarque_HypothesesVW}
Note that the hypothesis $\vw\in L^\infty_{t,x}(\Omega)$ is crucial at this step. It can be relaxed assuming for example $\vw\in L^p_{t}L^q_x(\Omega)$ with $\frac{10}{\tau_0}-1-\frac{2}{p}-\frac{3}{q}>0$ where $\frac{5}{1-\alpha}<\tau_0<\frac{11}{2}$ and $0<\alpha<\frac{1}{12}$ is the exponent of the expected H\"older regularity.
\end{rema}
%%%%%%%%%%%%%%%%%%%%%%%%%%%%%%%%%%%%%%
\textbf{Proof.} It is worth noting here that the structure of this estimate follows closely the one of the local energy inequality given in  (\ref{EstimationEnergieLocale}) and in order to deduce this control, we will start estimating the terms of the right-hand side of (\ref{EstimationEnergieLocale}).
\begin{itemize}
\item Indeed, by the point \emph{4)} of Lemma \ref{testfonc} and using the quantities introduced in Definition \ref{Def_Invariants} we have, for the first term of the right-hand side of (\ref{EstimationEnergieLocale}):
\begin{align*}
\int_{t-\rho^2}^{t+\rho^2}\int_{B_\rho}(\partial_t\phi+\Delta \phi)(|\vu|^2+|\vb|^2+|\vw|^2 )dxds &\le \frac{r^2}{\rho^5}\int_{t-\rho^2}^{t+\rho^2}\int_{B_\rho}(|\vu|^2+|\vb|^2+|\vw|^2 )dxds\nonumber\\
&\le C\frac{r^2}{\rho ^2}( \mathcal{A}_\rho+\mathcal{B}_\rho+ \mathcal{C}_\rho).
\end{align*}
\item For the second term of the right-hand side of (\ref{EstimationEnergieLocale}) we have:
\begin{equation}\label{Estimation_Energie2terme0}
\int_{t-\rho^2}^{t+\rho^2}\int_{B_\rho} [(|\vu|^2+2p)\vb\cdot \vn\phi]dxds\le\int_{t-\rho^2}^{t+\rho^2}\int_{B_\rho}|\vu|^2(\vb\cdot\vn \phi)dxds+C\int_{t-\rho^2}^{t+\rho^2}\int_{B_\rho} |p| |\vb|| \vn\phi|dxds,
\end{equation}
and we will study the two previous terms separately. For the first term of the right-hand side above we introduce the quantity $(|\vu|^2)_\rho$ as the average 
\begin{equation}\label{Def_Moyenneutile}
(|\vu|^2)_\rho=\displaystyle{\frac{1}{|B(x,\rho)|}\int_{B(x,\rho)}|\vu(t,y)|^2dy},
\end{equation} 
and since $\vb$ is divergence free, for any test function $\psi$ compactly supported within $B(x,\rho)$, we have $\displaystyle{\int _{B(x,\rho)} (|\vu|^2)_\rho(\vb \cdot \grad)\psi dx =0}$. Then, since the test function $\phi$ is supported in the parabolic ball $Q_\rho$ (by Lemma \ref{testfonc}) and using Hölder's inequality, it follows that 
\begin{align*}
&\int_{t-\rho^2}^{t+\rho^2}\int_{B_\rho}|\vu|^2(\vb \cdot \grad)\phi dx ds=\int_{t-\rho^2}^{t+\rho^2}\int_{B_\rho} [|\vu|^2-(|\vu|^2)_\rho](\vb\cdot\grad)\phi dx ds\\
&\le \frac{C}{r^2}\int_{t-\rho^2}^{t+\rho^2}\int_{B_\rho} \left||\vu|^2-(|\vu|^2)_\rho\right||\vb|dx ds\le \frac{C}{r^2}\int_{t-\rho^2}^{t+\rho^2}\||\vu|^2-(|\vu|^2)_\rho\|_{L^{\frac{3}{2}}(B_\rho)}\|\vb(s,\cdot)\|_{L^{3}(B_\rho)}ds,
\end{align*}
where we used the fact that $\|\vn \phi\|_{L^\infty}\leq \frac{C}{r^2}$ (by the point \emph{3)} of Lemma \ref{testfonc}). Thus, by the Poincaré inequality  and using the Hölder inequality (in space and time variable), we obtain
\begin{align*}
&\le \frac{C}{r^2}\int_{t-\rho^2}^{t+\rho^2}\|\grad (|\vu(s,\cdot)|^2)\|_{L^{1}(B_\rho)}\|\vb(s,\cdot)\|_{L^{3}(B_\rho)}ds\le \frac{C}{r^2}\int_{t-\rho^2}^{t+\rho^2} \|\vu\|_{L^2(B_\rho)}\|\grad\otimes\vu\|_{L^2(B_\rho)}\|\vb(s,\cdot)\|_{L^{3}(B_\rho)}ds\\
&\le \frac{C}{r^2} \|\vu\|_{L_t^6L_x^2(Q_\rho)}\|\grad \otimes \vu\|_{L_{t,x}^2(Q_\rho)} \|\vb\|_{L_{t,x}^3(Q_\rho)}\leq C\frac{\rho^2}{r^2}\mathcal{A}^{\frac{1}{2}}\alpha_\rho^{\frac12}\zeta_\rho^{\frac13},
\end{align*}
since by the Definition \ref{Def_Invariants} we have $\|\vu\|_{L_t^6L_x^2(Q_\rho)}\le C \rho ^{\frac{1}{3}}\|\vu\|_{L_t^{\infty}L_x^2(Q_\rho)}\le C\rho ^{\frac{5}{6}}\mathcal{A}^{\frac{1}{2}}$, 
$\|\grad \otimes \vu\|_{L_{t,x}^2(Q_\rho)}\leq\rho^{\frac12}\alpha_\rho^{\frac12}$ and $\|\vb\|_{L_{t,x}^3(Q_\rho)}\leq \rho^{\frac23}\zeta_\rho^{\frac13}$. Using the second inequality of the Lemma \ref{Lem_norml3} we obtain
\begin{equation}\label{Estimation_Energie2terme0A}
\int_{t-\rho^2}^{t+\rho^2}\int_{B_\rho}|\vu|^2(\vu \cdot \grad)\phi dxds\le C\frac{\rho^2}{r^2} \mathcal{A}^{\frac{1}{2}}\alpha_\rho^{\frac12}(\mathcal{B}_\rho+\beta_\rho)^{\frac12}\leq C\frac{\rho^2}{r^2} \alpha_\rho^{\frac{1}{2}}
 (\mathcal{A}_\rho+\mathcal{B}_\rho+\beta_\rho),
 \end{equation}
and this control ends the study of the first term of the right-hand side of (\ref{Estimation_Energie2terme}). For the second term of (\ref{Estimation_Energie2terme}), we simply write (using the properties of the function $\phi$ given in Lemma \ref{testfonc} as well as the quantities given in Definition \ref{Def_Invariants} and Lemma \ref{Lem_norml3}):
\begin{equation}\label{Estimation_Energie2terme0B}
\int_{t-\rho^2}^{t+\rho^2}\int_{B_\rho} |p| |\vb|| \vn\phi|dxds \le \frac{C}{r^2}\|p\|_{L^{\frac{3}{2}}(Q_\rho)}\|\vb\|_{L^{3}(Q_\rho)}\leq \frac{C}{r^2}(\rho^{\frac43}\mathcal{P}_\rho^{\frac{2}{3}})(\rho^{\frac23}\zeta_\rho^{\frac13})\le C\frac{\rho ^2}{r^2}\mathcal{P}_\rho^{\frac{2}{3}}( \mathcal{B}_\rho+\beta_\rho)^{\frac{1}{2}}.
\end{equation}
With estimates (\ref{Estimation_Energie2terme0A}) and (\ref{Estimation_Energie2terme0B}), coming back to (\ref{Estimation_Energie2terme0}) we finally obtain 
\begin{equation}\label{Estimation_Energie2terme}
\int_{t-\rho^2}^{t+\rho^2}\int_{B_\rho} [(|\vu|^2+2p)\vb\cdot \vn\phi]dxds\leq C\frac{\rho^2}{r^2} \alpha_\rho^{\frac{1}{2}}
 (\mathcal{A}_\rho+\mathcal{B}_\rho+\beta_\rho)+ C\frac{\rho ^2}{r^2}\mathcal{P}_\rho^{\frac{2}{3}}( \mathcal{B}_\rho+\beta_\rho)^{\frac{1}{2}}.
\end{equation}
\item The third term of (\ref{EstimationEnergieLocale}) can be treated in a completely symmetric manner and we have the estimate:
$$\int_{t-\rho^2}^{t+\rho^2}\int_{B_\rho} [(|\vb|^2+2p)\vu\cdot \vn\phi]dxds\leq C\frac{\rho^2}{r^2} \beta_\rho^{\frac{1}{2}}(\mathcal{B}_\rho+\mathcal{A}_\rho+\alpha_\rho)+ C\frac{\rho ^2}{r^2}\mathcal{P}_\rho^{\frac{2}{3}}( \mathcal{A}_\rho+\alpha_\rho)^{\frac{1}{2}}.$$
\item For the fourth term of (\ref{EstimationEnergieLocale}) we have
$$\int_{t-\rho^2}^{t+\rho^2}\int_{B_\rho}(\rot\vw) \cdot[\phi(\vu+\vb)]dxds=\int_{t-\rho^2}^{t+\rho^2}\int_{B_\rho}(\rot\vw) \cdot(\phi\vu)dxds+\int_{t-\rho^2}^{t+\rho^2}\int_{B_\rho}(\rot\vw) \cdot(\phi\vb)dxds,$$
and due to the symmetry of the information available it is enough to study one of the terms above. We thus write, by the properties of the function $\phi$ given in Lemma \ref{testfonc}:
\begin{eqnarray*}
\int_{t-\rho^2}^{t+\rho^2}\int_{B_\rho}(\rot \vw)\cdot(\phi \vu)dxds&\le&\frac{C}{r}\|\rot \vw\|_{L_{t,x}^2(Q_\rho)} \|\vu\|_{L_{t,x}^2(Q_\rho)}  \le\frac{C}{r}\|\grad \otimes \vw\|_{L_{t,x}^2(Q_\rho)} \|\vu\|_{L_{t,x}^2(Q_\rho)}\\
&\leq &\frac{C}{r}\|\grad \otimes \vw\|_{L_{t,x}^2(Q_\rho)} \rho\|\vu\|_{L_{t}^\infty L_x^2(Q_\rho)}\leq \frac{C}{r}(\rho \gamma_\rho)^{\frac{1}{2}}\rho (\rho \mathcal{A}_\rho)^{\frac{1}{2}}=C\frac{\rho^2}{r}\gamma_\rho^{\frac{1}{2}}\mathcal{A}_\rho^{\frac{1}{2}},
\end{eqnarray*}
where we used the fact that $\|\vu\|_{L_{t,x}^2(Q_\rho)}\le C\rho \|\vu\|_{L_t^{\infty}L_x^2(Q_\rho)}$ and the Definition \ref{Def_Invariants}. Thus, with the second term involving $(\rot\vw) \cdot(\phi\vb)$ we finally obtain the estimate:
\begin{align}\label{Estimation_Energie4terme}
\int_{t-\rho^2}^{t+\rho^2}\int_{B_\rho}(\rot \vw)\cdot[\phi(\vu+\vb)]dxds 
\le C\frac{\rho^{2}}{r}\gamma_{\rho}^{\frac{1}{2}}(\mathcal{A}_{\rho}^{\frac{1}{2}}+\mathcal{B}_{\rho}^{\frac{1}{2}}).
\end{align}
\item For the term related with $\vf, \vg$ in (\ref{EstimationEnergieLocale}) we have by the properties of the function $\phi$ given in Lemma \ref{testfonc}:
\begin{align*}
\int_{t-\rho^2}^{t+\rho^2}\int_{B_\rho}[\vf\cdot (\phi \vu)+ \vg\cdot (\phi \vb)]dx ds &\le C \frac{1}{r}  \int_{t-\rho^2}^{t+\rho^2}\int_{B_\rho}|\vf| |\vu|+|\vg||\vb|dx ds\\
&\le\frac{C}{r} \|\vf\|_{L_{t,x}^{\frac{10}{7}}(Q_\rho)} \|\vu\|_{L_{t,x}^{\frac{10}{3}}(Q_\rho)}+\|\vg\|_{L_{t,x}^{\frac{10}{7}}(Q_\rho)} \|\vb\|_{L_{t,x}^{\frac{10}{3}}(Q_\rho)}.
\end{align*}
Recalling the control $\|\vu\|_{L_{t,x}^{\frac{10}{3}} (Q_\rho)}\leq C(\|\vu\|_{L_t^{\infty}L_x^{2} (Q_\rho)}+\|\vn\otimes \vu\|_{L_{t,x}^{2} (Q_\rho)})$ and since we have the identities $\|\vu\|_{L_t^{\infty}L_x^{2} (Q_\rho)}=\rho^{\frac12}\mathcal{A}_\rho^{\frac12}$, $\|\vn\otimes \vu\|_{L_{t,x}^{2} (Q_\rho)}=\rho^{\frac12}\mathcal{B}_\rho^{\frac12}$, $\rho^{\frac{1}{2}}\mathcal{F}_\rho^{\frac{7}{10}}=\|\vf\|_{L_{t,x}^{\frac{10}{7}}(Q_\rho)}$ and $\rho^{\frac{1}{2}}\mathcal{G}\rho^{\frac{7}{10}}=\|\vg\|_{L_{t,x}^{\frac{10}{7}}(Q_\rho)}$, we obtain:
$$\int_{t-\rho^2}^{t+\rho^2}\int_{B_\rho}[\vf\cdot (\phi \vu)+ \vg\cdot (\phi \vb)]dx ds \le C\frac{\rho}{r}\left(\mathcal{F}_{\rho}^{\frac{7}{10}}(\mathcal{A}_\rho+\alpha_\rho)^{\frac{1}{2}}+\mathcal{G}_{\rho}^{\frac{7}{10}}(\mathcal{B}_\rho+\beta_\rho)^{\frac{1}{2}}\right).$$
\item For the sixth term of (\ref{EstimationEnergieLocale}) we have, by the properties of the function $\phi$ given in Lemma \ref{testfonc}, by the H\"older inequalities and by the Definition \ref{Def_Invariants}
$$ \int_{t-\rho^2}^{t+\rho^2}\int_{B_\rho} \left|[div(\vw)(\vn\phi\cdot \vw)](s,x)\right|dxds\le \frac{C}{r^2}\|div(\vw)\|_{L_{t,x}^2(Q_\rho)} \|\vw\|_{L_{t,x}^2(Q_\rho)}\le \frac{C}{r^2}(\rho\mathcal{W}_{\rho})^{\frac{1}{2}}\rho^{\frac{5}{2}}\|\vw\|_{L^\infty_{t,x}(\Omega)},$$
from which we obtain\footnote{Note that for the term $\vw$, a $L^p_tL^q_x$-norm can be used here instead of the $L^\infty_tL^\infty_x$-norm. See Remarks \ref{Remarque_HypothesesVW} and \ref{Remarque_HypothesesVW1}.}:
$$ \int_{t-\rho^2}^{t+\rho^2}\int_{B_\rho} |div(\vw)(\vw\cdot\grad) \phi| dx ds 
 \le C\frac{\rho^3}{r^2}\mathcal{W}_{\rho}^{\frac{1}{2}}\|\vw\|_{L^\infty_{t,x}(\Omega)}.$$
\item For the seventh term of (\ref{EstimationEnergieLocale}) we need to study the following quantity
\begin{equation}\label{Estimation_Energie7terme0}
\int_{t-\rho^2}^{t+\rho^2}\int_{B_\rho}\left|\bigg[[(\va\cdot \vn)(\vu+\vb)]\cdot (\phi (\vu+\vb))\bigg](s,x)\right|dxds.
\end{equation}
By symmetry of the available information over the vector fields $\vu$ and $\vb$, it is enough to consider the term $\displaystyle{\int_{t-\rho^2}^{t+\rho^2}\int_{B_\rho}\left|[(\va\cdot \vn)\vu]\cdot (\phi \vu)\right|dxds}$ and we write, by the H\"older inequalities:
$$\int_{t-\rho^2}^{t+\rho^2}\int_{B_\rho}\left|[(\va\cdot \vn)\vu]\cdot (\phi \vu)\right|dxds\leq \frac{C}{r}\|\va\|_{L^6_{t,x}(Q_\rho)}\|\vn \otimes \vu\|_{L^2_{t,x}(Q_\rho)}\|\vu\|_{L^3_{t,x}(Q_\rho)},$$
but since we have $\|\va\|_{L^6_{t,x}(Q_\rho)}\leq \|\va\|_{L^6_{t,x}(\Omega)}$, $(\rho\alpha_\rho)^{\frac12}=\|\vn \otimes \vu\|_{L^2_{t,x}(Q_\rho)}$ and since we have $\|\vu\|_{L^3_{t,x}(Q_\rho)} = \rho^{\frac23}\lambda_\rho^{\frac13}\leq C \rho^{\frac23} (\mathcal{A}_\rho+\alpha_\rho)^{\frac{1}{2}}$ (by Lemma \ref{Lem_norml3}) we obtain: 
$$\int_{t-\rho^2}^{t+\rho^2}\int_{B_\rho}\left|[(\va\cdot \vn)\vu]\cdot (\phi \vu)\right|dxds\leq C\frac{\rho}{r}\alpha^{\frac12}(\mathcal{A}_\rho+\alpha_\rho)^{\frac12}\rho^{\frac16}\|\va\|_{L^6_{t,x}(\Omega)}.$$
Performing the same computations for the remaining terms of (\ref{Estimation_Energie7terme0}) we have
\begin{eqnarray*}
\int_{t-\rho^2}^{t+\rho^2}\int_{B_\rho}\left|\bigg[[(\va\cdot \vn)(\vu+\vb)]\cdot (\phi (\vu+\vb))\bigg](s,x)\right|dxds\leq C\frac{\rho}{r}(\alpha^{\frac12}+\beta_\rho^{\frac{1}{2}})\notag\\
\times\bigg((\mathcal{A}_\rho+\alpha_\rho)^{\frac12}
+(\mathcal{B}_\rho+\beta_\rho)^{\frac12}\bigg)\rho^{\frac16}\|\va\|_{L^6_{t,x}(\Omega)}.
\end{eqnarray*}
\item The eighth term of (\ref{EstimationEnergieLocale}) is $\displaystyle{\int_{t-\rho^2}^{t+\rho^2}\int_{B_\rho}\bigg[\big((\vu+\vb)\cdot \vn\big)\big(\phi(\vu+\vb)\big)\cdot\va\bigg] (s,x)dxds}$ and again, it is enough to study the following generic term which contains the term $\big(\vu\cdot \vn\big)\big(\phi\vu\big)\cdot\va$ and we have
\begin{eqnarray*}
\int_{t-\rho^2}^{t+\rho^2}\int_{B_\rho}\bigg[\big(\vu\cdot \vn\big)\big(\phi\vu\big)\cdot\va\bigg]dxds\leq \|\vu\|_{L^3_{t,x}(Q_\rho)}\left(\|\vn\phi\|_{L^\infty_{t,x}}\|\vu\|_{L^2_{t,x}(Q_\rho)}+\|\phi\|_{L^\infty_{t,x}}\|\vn\otimes\vu\|_{L^2_{t,x}(Q_\rho)}\right)\|\va\|_{L^6_{t,x}(Q_\rho)}\\
\leq  \|\vu\|_{L^3_{t,x}(Q_\rho)}\left(\frac{C}{r^2}\rho\|\vu\|_{L^2_{t}L^6_{t}(Q_\rho)}+\frac{C}{r}\|\vu\|_{L^2_{t}\dot{H}^1_x(Q_\rho)}\right)\|\va\|_{L^6_{t,x}(Q_\rho)}
\end{eqnarray*}
$$\leq  \rho^{\frac23} (\mathcal{A}_\rho+\alpha_\rho)^{\frac{1}{2}}\left(\frac{C}{r^2}\rho(\rho\alpha_\rho)^{\frac12}+\frac{C}{r}(\rho\alpha_\rho)^{\frac12}\right)\|\va\|_{L^6_{t,x}(\Omega)}
\leq  C\left(\frac{\rho^{2}}{r^2} +\frac{\rho}{r}\right) \alpha_\rho^{\frac{1}{2}}(\mathcal{A}_\rho+\alpha_\rho)^{\frac{1}{2}}\rho^{\frac16}\|\va\|_{L^6_{t,x}(\Omega)},$$
where we used the properties of the function $\phi$ given in Lemma \ref{testfonc}, the Definition \ref{Def_Invariants} and the Lemma \ref{Lem_norml3}. Thus, considering the remaining terms we can write
\begin{eqnarray*}
\int_{t-\rho^2}^{t+\rho^2}\int_{B_\rho}\bigg[\big((\vu+\vb)\cdot \vn\big)\big(\phi(\vu+\vb)\big)\cdot\va\bigg]dxds&\leq &C\bigg((\mathcal{A}_\rho+\alpha_\rho)^{\frac{1}{2}}+(\mathcal{B}_\rho+\beta_\rho)^{\frac{1}{2}}\bigg) \left(\frac{\rho^{2}}{r^2} +\frac{\rho}{r}\right) (\alpha_\rho^{\frac{1}{2}}+\beta_\rho^{\frac{1}{2}})\notag\\
&&\times  \rho^\frac16\|\va\|_{L^6_{t,x}(\Omega)}.
\end{eqnarray*}
\item For the ninth term of (\ref{EstimationEnergieLocale}) we have to consider the quantity $\displaystyle{\int_{t-\rho^2}^{t+\rho^2}\int_{B_\rho}[|\vw|^2(\vu+\vb)\cdot \grad\phi](s,x) dxds}$ which has the same structure of the first term of the right-hand side of  (\ref{Estimation_Energie2terme0}) and thus, by the same arguments we obtain 
$$\int_{t-\rho^2}^{t+\rho^2}\int_{B_\rho}[|\vw|^2(\vu+\vb)\cdot \grad\phi]dxds \le\frac{\rho^2}{r^2}\mathcal{C}_\rho^{\frac{1}{2}} \gamma_{\rho}^{\frac{1}{2}}
\left((\mathcal{A}_\rho+\alpha_\rho)^{\frac{1}{2}}+(\mathcal{B}_\rho+\beta_\rho)^{\frac{1}{2}}\right).$$
\item The last term of (\ref{EstimationEnergieLocale}) is given by the expression $\displaystyle{\int_{t-\rho^2}^{t+\rho^2}\int_{B_\rho}[\rot(\vu+\vb)]\cdot(\phi \vw)dxds}$ and we remark that it is of the same structure of the term (\ref{Estimation_Energie4terme}), so we obtain
$$\int_{t-\rho^2}^{t+\rho^2}\int_{B_\rho}|\rot(\vu+\vb)\cdot(\phi \vw)|dxds\le C\frac{\rho^2}{r}(\alpha_{\rho}+\beta_{\rho})^{\frac{1}{2}}\mathcal{C}_{\rho}^{\frac{1}{2}}.$$
\end{itemize}
Once we have estimated all these terms, in order to obtain (\ref{Formula_FirstEstimate}) it is enough to gather them: doing so we obtain an uniform estimate with respect to the radius $r$ and to end the proof we remark that the left-hand side of the energy inequality is controlled (using the quantities given in Definition \ref{Def_Invariants}) by the left-hand side of (\ref{Formula_FirstEstimate}).\hfill$\blacksquare$ 
%%%%%%%%%%%%%%%%%%%%%%%%%%%%%%%%%%%%%%
\subsection{A second estimate}
The control obtained in the previous section is crucial but it is not enough to our purposes as we need to obtain a deeper control over the pressure. For this
%%%%%%%%%%%%%%%%%%%%%%%%%%%%%%%%%%%%%%
 \begin{lem}\label{pre}
For some $0< \sigma <\frac{1}{2} $ and for a parabolic ball $Q_\sigma$ of the form (\ref{Def_BoulesQ}), we have the following estimate on the pressure
$$\|p\|_{L_{t,x}^{\frac{3}{2}}(Q_\sigma)}\leq C\sigma^{\frac{1}{3}}\|\vu\|_{L^{\infty}_tL^{2}_x (Q_1)} \|\vb\|_{L_{t}^{2}\dot{H}^1_x (Q_1)}
+C \sigma^2\big(\|\vu\|_{L^2_t\dot{H}^1_x(Q_1)}+\|\vb\|_{L^2_t\dot{H}^1_x(Q_1)}\big) \|\va\|_{L^{6}_{t,x}(Q_1)}
 +\sigma^2\|p\|_{L^{\frac{3}{2}}_{t,x}(Q_1)}$$
\end{lem}
\begin{rema}
For the time being we assume the controls of the right-hand side of the previous estimate. We will see later on, by a suitable change of variables, how to recover the information over the balls $Q_r\subset \Omega$.
\end{rema}
%%%%%%%%%%%%%%%%%%%%%%%%%%%%%%%%%%%%%%
\textbf{Proof.}
First, we introduce a smooth function $\eta:\R\longrightarrow [0,1]$
supported by the ball $B_1$ such that $\eta =1$ on the ball $B_{\frac{3}{5}}$
and $\eta=0$ outside the ball $B_{\frac{4}{5}}$. By a straightforward calculation we have the identity $-\Delta(\eta p)=-\eta \Delta p+(\Delta \eta)p-2\sum_{i=1}^3\partial_i((\partial_i\eta)p)$ and we thus have
\begin{equation}\label{EstimationPression1}
\|p\|_{L_{t,x}^{\frac{3}{2}}(Q_\sigma)}\le \underbrace{\left\|\frac{(-\eta \Delta p)}{(-\Delta)}\right\|_{L^{\frac{3}{2}}_{t,x}(Q_\sigma)}}_{(a)}+\underbrace{\left\|\frac{(\Delta \eta) p}{(-\Delta)}\right\|_{L^{\frac{3}{2}}_{t,x}(Q_\sigma)}}_{(b)}+2\sum_{i=1}^3\underbrace{\left\|\frac{\partial_i((\partial_i \eta)p)}{(-\Delta)}\right\|_{L^{\frac{3}{2}}_{t,x}(Q_\sigma)}}_{(c)}.
\end{equation}
\begin{itemize}
\item For the first term of (\ref{EstimationPression1}) above, we use the expression of the pressure given in (\ref{FormulePressionIntro}) which allows us to write $2\Delta p=-div((\vb\cdot \vn)\vu)-div((\vu\cdot \vn)\vb)-div(div((\vu+\vb) \otimes \va+\va \otimes (\vu+\vb)))$ and, due to the fact that $div(\vu)=div(\vb)=div(\va)=0$, we obtain the expression
$$\Delta p=-\sum^3_{i,j= 1}\partial_i \partial_j (u_i b_j)+\sum^3_{i,j= 1}\partial_i \partial_j\left((u_i+b_i)a_j+a_i(u_j+b_j)\right),$$
from which one gets 
\begin{equation}\label{EstimationPression10}
\left\|\frac{(-\eta \Delta p)}{(-\Delta)}\right\|_{L^{\frac{3}{2}}_{t,x}(Q_\sigma)}\leq \sum^3_{i,j= 1}\underbrace{\left\|\frac{\eta\partial_i \partial_j (u_i b_j)}{(-\Delta)}\right\|_{L^{\frac{3}{2}}_{t,x}(Q_\sigma)}}_{(a.1)}+\underbrace{\left\|\frac{\eta\partial_i \partial_j\left((u_i+b_i)a_j+a_i(u_j+b_j)\right)}{(-\Delta)}\right\|_{L^{\frac{3}{2}}_{t,x}(Q_\sigma)}}_{(a.2)}.
\end{equation}
In order to study the term $(a.1)$ above, we introduce the quantity $\mathfrak{U}_{i,j} = u_i (b_j - (b_j)_1)$ where $ (b_j)_1$ is the average of $b_j$ over the ball of radius $1$ (recall the definition (\ref{Def_Moyenneutile})) and since $\vu$ is divergence free we have the identity $\displaystyle{\sum^3_{i,j= 1}}\partial_i \partial_j (u_i b_j)=\displaystyle{\sum^3_{i,j= 1}}\partial_i \partial_j \mathfrak{U}_{i,j}$. Noting now that we also have the identity $\eta\partial_i \partial_j \mathfrak{U}_{i,j}=\partial_i \partial_j(\eta \mathfrak{U}_{i,j}) - \partial_i \big((\partial_j \eta) \mathfrak{U}_{i,j} \big) - \partial_j \big((\partial_i \eta) \mathfrak{U}_{i,j} \big) + 2(\partial_i \partial_j \eta) \mathfrak{U}_{i,j}$, we obtain
\begin{eqnarray}
\left\|\frac{\eta\partial_i \partial_j (u_i b_j)}{(-\Delta)}\right\|_{L^{\frac{3}{2}}_{t,x}(Q_\sigma)}&\leq& \left\|\frac{\partial_i \partial_j(\eta \mathfrak{U}_{i,j})}{(-\Delta)}\right\|_{L^{\frac{3}{2}}_{t,x}(Q_\sigma)}+\left\|\frac{\partial_i \big((\partial_j \eta) \mathfrak{U}_{i,j} \big)}{(-\Delta)}\right\|_{L^{\frac{3}{2}}_{t,x}(Q_\sigma)}\label{EstimationPression2InegaliteA1}\\
&&+\left\|\frac{\partial_j \big((\partial_i \eta) \mathfrak{U}_{i,j} \big) }{(-\Delta)}\right\|_{L^{\frac{3}{2}}_{t,x}(Q_\sigma)}+C\left\|\frac{(\partial_i \partial_j \eta) \mathfrak{U}_{i,j}}{(-\Delta)}\right\|_{L^{\frac{3}{2}}_{t,x}(Q_\sigma)}.\notag
\end{eqnarray}
The first term of the right-hand side above is easy to control, indeed denoting by $\mathcal{R}_i=\frac{\partial_i}{\sqrt{-\Delta}}$ the usual Riesz transforms on $\mathbb{R}^3$, by the boundedness of these operators in Lebesgue spaces and using the support properties of the auxiliary function $\eta$, we have (recalling that $\mathfrak{U}_{i,j} = u_i (b_j - (b_j)_1)$):
\begin{eqnarray*}
\left\|\frac{\partial_i \partial_j}{(-\Delta)} \eta \mathfrak{U}_{i,j}(t,\cdot) \right\|_{L^{\frac32} (B_\sigma)} &\leq &\|\mathcal{R}_i \mathcal{R}_j (\eta \mathfrak{U}_{i,j} )(t,\cdot) \|_{L^{\frac32} (\R)}  \leq  C\|\eta \mathfrak{U}_{i,j}(t,\cdot)\|_{L^{\frac32} (B_1)}\\ 
&\leq &C\|u_i(t,\cdot) \|_{L^{2}(B_1)} \|b_j(t,\cdot) - (b_j)_1\|_{L^{6} (B_1)}\leq C\|\vu(t,\cdot)\|_{L^{2} (B_1)} \|\vn \otimes \vb(t,\cdot)\|_{L^{2} (B_1)},
\end{eqnarray*}
where we used H\"older and Poincaré inequalities in the last line. Now taking the $L^{\frac32}$-norm in the time variable of the previous inequality we obtain
\begin{equation}\label{EstimationPression2InegaliteA101}
\left\|\frac{\partial_i \partial_j}{(-\Delta)} \eta \mathfrak{U}_{i,j} \right\|_{L_{t,x}^{\frac32} (Q_\sigma)}\leq C\sigma^{\frac{1}{3}}\|\vu\|_{L^{\infty}_tL^{2}_x (Q_1)} \|\vn \otimes \vb\|_{L_{t,x}^{2} (Q_1)}. 
\end{equation}
The second and the third term of the right-hand side of (\ref{EstimationPression2InegaliteA1}) are treated in a similar manner., so we will only consider one of them. Since $\partial_i \eta$ vanishes on $B_{\frac35} \cup B^c_{\frac45}$ and since $B_\sigma \subset B_{\frac12}\subset B_{\frac{3}{5}}$, with the integral representation of the operator $\frac{\partial_i}{(- \Delta )}$ we have for the second term of (\ref{EstimationPression2InegaliteA1}) the inequalities (taking into account only the space variable):
$$\left\|\frac{\partial_i}{(- \Delta )}\big((\partial_j\eta)\mathfrak{U}_{i,j}\big)(t,\cdot)\right\|_{L^{\frac32}(B_\sigma)} \leq C\sigma^2\left\|\frac{\partial_i}{(- \Delta )}\big((\partial_j\eta)\mathfrak{U}_{i,j}\big)(t,\cdot)\right\|_{L^{\infty}(B_\sigma)}$$
\begin{eqnarray}
&\leq &C \, \sigma^2 \left\|\int_{\{\frac35<|y|<\frac45\}} \frac{x_i - y_i }{|x-y|^3} \big((\partial_j \eta) \mathfrak{U}_{i,j} \big)(t,y)  \, dy\right\|_{L^{\infty}(B_\sigma)} \leq C \, \sigma^2 \|\mathfrak{U}_{i,j}(t,\cdot)\|_{L^{1} (B_1)}\label{KernelEstimate1}\\
& \leq &C \, \sigma^2 \|u_i(t,\cdot) \|_{L^{2} (B_1)} \|b_j(t,\cdot) - (b_j)_1\|_{L^{2}(B_1)} \leq C \sigma^2  \|\vu(t,\cdot)\|_{L^{2} (B_1)}\|\vn \otimes \vb(t,\cdot)\|_{L^{2} (B_1)},\notag
\end{eqnarray}
where we used the same ideas as previously. Taking the $L^{\frac32}$-norm in the time variable, we obtain
\begin{equation}\label{EstimationPression2InegaliteA102}
\left\|\frac{\partial_i}{(- \Delta )}\big((\partial_j\eta)\mathfrak{U}_{i,j}\big)\right\|_{L^{\frac32}_{t,x}(Q_\sigma)} \leq C\sigma^{\frac{7}{3}}\|\vu\|_{L^{\infty}_tL^{2}_x (Q_1)} \|\vb\|_{L_{t}^{2}\dot{H}^1_x (Q_1)}\leq C\sigma^{\frac{1}{3}}\|\vu\|_{L^{\infty}_tL^{2}_x (Q_1)} \|\vb\|_{L_{t}^{2}\dot{H}^1_x (Q_1)}. 
\end{equation}
(since $\sigma^{\frac{7}{3}}\leq \sigma^{\frac{1}{3}}$ as we have $0<\sigma<\frac12$). For the last term of (\ref{EstimationPression2InegaliteA1}), we recall that the convolution kernel associated to the operator $\frac{1}{(-\Delta)}$ is $\frac{C}{|x|}$, and thus following the same ideas we have the inequality
\begin{equation}\label{EstimationPression2InegaliteA103}
\left\|\frac{(\partial_i \partial_j \eta) N_{i,j}}{(-\Delta)}\right\|_{L_{t,x}^{\frac 32}(Q_\sigma)}\leq C\sigma^{\frac{1}{3}}\|\vu\|_{L^{\infty}_tL^{2}_x (Q_1)} \|\vn \otimes \vb\|_{L_{t,x}^{2} (Q_1)}.
\end{equation}
Thus, gathering the estimates (\ref{EstimationPression2InegaliteA101}), (\ref{EstimationPression2InegaliteA102}) and (\ref{EstimationPression2InegaliteA103}) and coming back to (\ref{EstimationPression2InegaliteA1}) we finally obtain
\begin{equation}\label{EstimationPression2InegaliteA1Final}
(a.1)=\left\|\frac{\eta\partial_i \partial_j (u_i b_j)}{(-\Delta)}\right\|_{L^{\frac{3}{2}}_{t,x}(Q_\sigma)}\leq C\sigma^{\frac{1}{3}}\|\vu\|_{L^{\infty}_tL^{2}_x (Q_1)} \|\vb\|_{L_{t}^{2}\dot{H}^1_x (Q_1)}.
\end{equation}
We study now the term $(a.2)$ of (\ref{EstimationPression10}). Due to the symmetry of the quantity $\eta\partial_i\partial_j\left((u_i+b_i)a_j+a_i(u_j+b_j)\right)$ it is enough to treat one term of the form $\eta\partial_i\partial_j(u_ia_j)$ for which we use as before the identity
$\eta\partial_i \partial_j (u_ia_j)=\partial_i \partial_j(\eta (u_ia_j)) - \partial_i \big((\partial_j \eta) (u_ia_j) \big) - \partial_j \big((\partial_i \eta) (u_ia_j) \big) + 2(\partial_i \partial_j \eta) (u_ia_j)$ and we have 
\begin{eqnarray}
\left\|\frac{\eta\partial_i \partial_j (u_i a_j)}{(-\Delta)}\right\|_{L^{\frac{3}{2}}_{t,x}(Q_\sigma)}&\leq& \left\|\frac{\partial_i \partial_j (\eta(u_i a_j))}{(-\Delta)}\right\|_{L^{\frac{3}{2}}_{t,x}(Q_\sigma)}+\left\|\frac{\partial_i ((\partial_j\eta) (u_i a_j))}{(-\Delta)}\right\|_{L^{\frac{3}{2}}_{t,x}(Q_\sigma)}\notag\\
&&+\left\|\frac{\partial_j ((\partial_i\eta) (u_i a_j))}{(-\Delta)}\right\|_{L^{\frac{3}{2}}_{t,x}(Q_\sigma)}+2\left\|\frac{ (\partial_i\partial_j\eta) (u_i a_j)}{(-\Delta)}\right\|_{L^{\frac{3}{2}}_{t,x}(Q_\sigma)}.\label{EstimationPression2InegaliteA2}
\end{eqnarray}
For the first term of the right-hand side above, introducing the Riesz transforms and using the support properties of the localizing function $\eta$ we have:
$$ \left\|\frac{\partial_i \partial_j (\eta(u_i a_j))}{(-\Delta)}\right\|_{L^{\frac{3}{2}}(B_\sigma)}= \left\|\mathcal{R}_i\mathcal{R}_j (\eta(u_i a_j))\right\|_{L^{\frac{3}{2}}(B_\sigma)}\leq \left\|\eta(u_i a_j)\right\|_{L^{\frac{3}{2}}(B_1)},$$
now taking the $L^\frac{3}{2}$-norm in the time variable and applying the H\"older inequalities (in space and then in time) we have 
\begin{equation}\label{EstimationPression2InegaliteA201}
\left\|\frac{\partial_i \partial_j (\eta(u_i a_j))}{(-\Delta)}\right\|_{L^{\frac{3}{2}}_{t,x}(Q_\sigma)}\leq C\|u_i\|_{L^{2}_tL^{6}_x(Q_1)}\|a_j\|_{L^{6}_tL^{2}_x(Q_1)}\leq C\|u_i\|_{L^{2}_t\dot{H}^{1}_x(Q_1)}\|a_j\|_{L^{6}_tL^{6}_x(Q_1)},
\end{equation}
where in the last estimate we used the local inclusion between Lebesgue spaces. Now, just as before (when studying (\ref{EstimationPression2InegaliteA1})), the second and the third term of (\ref{EstimationPression2InegaliteA2}) can be treated in a similar manner and we will just study the second term and we have, following the same ideas displayed in (\ref{KernelEstimate1}):
$$\left\|\frac{\partial_i ((\partial_j\eta) (u_i a_j))}{(-\Delta)}\right\|_{L^{\frac{3}{2}}(B_\sigma)}\leq C\sigma^2\|u_i a_j\|_{L^1(B_1)}\leq C\sigma^2\|u_i\|_{L^6(B_1)} \|a_j\|_{L^{\frac{6}{5}}(B_1)}\leq C\sigma^2\|u_i\|_{\dot{H}^1(B_1)} \|a_j\|_{L^{6}(B_1)},$$
and with an integration in the time variable applying the H\"older inequalities it comes
\begin{equation}\label{EstimationPression2InegaliteA202}
\left\|\frac{\partial_i ((\partial_j\eta) (u_i a_j))}{(-\Delta)}\right\|_{L^{\frac{3}{2}}_{t,x}(Q_\sigma)}\leq C\sigma^2\|u_i\|_{L^2_t\dot{H}^1_x(Q_1)} \|a_j\|_{L^{6}_tL^{6}_x(Q_1)}.
\end{equation}
For the last term of (\ref{EstimationPression2InegaliteA2}) we proceed in a similar manner noting that the convolution kernel associated to the operator $\frac{1}{(-\Delta)}$ is $\frac{C}{|x|}$ and due to the support properties of the localizing function $\eta$ we can write $\left\|\frac{ (\partial_i\partial_j\eta) (u_i a_j)}{(-\Delta)}\right\|_{L^{\frac{3}{2}}(B_\sigma)}\leq C \sigma^2\|u_i a_j\|_{L^1(B_\sigma)}$ from which we easily deduce the estimate
\begin{equation}\label{EstimationPression2InegaliteA203}
\left\|\frac{ (\partial_i\partial_j\eta) (u_i a_j)}{(-\Delta)}\right\|_{L^{\frac{3}{2}}_{t,x}(Q_\sigma)}\leq C \sigma^2\|u_i\|_{L^2_t\dot{H}^1_x(Q_1)} \|a_j\|_{L^{6}_tL^{6}_x(Q_1)}.
\end{equation}
Thus, gathering the estimates (\ref{EstimationPression2InegaliteA201}), (\ref{EstimationPression2InegaliteA202}) and (\ref{EstimationPression2InegaliteA203}) and coming back to the inequality (\ref{EstimationPression2InegaliteA2}) we obtain:
$$\left\|\frac{\eta\partial_i \partial_j (u_i a_j)}{(-\Delta)}\right\|_{L^{\frac{3}{2}}_{t,x}(Q_\sigma)}\leq C \sigma^2\|\vu\|_{L^2_t\dot{H}^1_x(Q_1)}\|\va\|_{L^{6}_tL^{6}_x(Q_1)}.$$
Now, considering the terms of the form $\eta\partial_i\partial_j(b_ia_j)$ we have
\begin{equation}\label{EstimationPression2InegaliteA2Final}
\left\|\frac{\eta\partial_i \partial_j\left((u_i+b_i)a_j+a_i(u_j+b_j)\right)}{(-\Delta)}\right\|_{L^{\frac{3}{2}}_{t,x}(Q_\sigma)}\leq C \sigma^2\big(\|\vu\|_{L^2_t\dot{H}^1_x(Q_1)}+\|\vb\|_{L^2_t\dot{H}^1_x(Q_1)}\big) \|\va\|_{L^{6}_tL^{6}_x(Q_1)}.
\end{equation}
With the previous estimates for the terms $(a.1)$ and $(a.2)$  given in (\ref{EstimationPression2InegaliteA1Final}) and (\ref{EstimationPression2InegaliteA2Final}), respectively, and coming back to the expression (\ref{EstimationPression10}) we obtain
\begin{equation}\label{EstimationPression10Final}
\begin{split}
\left\|\frac{(-\eta \Delta p)}{(-\Delta)}\right\|_{L^{\frac{3}{2}}_{t,x}(Q_\sigma)}&\leq C\sigma^{\frac{1}{3}}\|\vu\|_{L^{\infty}_tL^{2}_x (Q_1)} \|\vb\|_{L_{t}^{2}\dot{H}^1_x (Q_1)}\\
&+C \sigma^2\big(\|\vu\|_{L^2_t\dot{H}^1_x(Q_1)}+\|\vb\|_{L^2_t\dot{H}^1_x(Q_1)}\big) \|\va\|_{L^{6}_tL^{6}_x(Q_1)}.
\end{split}
\end{equation}
\item We can now study the term $(b)$ of (\ref{EstimationPression1}) and we have (proceeding just like in (\ref{KernelEstimate1}) with the kernel of the operator $\frac{1}{(-\Delta)}$ and the support properties of $\eta$): $\left\|\frac{(\Delta \eta) p}{(-\Delta)}\right\|_{L^{\frac{3}{2}}(B_\sigma)}\leq C\sigma^2\|p\|_{L^1(B_1)}\leq C\sigma^2\|p\|_{L^{\frac{3}{2}}(B_1)}$ and taking the $L^\frac{3}{2}$-norm in the time variable it comes
\begin{equation}\label{EstimationPression20Final}
\left\|\frac{(\Delta \eta) p}{(-\Delta)}\right\|_{L^{\frac{3}{2}}_{t,x}(Q_\sigma)}\leq C\sigma^2\|p\|_{L^{\frac{3}{2}}_{t,x}(Q_1)}.
\end{equation}
\item The last term of (\ref{EstimationPression1}) can be easily treated by following the same ideas displayed previously and we obtain
\begin{equation}\label{EstimationPression30Final}
\left\|\frac{\partial_i((\partial_i \eta)p)}{(-\Delta)}\right\|_{L^{\frac{3}{2}}_{t,x}(Q_\sigma)}\leq C\sigma^2\|p\|_{L^{\frac{3}{2}}_{t,x}(Q_1)}.
\end{equation}
\end{itemize}
To end the proof of the Lemma, it is enough to use the estimates (\ref{EstimationPression10Final}), (\ref{EstimationPression20Final}) and (\ref{EstimationPression30Final}) in (\ref{EstimationPression1}) to obtain the wished inequality. \hfill $\blacksquare$\\

Now, using a scaling argument and the control given in the last lemma, we have the
following proposition. 
%%%%%%%%%%%%%%%%%%%%%%%%%%%%%%%%%%%%%%
\begin{propo}[Second estimate]\label{Propo_SecondEstimate} With the quantities
defined in Definition \ref{Def_Invariants}, under the hypotheses of Theorem
 \ref{HolderRegularity_theorem}, and for $0<r<\frac{\rho}{2}\le 1$, we have the estimate
\begin{align}
\mathcal{P}_r^{\frac{2}{3}}\le C \biggl(\left(\frac{\rho}{r}\right)\mathcal{A}_\rho^{\frac{1}{2}}\beta_\rho^{\frac{1}{2}}+\left(\frac{r}{\rho}\right)^{\frac23}\left(\alpha_\rho^{\frac{1}{2}}+\beta_\rho^{\frac{1}{2}}\right)\|\va\|_{L_{t,x}^6(Q_\rho)}+ \left(\frac{r}{\rho}\right)^{\frac23}\mathcal{P}_\rho^{\frac{2}{3}} \biggr)\label{estimate2}
\end{align}
\end{propo}
%%%%%%%%%%%%%%%%%%%%%%%%%%%%%%%%%%%%%%
\textbf{Proof.} Set $\sigma=\frac{r}{\rho}$ and consider the following functions 
$$p_\rho(t,x)=p(\rho^2t,\rho x),\quad \vu_\rho(t,x)=\vu(\rho^2t,\rho x), \quad b_\rho(t,x)=\vb(\rho^2t,\rho x)\quad \mbox{ and }\quad\va_\rho(t,x)= \va(\rho^2t,\rho x),$$ 
thus, by Lemma \ref{pre} and using the rescaled function above we obtain
\begin{align*}
\rho^{-\frac{10}{3}}\|p\|_{L_{t,x}^{\frac{3}{2}}(Q_r)}
&\le C\biggl(\left(
\frac{r}{\rho}\right)^{\frac{1}{3}}\left(\rho^{-\frac{3}{2}}\|\vu\|_{L_t^{\infty} L_x^2(Q_\rho)}\rho^{-\frac{3}{2}}\|\vb\|_{L_{t}^2\dot{H}^1_x(Q_\rho)}\right)\\
&+\left(\frac{r}{\rho}\right)^{2}\biggl(\rho^{-\frac{3}{2}} \|\vu\|_{L_{t}^2 \dot{H}_{x}^{1}(Q_\rho)}+\rho^{-\frac{3}{2}} \|\vb\|_{L_{t}^2 \dot H_{x}^{1}(Q_\rho)}\biggr)\rho^{-\frac56} \|\va\|_{L_{t,x}^6(Q_\rho)}+ \left(\frac{r}{\rho}\right)^{2} \rho ^{-\frac{10}{3}}
 \|p\|_{L_{t,x}^{\frac{3}{2}}(Q_\rho)} \biggr).
\end{align*}
Now, recalling that, by the Definition \ref{Def_Invariants} (see also Remark \ref{Remarque_DefinitionScaled}) we have the notation $r^{\frac{4}{3}}\mathcal{P}_r^{\frac{2}{3}}=\|p\|_{L_{t,x}^{\frac{3}{2}}(Q_r)}$,  $\rho^{\frac12}\mathcal{A}_\rho^{\frac12}= \|\vu\|_{L_{t}^\infty L_{x}^{2}(Q_\rho)}$, $\rho^{\frac12}\alpha_\rho^{\frac12}= \|\vu\|_{L_{t}^2 \dot{H}_{x}^{1}(Q_\rho)}$ and $\rho^{\frac12}\beta_\rho^{\frac12}= \|\vb\|_{L_{t}^2 \dot{H}_{x}^{1}(Q_\rho)}$, thus we can write
$$\frac{r^{\frac{4}{3}}}{\rho^{\frac{10}{3}} } \mathcal{P}_r^{\frac{2}{3}}\le C \biggl(\left(\frac{r}{\rho}\right)^{\frac{1}{3}} \left( \rho^{-2}\mathcal{A}_\rho^{\frac{1}{2}}\beta_\rho^{\frac{1}{2}}\right)+\left(\frac{r}{\rho}\right)^{2} \rho^{-\frac{11}{6}}\left(\alpha_\rho^{\frac{1}{2}}+\beta_\rho^{\frac{1}{2}}\right)\|\va\|_{L_{t,x}^6(Q_\rho)}+ \left(\frac{r}{\rho}\right)^{2} \rho ^{-2} \mathcal{P}_\rho^{\frac{2}{3}} \biggr),$$
and we obtain (as $\rho^{-\frac{11}{6}}\leq \rho^{-2}$ since $0<\rho<1$)
$$\mathcal{P}_r^{\frac{2}{3}}\le C \biggl(\left(\frac{\rho}{r}\right)\mathcal{A}_\rho^{\frac{1}{2}}\beta_\rho^{\frac{1}{2}}+\left(\frac{r}{\rho}\right)^{\frac23}\left(\alpha_\rho^{\frac{1}{2}}+\beta_\rho^{\frac{1}{2}}\right)\|\va\|_{L_{t,x}^6(Q_\rho)}+ \left(\frac{r}{\rho}\right)^{\frac23}\mathcal{P}_\rho^{\frac{2}{3}} \biggr),$$
which is the desired estimate. \hfill$\blacksquare$
%%%%%%%%%%%%%%%%%%%%%%%%%%%%%%%%%%%%%%
\section{Inductive Argument}\label{Secc_Inductive}
Once we have obtained the estimates (\ref{Formula_FirstEstimate}) and (\ref{estimate2}) it is possible to perform an inductive argument in order to obtain a (local, parabolic) Morrey information over the variables $\vu$, $\vb$ and $\vw$.
%%%%%%%%%%%%%%%%%%%%%%%%%%%%%%%%%%%%%%
\begin{propo}\label{PropositionFirstMorreySpace}
Let $(\vu, \vb, \vw, p)$ be a suitable solution of the magneto-micropolar equations (\ref{EquationMMP}) over the subset $\Omega$. Under the general assumptions of Theorem \ref{HolderRegularity_theorem}, there exists a positive constant $\epsilon^*$ which depends only on $\tau_a,\tau_b$, $\tau_c=\min\{\tau_{a}, \tau_b\}>\frac{5}{2-\alpha}>\frac{5}{3}$ with $0<\alpha<\frac{1}{12}$ and on $\tau_0$ such that if $(t_0,x_0)\in \Omega$ and
\begin{equation}\label{grad}
\limsup_{r\to 0}\frac{1}{r}\int \int _{]t_0-r^2, t_0+r^2[\times B(x_0,r)}|\grad \otimes \vu|^2+|\grad \otimes \vb|^2+|\grad \otimes \vw|^2dx ds < \epsilon^*,
\end{equation}
then there exists a parabolic neighborhood $Q_{R_1}$ of $(t_0,x_0)$ with
$0<R_1<4{\bf R}$ such that
\begin{equation}\label{ConclusionPropositionFirstMorreySpace}
\mathds{1}_{Q_{R_1}} \vu\in \mathcal{M}_{t,x}^{3,\tau_0},\quad\mathds{1}_{Q_{R_1}}\vb\in \mathcal{M}_{t,x}^{3,\tau_0},\quad\mathds{1}_{Q_{R_1}}\vw\in \mathcal{M}_{t,x}^{3,\tau_0}.
\end{equation}
\end{propo}
%%%%%%%%%%%%%%%%%%%%%%%%%%%%%%%%%%%%%%
Note that the conclusion of this proposition is exaclty the first hypothesis of the Proposition \ref{HolderRegularityproposition}.\\

\noindent \textbf{Proof.} Recalling that from the global hypothesis of Theorem \ref{HolderRegularity_theorem} we have a local control over the set $\Omega$, thus as we want to obtain a local information and since we assumed $Q_{\bf R}(t_0,x_0)\subset \Omega$ and by the definition of Morrey
spaces, we only need to prove that there exists a radius $R_1$ small enough such that
for all $0<r<R_1$ and for all $(t,x)\in Q_{R_1}(t_0,x_0)$ we have the following control
\begin{equation}\label{morrey1}
\iint_{Q_r} |\vu|^3+|\vb|^3+|\vw|^3 dy ds \le C r^{5(1-\frac{3}{\tau_0})}.
\end{equation}
In order to obtain this estimate, we will implement an inductive argument using the averaged quantities introduced in the Definition \ref{Def_Invariants}. Indeed, using the Lemma
\ref{Lem_norml3}, we can write 
\begin{equation*}
\iint_{Q_r} |\vu|^3+|\vb|^3+|\vw|^3 dy ds=r^2(\lambda_r+\zeta_r+\sigma_r) \le C r^2(\mathcal{A}_r+\mathcal{B}_r+\mathcal{C}_r+\alpha_r+\beta_r+ \gamma_r )^{\frac{3}{2}}.
\end{equation*}
Then in order to obtain the control (\ref{morrey1}) for all small $0<r<R_1$, and all point $(t,x)\in Q_{R_1}$, it is enough to show the estimate:
\begin{equation*}
\mathcal{A}_r+\mathcal{B}_r+\mathcal{C}_r+\alpha_r+\beta_r+ \gamma_r \le r^{2(1-\frac{5}{\tau_0})}.
\end{equation*}
Let us introduce the following quantities:
\begin{equation}\label{defA}
\mathbf{A}_r=\frac{1}{r^{2(1-\frac{5}{\tau_0})}}( \mathcal{A}_r+\mathcal{B}_r+\mathcal{C}_r+\alpha_r+\beta_r+ \gamma_r +\mathcal{W}_r)\quad\mbox{and}\quad \mathbf{H}_r=r^{\frac{10}{\tau_0}}\mathcal{H}_r.
\end{equation}
Note that the introduction of the quantity $\mathcal{W}_r$ in the first term above is reminiscent from the estimate (\ref{Formula_FirstEstimate}) obtained previously. Thus to prove (\ref{morrey1}) we only need to show that there exists $0<\kappa<1$ and some $0<R_1<{\bf R}$ such that for all $n\in \mathbb{N}$ and $(t,x)\in Q_{R_1}$, we have 
\begin{equation}\label{iterative}
\mathbf{A}_{\kappa^n R_1}\le C,
\end{equation}
and the idea is to use an inductive argument that ensures that we have these estimates above for all radius of the following type $\kappa^n R_1>0$. 
Remark that due to the definition of the quantity $\mathbf{A}_r$ given in (\ref{defA}), we will also obtain some information over the gradients of $\vu,\vb$ and $\vw$ (see Corollary \ref{corolarioMorrey} below).\\

In order to simplify the arguments, we shall need to introduce the following quantities
\begin{equation}\label{defB}
\mathbf{B}_r=(\alpha_r+\beta_r+ \gamma_r+\mathcal{W}_r),\quad \mathbf{P}_r=\frac{1}{r^{\frac{3}{2}(1-\frac{5}{\tau_0})}}\mathcal{P}_r,\quad \mathbf{D}_r=\frac{1}{3-\frac{5}{\tau_c}}(\mathcal{F}_r^{\frac{7}{10}}+\mathcal{G}_r^{\frac{7}{10}}),
\end{equation}
for some $\tau_c>0$ such that $2+\frac{5}{\tau_0}-\frac{5}{\tau_c}>0$. Our starting point is the estimate (\ref{Formula_FirstEstimate}) obtained previously:
\begin{align}
\mathcal{A}_r+\mathcal{B}_r+\mathcal{C}_r&+\alpha_r+\beta_r+ \gamma_r +\mathcal{W}_r+ r^2\mathcal{H}_r \le C\underbrace{\frac{r^2}{\rho ^2}( \mathcal{A}_\rho+\mathcal{B}_\rho+ \mathcal{C}_\rho)}_{(1)}+C\underbrace{\frac{\rho^2}{r^2} \alpha_\rho^{\frac{1}{2}}(\mathcal{A}_\rho+\mathcal{B}_\rho+\beta_\rho)}_{(2)}\notag\\
&+C\underbrace{\frac{\rho ^2}{r^2}\mathcal{P}_\rho^{\frac{2}{3}}( \mathcal{B}_\rho+\beta_\rho)^{\frac{1}{2}}}_{(3)}+C\underbrace{\frac{\rho^2}{r^2} \beta_\rho^{\frac{1}{2}}(\mathcal{B}_\rho+\mathcal{A}_\rho+\alpha_\rho)}_{(4)}+ C\underbrace{\frac{\rho ^2}{r^2}\mathcal{P}_\rho^{\frac{2}{3}}(\mathcal{A}_\rho+\alpha_\rho)^{\frac{1}{2}}}_{(5)}+ C\underbrace{\frac{\rho^{2}}{r}\gamma_{\rho}^{\frac{1}{2}}(\mathcal{A}_{\rho}^{\frac{1}{2}}+\mathcal{B}_{\rho}^{\frac{1}{2}})}_{(6)}\notag\\
&+C\underbrace{\frac{\rho}{r}\left(\mathcal{F}_{\rho}^{\frac{7}{10}}(\mathcal{A}_\rho+\alpha_\rho)^{\frac{1}{2}}+\mathcal{G}_{\rho}^{\frac{7}{10}}(\mathcal{B}_\rho+\beta_\rho)^{\frac{1}{2}}\right)}_{(7)}+C\underbrace{\frac{\rho^3}{r^2}\mathcal{W}_{\rho}^{\frac{1}{2}}\|\vw\|_{L^\infty_{t,x}(\Omega)}}_{(8)}\label{EstimationPourIteration}\\
&+C\underbrace{\bigg((\mathcal{A}_\rho+\alpha_\rho)^{\frac{1}{2}}+(\mathcal{B}_\rho+\beta_\rho)^{\frac{1}{2}}\bigg) \left(\frac{\rho^{2}}{r^2} +\frac{\rho}{r}\right) (\alpha_\rho^{\frac{1}{2}}+\beta_\rho^{\frac{1}{2}})\times  \rho^\frac16\|\va\|_{L^6_{t,x}(\Omega)}}_{(9)}\notag\\
&+C\underbrace{\frac{\rho^2}{r^2}\mathcal{C}_\rho^{\frac{1}{2}} \gamma_{\rho}^{\frac{1}{2}}\left((\mathcal{A}_\rho+\alpha_\rho)^{\frac{1}{2}}+(\mathcal{B}_\rho+\beta_\rho)^{\frac{1}{2}}\right)}_{(10)}+C\underbrace{\frac{\rho^2}{r}(\alpha_{\rho}+\beta_{\rho})^{\frac{1}{2}}\mathcal{C}_{\rho}^{\frac{1}{2}}}_{(11)}.\notag
\end{align}
Multiplying both sides of the inequality (\ref{EstimationPourIteration}) by $\frac{1}{r^{2(1-\frac{5}{\tau_0})}}$, using the formula (\ref{defA}), we obtain in the left-hand side
$$\frac{1}{r^{2(1-\frac{5}{\tau_0})}}\left(\mathcal{A}_r+\mathcal{B}_r+\mathcal{C}_r+\alpha_r+\beta_r+ \gamma_r +\mathcal{W}_r+ r^2\mathcal{H}_r\right)=\mathbf{A}_r+\mathbf{H}_r.$$
Now we will study each term of the right-hand side above multiplied by $\frac{1}{r^{2(1-\frac{5}{\tau_0})}}$:
\begin{itemize}
\item For the term $(1)$ above we have, using the definition of the quantity $\mathbf{A}_\rho$ given in (\ref{defA}):
$$\frac{1}{r^{2(1-\frac{5}{\tau_0})}}\left(\frac{r^2}{\rho ^2}(\mathcal{A}_\rho+\mathcal{B}_\rho+ \mathcal{C}_\rho)\right) \le \frac{1}{r^{2(1-\frac{5}{\tau_0})}}\frac{r^2}{\rho ^2}\rho^{2(1-\frac{5}{\tau_0})}\mathbf{A}_\rho =\left(\frac{r}{\rho}\right)^{\frac{10}{\tau_0}}\mathbf{A}_\rho.$$
\item For the term $(2)$ of (\ref{EstimationPourIteration}), by the definition of $\mathbf{A}_\rho $ and $\mathbf{B}_\rho $ given in (\ref{defA}) and (\ref{defB}) respectively, we can write
$$\frac{1}{r^{2(1-\frac{5}{\tau_0})}}\left(\frac{\rho^2}{r^2} \alpha_\rho^{\frac{1}{2}}(\mathcal{A}_\rho+\mathcal{B}_\rho+\beta_\rho)\right)\leq \frac{1}{r^{2(1-\frac{5}{\tau_0})}}\left(\frac{\rho^2}{r^2} \mathbf{B}_\rho^{\frac{1}{2}}\rho^{2(1-\frac{5}{\tau_0})}\mathbf{A}_\rho \right)=\left(\frac{\rho}{r}\right)^{4-\frac{10}{\tau_0}}\mathbf{A}_\rho \mathbf{B}_\rho^{\frac{1}{2}}.$$
\item For the term $(3)$ of (\ref{EstimationPourIteration}), using the expressions of $\mathbf{A}_\rho $ and $\mathbf{P}_\rho $ given in (\ref{defA}) and (\ref{defB}) respectively, we have
$$\frac{1}{r^{2(1-\frac{5}{\tau_0})}}\left(\frac{\rho ^2}{r^2}\mathcal{P}_\rho^{\frac{2}{3}}( \mathcal{B}_\rho+\beta_\rho)^{\frac{1}{2}}\right)\leq \frac{1}{r^{2(1-\frac{5}{\tau_0})}}\frac{\rho ^2}{r^2}(\rho^{\frac32(1-\frac{5}{\tau_0})}\mathbf{P}_\rho)^{\frac{2}{3}}(\rho^{2(1-\frac{5}{\tau_0})}\mathbf{A}_\rho)^{\frac{1}{2}}=\left(\frac{\rho}{r}\right)^{4-\frac{10}{\tau_0}}\mathbf{P}_\rho^{\frac{2}{3}}\mathbf{A}_\rho^{\frac{1}{2}}.$$
\item The term $(4)$ of (\ref{EstimationPourIteration}) can be treated in the same manner as the term $(2)$ and we obtain
$$\frac{1}{r^{2(1-\frac{5}{\tau_0})}}\left(\frac{\rho^2}{r^2} \beta_\rho^{\frac{1}{2}}(\mathcal{B}_\rho+\mathcal{A}_\rho+\alpha_\rho)\right)\leq \left(\frac{\rho}{r}\right)^{4-\frac{10}{\tau_0}}\mathbf{A}_\rho \mathbf{B}_\rho^{\frac{1}{2}}.$$
\item The term $(5)$ of (\ref{EstimationPourIteration}) can be treated in the same manner as the term $(3)$ and we obtain
$$\frac{1}{r^{2(1-\frac{5}{\tau_0})}}\left(\frac{\rho ^2}{r^2}\mathcal{P}_\rho^{\frac{2}{3}}( \mathcal{A}_\rho+\alpha_\rho)^{\frac{1}{2}}\right)\leq \left(\frac{\rho}{r}\right)^{4-\frac{10}{\tau_0}}\mathbf{P}_\rho^{\frac{2}{3}}\mathbf{A}_\rho^{\frac{1}{2}}.$$
\item By the definition of $\mathbf{A}_\rho$ and $\mathbf{B}_\rho$ given in (\ref{defA}) and (\ref{defB}) respectively, the term $(6)$ of (\ref{EstimationPourIteration}) can be rewritten as follows
$$\frac{1}{r^{2(1-\frac{5}{\tau_0})}}\left(\frac{\rho^{2}}{r}\gamma_{\rho}^{\frac{1}{2}}(\mathcal{A}_{\rho}^{\frac{1}{2}}+\mathcal{B}_{\rho}^{\frac{1}{2}})\right)\leq \frac{1}{r^{2(1-\frac{5}{\tau_0})}}\left(\frac{\rho^2}{r}\mathbf{B}_\rho^{\frac12}\rho^{(1-\frac{5}{\tau_0})}\mathbf{A}_\rho^{\frac12}\right)\leq \left(\frac{\rho}{r}\right)^{3-\frac{10}{\tau_0}}\rho^{\frac{5}{\tau_0}}\mathbf{A}_\rho^{\frac12}\mathbf{B}_\rho^{\frac12}.$$
\item The term $(7)$ of (\ref{EstimationPourIteration}) is estimate using the definition of $\mathbf{D}_\rho$ given in (\ref{defB}):
$$\frac{1}{r^{2(1-\frac{5}{\tau_0})}}\left(\frac{\rho}{r}\left(\mathcal{F}_{\rho}^{\frac{7}{10}}(\mathcal{A}_\rho+\alpha_\rho)^{\frac{1}{2}}+\mathcal{G}_{\rho}^{\frac{7}{10}}(\mathcal{B}_\rho+\beta_\rho)^{\frac{1}{2}}\right)\right)\leq C \left(\frac{\rho}{r}\right)^{3-\frac{10}{\tau_0}}\rho^{2+\frac{5}{\tau_0}-\frac{5}{\tau_c}}\mathbf{D}_\rho\mathbf{A}_\rho^{\frac{1}{2}}.$$
\item For the term $(8)$ of (\ref{EstimationPourIteration}) we use the definition of $\mathbf{B}_\rho$ given in (\ref{defB}) to obtain:
$$\frac{1}{r^{2(1-\frac{5}{\tau_0})}}\left(\frac{\rho^3}{r^2}\mathcal{W}_{\rho}^{\frac{1}{2}}\|\vw\|_{L^\infty_{t,x}(\Omega)}\right)\leq \left(\frac{\rho}{r}\right)^{4-\frac{10}{\tau_0}} \rho^{\frac{10}{\tau_0}-1}\mathbf{B}_{\rho}^{\frac{1}{2}}\|\vw\|_{L^\infty_{t,x}(\Omega)}.$$
\begin{rema}\label{Remarque_HypothesesVW1}
Note that, following Remark \ref{Remarque_HypothesesVW}, if we assume $\vw \in L^p_{t}L^q_x(\Omega)$ with $\frac{10}{\tau_0}-1-\frac{2}{p}-\frac{3}{q}>0$ (which is possible since $\frac{11}{2}>\tau_0>\frac{5}{1-\alpha}$), then the previous bound is $\left(\frac{\rho}{r}\right)^{4-\frac{10}{\tau_0}} \rho^{\frac{10}{\tau_0}-1-\frac{2}{p}-\frac{3}{q}}\mathbf{B}_{\rho}^{\frac{1}{2}}\|\vw\|_{L^p_{t}L^q_x(\Omega)}$.
\end{rema}
\item Since we have $(\alpha^{\frac12}+\beta_\rho^{\frac{1}{2}})\leq C\mathbf{B}_{\rho}^{\frac12}$, $(\mathcal{A}_\rho+\alpha_\rho)^{\frac12}\leq C\rho^{(1-\frac{5}{\tau_0})}\mathbf{A}_{\rho}^{\frac12}$ and $(\mathcal{B}_\rho+\beta_\rho)^{\frac12}\leq C\rho^{(1-\frac{5}{\tau_0})}\mathbf{A}_{\rho}^{\frac12}$, thus for the term $(9)$ of (\ref{EstimationPourIteration}) we write
\begin{eqnarray*}
\frac{1}{r^{2(1-\frac{5}{\tau_0})}}\bigg((\mathcal{A}_\rho+\alpha_\rho)^{\frac{1}{2}}+(\mathcal{B}_\rho+\beta_\rho)^{\frac{1}{2}}\bigg) \left(\frac{\rho^{2}}{r^2} +\frac{\rho}{r}\right) (\alpha_\rho^{\frac{1}{2}}+\beta_\rho^{\frac{1}{2}})\times  \rho^\frac16\|\va\|_{L^6_{t,x}(\Omega)}\leq \frac{C}{r^{2(1-\frac{5}{\tau_0})}}\big(\rho^{(1-\frac{5}{\tau_0})}\mathbf{A}_{\rho}^{\frac12}\big)\\
\times  \left(\frac{\rho^{2}}{r^2} +\frac{\rho}{r}\right)\mathbf{B}_{\rho}^{\frac12}\rho^\frac16\|\va\|_{L^6_{t,x}(\Omega)}^{\frac12},
\end{eqnarray*}
from which we deduce:
$$\leq C\left(\left(\frac{\rho}{r}\right)^{4-\frac{10}{\tau_0}}+\left(\frac{\rho}{r}\right)^{3-\frac{10}{\tau_0}}\right)\rho^{\frac{5}{\tau_0}-\frac56}\mathbf{A}_\rho^{\frac12}\mathbf{B}_\rho^{\frac12}\|\va\|_{L^6_{t,x}(\Omega)}.$$
\item The term $(10)$ of (\ref{EstimationPourIteration}) is treated as follows: recalling that $\gamma_\rho\leq \mathbf{B}_\rho$ by (\ref{defB}) and since we have $\mathcal{C}_\rho^{\frac12}\leq \rho^{(1-\frac{5}{\tau_0})}\mathbf{A}_{\rho}^{\frac12}$, $(\mathcal{A}_\rho+\alpha_\rho)^{\frac{1}{2}}\leq \rho^{(1-\frac{5}{\tau_0})}\mathbf{A}_{\rho}^{\frac12}$ and $(\mathcal{B}_\rho+\beta_\rho)^{\frac{1}{2}}\leq \rho^{(1-\frac{5}{\tau_0})}\mathbf{A}_{\rho}^{\frac12}$ by (\ref{defA}), then we can write
$$\frac{1}{r^{2(1-\frac{5}{\tau_0})}}\bigg(\frac{\rho^2}{r^2}\mathcal{C}_\rho^{\frac{1}{2}} \gamma_{\rho}^{\frac{1}{2}}\left((\mathcal{A}_\rho+\alpha_\rho)^{\frac{1}{2}}+(\mathcal{B}_\rho+\beta_\rho)^{\frac{1}{2}}\right)\bigg)\leq C\left(\frac{\rho}{r}\right)^{4-\frac{10}{\tau_0}}\mathbf{A}_\rho\mathbf{B}_\rho^{\frac12}.$$
\item The last term of (\ref{EstimationPourIteration}) is easy to estimate as we have $(\alpha_{\rho}+\beta_{\rho})^{\frac{1}{2}}\leq C\mathbf{B}_\rho^{\frac12}$ and $\mathcal{C}_{\rho}^{\frac{1}{2}}\leq \rho^{(1-\frac{5}{\tau_0})}\mathbf{A}_\rho^{\frac12}$, then we have
$$\frac{1}{r^{2(1-\frac{5}{\tau_0})}}\bigg(\frac{\rho^2}{r}(\alpha_{\rho}+\beta_{\rho})^{\frac{1}{2}}\mathcal{C}_{\rho}^{\frac{1}{2}}\bigg)\leq \left(\frac{\rho}{r}\right)^{3-\frac{10}{\tau_0}}\rho^{\frac{5}{\tau_0}}\mathbf{A}_\rho^{\frac12}\mathbf{B}_\rho^{\frac12}.$$
\end{itemize}
Once we have all these estimates for the right-hand side of (\ref{EstimationPourIteration}) we finally obtain the following control
\begin{align}
\mathbf{A}_r +\mathbf{H}_r &\le C\Bigg(\left(\frac{r}{\rho}\right)^{\frac{10}{\tau_0}}\mathbf{A}_\rho+\left(\frac{\rho}{r}\right)^{4-\frac{10}{\tau_0}}\mathbf{A}_\rho \mathbf{B}_\rho^{\frac{1}{2}}+\left(\frac{\rho}{r}\right)^{4-\frac{10}{\tau_0}}\mathbf{P}_\rho^{\frac{2}{3}}\mathbf{A}_\rho^{\frac{1}{2}}\notag\\
&+C\left(\left(\frac{\rho}{r}\right)^{4-\frac{10}{\tau_0}}+\left(\frac{\rho}{r}\right)^{3-\frac{10}{\tau_0}}\right)\rho^{\frac{5}{\tau_0}-\frac56}\mathbf{A}_\rho^{\frac12}\mathbf{B}_\rho^{\frac12}\|\va\|_{L^6_{t,x}(\Omega)}+\left(\frac{\rho}{r}\right)^{3-\frac{10}{\tau_0}}\rho^{\frac{5}{\tau_0}}\mathbf{A}_\rho^{\frac12}\mathbf{B}_\rho^{\frac12}\notag\\
&+\left(\frac{\rho}{r}\right)^{3-\frac{10}{\tau_0}}\rho^{2+\frac{5}{\tau_0}-\frac{5}{\tau_c}}\mathbf{D}_\rho\mathbf{A}_\rho^{\frac{1}{2}}+\left(\frac{\rho}{r}\right)^{4-\frac{10}{\tau_0}} \rho^{\frac{10}{\tau_0}-1}\mathbf{B}_{\rho}^{\frac{1}{2}}\|\vw\|_{L^\infty_{t,x}(\Omega)}\Bigg).\label{estimateA}
\end{align}
%%%%%%%%%%%%%%%%%%%%%%%%%%%%%%%%%%%%%% 
Now, we study the estimate for the pressure (\ref{estimate2}) which is given by the control
$$\mathcal{P}_r^{\frac{2}{3}}\le C \biggl(\left(\frac{\rho}{r}\right)\mathcal{A}_\rho^{\frac{1}{2}}\beta_\rho^{\frac{1}{2}}+\left(\frac{r}{\rho}\right)^{\frac23}\left(\alpha_\rho^{\frac{1}{2}}+\beta_\rho^{\frac{1}{2}}\right)\|\va\|_{L_{t,x}^6(Q_\rho)}+ \left(\frac{r}{\rho}\right)^{\frac23}\mathcal{P}_\rho^{\frac{2}{3}} \biggr),$$
and in the same spirit as before, we will introduce the quantity $\mathbf{P}_r=\frac{1}{r^{\frac{3}{2}(1-\frac{5}{\tau_0})}}\mathcal{P}_r$ given in (\ref{defB}) in the left-hand side above. To this end, we will first rise the inequality above to the power $\frac{3}{2}$ and then we will multiply both sides by $\frac{1}{r^{\frac{3}{2}(1-\frac{5}{\tau_0})}}$ and we have
$$\mathbf{P}_r=\frac{1}{r^{\frac{3}{2}(1-\frac{5}{\tau_0})}}\mathcal{P}_r\le \frac{C}{r^{\frac{3}{2}(1-\frac{5}{\tau_0})}} \biggl(\left(\frac{\rho}{r}\right)^{\frac32}\mathcal{A}_\rho^{\frac{3}{4}}\beta_\rho^{\frac{3}{4}}+\left(\frac{r}{\rho}\right)\left(\alpha_\rho^{\frac{1}{2}}+\beta_\rho^{\frac{1}{2}}\right)^{\frac32}\|\va\|_{L_{t,x}^6(Q_\rho)}^{\frac32}+ \left(\frac{r}{\rho}\right)\mathcal{P}_\rho \biggr).$$
We remark now that we have (by the definition of $\mathbf{A}_\rho$ given in (\ref{defA})):
$$\frac{1}{r^{\frac{3}{2}(1-\frac{5}{\tau_0})}}\left(\frac{\rho}{r}\right)^{\frac{3}{2}}\mathcal{A}_\rho^{\frac{3}{4}} \beta_\rho^{\frac{3}{4}}\leq\frac{1}{r^{\frac{3}{2}(1-\frac{5}{\tau_0})}}\left(\frac{\rho}{r}\right)^{\frac{3}{2}}\rho^{\frac{3}{2}(1-\frac{5}{\tau_0})}(\mathbf{A}_\rho\mathbf{B}_\rho)^{\frac{3}{4}}=\left(\frac{\rho}{r}\right)^{3-\frac{15}{2\tau_0}}(\mathbf{A}_\rho\mathbf{B}_\rho)^{\frac{3}{4}}.$$
Note that we also have (by the definition of $\mathbf{B}_\rho$ given in (\ref{defB}))
$$\frac{1}{r^{\frac{3}{2}(1-\frac{5}{\tau_0})}}\left(\frac{r}{\rho}\right)\left(\alpha_\rho^{\frac{1}{2}}+\beta_\rho^{\frac{1}{2}}\right)^{\frac32}\|\va\|_{L_{t,x}^6(Q_\rho)}^{\frac32}\leq C\left(\frac{\rho}{r}\right)^{\frac{1}{2}-\frac{15}{2\tau_0}}\rho^{\frac{15}{2\tau_0}-\frac32}\mathbf{B}_\rho^{\frac34}\|\va\|_{L_{t,x}^6(Q_\rho)}^{\frac32},$$
and finally we have by the definition of $\mathbf{P}_\rho$ given in (\ref{defB}): $\frac{1}{r^{\frac{3}{2}(1-\frac{5}{\tau_0})}} \left(\frac{r}{\rho}\right)\mathcal{P}_\rho= \left(\frac{\rho}{r}\right)^{\frac{1}{2}-\frac{15}{2\tau_0}}\mathbf{P}_\rho$. Then, gathering all these estimates we have 
\begin{align}\label{estimateP}
\mathbf{P}_r\le C\biggl(\left(\frac{\rho}{r}\right)^{3-\frac{15}{2\tau_0}}
(\mathbf{A}_\rho \mathbf{B}_\rho)^{\frac{3}{4}}+\left(\frac{\rho}{r}\right)^{\frac{1}{2}-\frac{15}{2\tau_0}}\left[\rho^{\frac{15}{2\tau_0}-\frac32}\mathbf{B}_\rho^{\frac34}\|\va\|_{L_{t,x}^6(Q_\rho)}^{\frac32}+\mathbf{P}_\rho\right]\biggr).
\end{align}
Now we fix $0<\kappa<1$ such that $r=\kappa \rho$. Then, we define a new expression that will help us to set up the inductive argument
\begin{equation}\label{theta}
\mathbf{\Theta}_r(t,x)=\mathbf{A}_r(t,s)+\mathbf{H}_r(t,s)+\left(\kappa^{\frac{15}{\tau_0}-\frac{15}{2}}\mathbf{P}_r(t,x)\right)^{\frac{4}{3}}.
\end{equation}
We will see how to obtain from (\ref{estimateA}) and (\ref{estimateP}) a
recursive equation in terms of $\mathbf{\Theta}_r$ from which we will deduce (\ref{iterative}). Indeed,
we have the following lemma.
\begin{lem}\label{thetaire}
For all $(t,x)\in Q_{2R_1}(t_0,x_0)$, for all $0<r<\frac{\rho}{2}$ and for all $\rho$ small enough we have
the inequality
$$\mathbf{\Theta}_r(t,x)\le \frac{1}{2}\mathbf{\Theta}_\rho(t,x)+\epsilon,$$
where $\epsilon$ is a small constant that depends on the information available on the forces $\vf$, $\vg$ and the perturbation $\va$.
\end{lem}
\textbf{Proof.} We will use the estimates (\ref{estimateA}) and (\ref{estimateP}) obtained previously. Indeed, introducing the quantity $\kappa=\frac{r}{\rho}$ we easily obtain:
\begin{eqnarray}
&&\mathbf{\Theta}_r=\mathbf{A}_r+\mathbf{H}_r+\left(\kappa^{\frac{15}{\tau_0}-\frac{15}{2}}\mathbf{P}_r\right)^{\frac{4}{3}}\leq C\Bigg(\underbrace{\kappa^{\frac{10}{\tau_0}}\mathbf{A}_\rho+\kappa^{\frac{10}{\tau_0}-4}\mathbf{A}_\rho \mathbf{B}_\rho^{\frac{1}{2}}}_{(1)}+\underbrace{\kappa^{\frac{10}{\tau_0}-4}\mathbf{P}_\rho^{\frac{2}{3}}\mathbf{A}_\rho^{\frac{1}{2}}}_{(2)}\notag\\
&&+\underbrace{\left(\kappa^{\frac{10}{\tau_0}-4}+\kappa^{\frac{10}{\tau_0}-3}\right)\rho^{\frac{5}{\tau_0}-\frac56}\mathbf{A}_\rho^{\frac12}\mathbf{B}_\rho^{\frac12}\|\va\|_{L^6_{t,x}(\Omega)}}_{(3)}+\underbrace{\kappa^{\frac{10}{\tau_0}-3}\rho^{\frac{5}{\tau_0}}\mathbf{A}_\rho^{\frac12}\mathbf{B}_\rho^{\frac12}}_{(4)}+\underbrace{\kappa^{\frac{10}{\tau_0}-3}\rho^{2+\frac{5}{\tau_0}-\frac{4}{\tau_c}}\mathbf{D}_\rho\mathbf{A}_\rho^{\frac{1}{2}}}_{(5)}\notag\\
&&+\underbrace{\kappa^{\frac{10}{\tau_0}-4} \rho^{\frac{10}{\tau_0}-1}\mathbf{B}_{\rho}^{\frac{1}{2}}\|\vw\|_{L^\infty_{t,x}(\Omega)}}_{(6)}\Bigg)+C\underbrace{\biggl(\kappa^{\frac{45}{2\tau_0}-\frac{21}{2}}(\mathbf{A}_\rho \mathbf{B}_\rho)^{\frac{3}{4}}+\kappa^{\frac{45}{2\tau_0}-8}\left[\rho^{\frac{15}{2\tau_0}-\frac32}\mathbf{B}_\rho^{\frac34}\|\va\|_{L_{t,x}^6(Q_\rho)}^{\frac32}+\mathbf{P}_\rho\right]\biggr)^{\frac43}}_{(7)}.\label{EstimationTheta}
\end{eqnarray}
We will now study each one of the previous terms. 
\begin{itemize}
\item The first term above can be easily treated as we obviously have $\mathbf{A}_\rho\leq \mathbf{\Theta}_\rho$, thus we write
$$\kappa^{\frac{10}{\tau_0}}\mathbf{A}_\rho+\kappa^{\frac{10}{\tau_0}-4}\mathbf{A}_\rho \mathbf{B}_\rho^{\frac{1}{2}}\leq \kappa^{\frac{10}{\tau_0}}\mathbf{\Theta}_\rho+\kappa^{\frac{10}{\tau_0}-4}
\mathbf{ \Theta_\rho} \mathbf{B}_\rho^{\frac{1}{2}}.$$
\item For the term $(2)$ of (\ref{EstimationTheta}) we write, by the Young inequalities
\begin{eqnarray*}
\kappa^{\frac{10}{\tau_{0}}-4} \Poq^{\frac{2}{3}} \Ao^{\frac{1}{2}}&=&\kappa^{\frac{10}{\tau_{0}}-4} \left(\kappa^{5(\frac{1}{\tau_{0}} -\frac{1}{2})} \Poq^{\frac{2}{3}} \times\kappa^{5(\frac{1}{2} - \frac{1}{\tau_{0}} )} \Ao^{\frac{1}{2}}\right)\leq\kappa^{\frac{10}{\tau_{0}}-4} \left( \kappa^{10(\frac{1}{2} - \frac{1}{\tau_{0}} )} \Ao + \kappa^{10(\frac{1}{\tau_{0}} - \frac{1}{2} )} \Poq^{\frac{4}{3}} \right)\notag \\
&\leq & \kappa \left(  \Ao + \left(\kappa^{\frac{15}{\tau_{0}} -\frac{15}{2})}\Poq \right)^{\frac{4}{3}} \right) \leq \kappa  \TTo.
\end{eqnarray*}
\item For the term $(3)$ of (\ref{EstimationTheta}),  we obtain by the Young inequalities (and noting that we have $\kappa^{\frac{10}{\tau_0}-3}\leq \kappa^{\frac{10}{\tau_0}-4}$ since $0<\kappa<1$):
\begin{eqnarray*}
\left(\kappa^{\frac{10}{\tau_0}-4}+\kappa^{\frac{10}{\tau_0}-3}\right)\rho^{\frac{5}{\tau_0}-\frac56}\mathbf{A}_\rho^{\frac12}\mathbf{B}_\rho^{\frac12}\|\va\|_{L^6_{t,x}(\Omega)}&\leq &C\kappa^{\frac{10}{\tau_0}-4}\rho^{\frac{5}{\tau_0}-\frac56}(\kappa^4\mathbf{A}_\rho)^{\frac12}(\kappa^{-4}\mathbf{B}_\rho)^{\frac12}\|\va\|_{L^6_{t,x}(\Omega)}\\
&\leq &C\kappa^{\frac{10}{\tau_0}}\rho^{\frac{5}{\tau_0}-\frac56}\mathbf{\Theta}_\rho\|\va\|_{L^6_{t,x}(\Omega)}+C\kappa^{\frac{10}{\tau_0}-8}\rho^{\frac{5}{\tau_0}-\frac{5}{6}}\mathbf{B}_\rho\|\va\|_{L^6_{t,x}(\Omega)}.
\end{eqnarray*}
\item The term $(4)$ of (\ref{EstimationTheta}) is treated as follows:
$$\kappa^{\frac{10}{\tau_0}-3}\rho^{\frac{5}{\tau_0}}\mathbf{A}_\rho^{\frac12}\mathbf{B}_\rho^{\frac12}\leq \kappa^{\frac{10}{\tau_0}-3}\rho^{\frac{5}{\tau_0}}(\kappa^{\frac32}\mathbf{A}_\rho)^{\frac12}(\kappa^{-\frac32}\mathbf{B}_\rho)^{\frac12}\leq \kappa^{\frac{10}{\tau_0}}\rho^{\frac{5}{\tau_0}}\mathbf{\Theta}_\rho+\kappa^{\frac{10}{\tau_0}-6}\rho^{\frac{5}{\tau_0}}\mathbf{B}_\rho.$$
\item For the term $(5)$ of (\ref{EstimationTheta}) we simple write:
$$\kappa^{\frac{10}{\tau_0}-3}\rho^{2+\frac{5}{\tau_0}-\frac{4}{\tau_c}}\mathbf{D}_\rho\mathbf{A}_\rho^{\frac{1}{2}}\leq  \kappa^{\frac{10}{\tau_0}-3}\rho^{2+\frac{5}{\tau_0}-\frac{5}{\tau_c}}(\mathbf{D}_\rho^2+\mathbf{A}_\rho)\leq  \kappa^{\frac{10}{\tau_0}-3}\rho^{2+\frac{5}{\tau_0}-\frac{5}{\tau_c}}(\mathbf{D}_\rho^2+\mathbf{\Theta}_\rho).$$
\item The term $(6)$ of (\ref{EstimationTheta}) needs no particular treatment.
\item For the last term of (\ref{EstimationTheta}), using the fact that $\left(\kappa^{\frac{15}{\tau_0}-\frac{15}{2}}\mathbf{P}_\rho\right)^{\frac{4}{3}}\leq \mathbf{\Theta_\rho}$ by the definition of $\mathbf{\Theta_\rho}$ given in (\ref{theta}), we write:
\begin{eqnarray*}
\biggl(\kappa^{\frac{45}{2\tau_0}-\frac{21}{2}}(\mathbf{A}_\rho \mathbf{B}_\rho)^{\frac{3}{4}}+\kappa^{\frac{45}{2\tau_0}-8}\left[\rho^{\frac{15}{2\tau_0}-\frac32}\mathbf{B}_\rho^{\frac34}\|\va\|_{L_{t,x}^6(Q_\rho)}^{\frac32}+\mathbf{P}_\rho\right]\biggr)^{\frac43}&\leq &C\bigg(\kappa^{\frac{30}{\tau_0}-14}\mathbf{\Theta}_\rho \mathbf{B}_\rho+\\
&&+\kappa^{\frac{40}{\tau_0}-10}\mathbf{B}_\rho\|\va\|_{L_{t,x}^6(Q_\rho)}^{2}+\kappa^{\frac{10}{\tau_0}-\frac{2}{3}} \mathbf{\Theta}_\rho\bigg),
\end{eqnarray*}
\end{itemize}
Gathering all these estimates we observe that from (\ref{EstimationTheta}) we can write
\begin{align}
\mathbf{\Theta}_r &\le C\Biggl(\kappa^{\frac{10}{\tau_0}}+\kappa^{\frac{10}{\tau_0}-4}\mathbf{B}_\rho^{\frac{1}{2}}+\kappa+\kappa^{\frac{10}{\tau_0}}\rho^{\frac{5}{\tau_0}-\frac56}\|\va\|_{L^6_{t,x}(\Omega)}+\kappa^{\frac{10}{\tau_0}-3}\rho^{2+\frac{5}{\tau_0}-\frac{5}{\tau_c}}+\kappa^{\frac{30}{\tau_0}-14}\mathbf{B}_\rho+\kappa^{\frac{10}{\tau_0}-\frac{2}{3}}\Biggr)\mathbf{\Theta}_\rho\label{thetaite1}\\
&+C\bigg(\kappa^{\frac{10}{\tau_0}-8}\rho^{\frac{5}{\tau_0}-\frac{5}{6}}\mathbf{B}_\rho\|\va\|_{L^6_{t,x}(\Omega)}+\kappa^{\frac{10}{\tau_0}-6}\rho^{\frac{5}{\tau_0}}\mathbf{B}_\rho+\kappa^{\frac{10}{\tau_0}-3}\rho^{2+\frac{5}{\tau_0}-\frac{5}{\tau_c}}\mathbf{D}_\rho^2\notag\\
&+\kappa^{\frac{10}{\tau_0}-4} \rho^{\frac{10}{\tau_0}-1}\mathbf{B}_{\rho}^{\frac{1}{2}}\|\vw\|_{L^\infty_{t,x}(\Omega)}+\kappa^{\frac{40}{\tau_0}-10}\mathbf{B}_\rho\|\va\|_{L_{t,x}^6(Q_\rho)}^{2}\bigg).\label{thetaite2}
\end{align}
We claim now that we have, for the term (\ref{thetaite1}) above the following control
$$C\bigg(\kappa^{\frac{10}{\tau_0}}+\kappa^{\frac{10}{\tau_0}-4}\mathbf{B}_\rho^{\frac{1}{2}}+\kappa+\kappa^{\frac{10}{\tau_0}}\rho^{\frac{5}{\tau_0}-\frac56}\|\va\|_{L^6_{t,x}(\Omega)}+\kappa^{\frac{10}{\tau_0}-3}\rho^{2+\frac{5}{\tau_0}-\frac{5}{\tau_c}}+\kappa^{\frac{30}{\tau_0}-14}\mathbf{B}_\rho+\kappa^{\frac{10}{\tau_0}-\frac{2}{3}}\bigg)\leq \frac12.$$
Indeed, we recall that $\kappa=\frac{r}{\rho}<1$ is a fixed small parameter and that $0<\rho<1$ is also a small parameter. Moreover we recall that due to the hypothesis (\ref{grad}), we have $\underset{\rho \to 0}{\limsup}\mathbf{B}_\rho \le \widetilde{\epsilon}$ where $\widetilde{\epsilon}>0$ is also very small. Then all the terms of the form $\kappa^{a}$, $\kappa^a \rho^b$ with $a,b>0$ and $\kappa^{-c}\mathbf{B}_\rho$ or $\kappa^{-c}\mathbf{B}_\rho^\frac{1}{2}$ with $c>0$ can be made very small. Note that the size of the perturbation term, reflected in the quantity $\kappa^{\frac{10}{\tau_0}}\rho^{\frac{5}{\tau_0}-\frac56}\|\va\|_{L^6_{t,x}(\Omega)}$ can be easily absorbed as $\rho$ can be very small (we have $\frac{5}{\tau_0}-\frac56>0$ as $\frac{5}{1-\alpha}<\tau_0<\frac{11}{2}$). We remark that since $\rho$ is small, we have that the term $\kappa^{\frac{10}{\tau_0}-3}\rho^{2+\frac{5}{\tau_0}-\frac{5}{\tau_c}}$ can also made small as we have $2+\frac{5}{\tau_0}-\frac{5}{\tau_c}>0$. Finally note that $\frac{10}{\tau_0}-\frac{2}{3}>0$ since we have the upper bound $\frac{11}{2}>\tau_0>\frac{5}{1-\alpha}$.\\ 

By essentially the same arguments if $\rho>0$ is small enough, we have the following control for (\ref{thetaite2}):
\begin{eqnarray*}
C\big(\kappa^{\frac{10}{\tau_0}-8}\rho^{\frac{5}{\tau_0}-\frac{5}{6}}\mathbf{B}_\rho\|\va\|_{L^6_{t,x}(\Omega)}+\kappa^{\frac{10}{\tau_0}-6}\rho^{\frac{5}{\tau_0}}\mathbf{B}_\rho+\kappa^{\frac{10}{\tau_0}-3}\rho^{2+\frac{5}{\tau_0}-\frac{5}{\tau_c}}\mathbf{D}_\rho^2\\
+\kappa^{\frac{10}{\tau_0}-4} \rho^{\frac{10}{\tau_0}-1}\mathbf{B}_{\rho}^{\frac{1}{2}}\|\vw\|_{L^\infty_{t,x}(\Omega)}+\kappa^{\frac{40}{\tau_0}-10}\mathbf{B}_\rho\|\va\|_{L_{t,x}^6(Q_\rho)}^{2}\big)<\epsilon,
\end{eqnarray*}
where $\epsilon>0$ can be made small (remark that the quantity $\|\vw\|_{L^\infty_{t,x}(\Omega)}$ can easily be absorbed for $\rho$ small enough as we have $\frac{10}{\tau_0}-1>0$ since $\tau_0<\frac{11}{2}$. Note that the condition $\vw\in L^p_{t}L^q_x(\Omega)$ with $\frac{10}{\tau_0}-1-\frac{2}{p}-\frac{3}{q}>0$ stated in Remark \ref{Remarque_HypothesesVW} will give a similar result. See also Remark \ref{Remarque_HypothesesVW1} for this particular point). With these last observations, then from the inequality (\ref{thetaite1})-(\ref{thetaite2}), we obtain $\mathbf{\Theta}_r\le \frac{1}{2}\mathbf{\Theta}_\rho+\epsilon$ which is the conclusion of the Lemma \ref{thetaire}. \hfill $\blacksquare$\\
 
With this lemma at hand, we continue the proof of the Proposition  \ref{PropositionFirstMorreySpace}. Indeed, for any radius $\rho$ such that $0<\rho<{\bf R}<1$ (and since we have $Q_{{\bf R}}(t_0,x_0)\subset \Omega$) by the set of hypotheses (\ref{HypothesesLocal1}) we have the bounds $\|\vu\|_{L_t^\infty L_x^2(Q_\rho(t_0,x_0))}\le \|\vu\|_{L_t^\infty L_x^2(\Omega)}<+\infty$, $\|\grad \otimes\vu\|_{L_{t,x}^2(Q_\rho(t_0,x_0))}\le \|\grad\otimes\vu\|_{L_{t,x}^2(\Omega)}<+\infty$ (and the same estimates for $\vb$ and $\vw$) and $ \|p\|_{L_{t,x}^{\frac{3}{2}}(Q_\rho(t_0,x_0))}\le \|p\|_{L_{t,x}^{\frac{3}{2}}(\Omega)}<+\infty$. Then, by the Definition \ref{Def_Invariants}, we have the uniform bounds $\underset{0<\rho<{\bf R}}{\sup}\biggl\{ \rho \mathcal{A}_\rho,\rho \alpha_\rho, \rho \mathcal{B}_\rho,\rho \beta_\rho, \rho\mathcal{C}_\rho,\rho  \gamma_\rho, \rho\mathcal{W}_\rho, \rho^2 \mathcal{H}_\rho, \rho^2\mathcal{P}_\rho\biggr\}<+\infty$ from which we can deduce by the definition of the quantities $\mathbf{A}_\rho(t_0,x_0)$, $\mathbf{H}_\rho(t_0,x_0)$
and $\mathbf{P}_\rho (t_0,x_0)$ given in (\ref{defA}) and (\ref{defB}), the uniform bounds
\begin{equation}\label{morrey5}
\begin{split}
\sup_{0<\rho<{\bf R}}\rho^{3-\frac{10}{\tau_0}}\mathbf{A}_\rho(t_0,x_0)<+\infty, \quad &\sup_{0<\rho<{\bf R}}\rho^{2-\frac{10}{\tau_0}}\mathbf{H}_\rho(t_0,x_0)<+\infty,\quad\\
\text{and} \quad &\sup_{0<\rho<{\bf R}}\rho^{5-\frac{3}{2}(1+\frac{5}{\tau_0})}\mathbf{P}_\rho (t_0,x_0)<+\infty.
\end{split}
\end{equation}
Note now, that there exists a $0<\kappa<\frac{1}{2}$ and a fixed $0<\rho_0<{\bf R}$ small such that, by (\ref{morrey5}), the quantities $\mathbf{A}_{\rho_0}$, $\mathbf{H}_{\rho_0}$ and $\mathbf{P}_{\rho_0}$ are bounded: indeed, recall that we have $\tau_0>\frac{5}{1-\alpha}>5$ (where $0<\alpha<\frac{1}{12}$) and this implies that all the powers of $\rho$ in the expression above are positive. As a consequence of this fact, by (\ref{theta}) the quantity $\mathbf{\Theta}_{\rho_0}$ is itself bounded. Remark also that, if $\rho_0$
is small enough, then the inequality (\ref{thetaite1}) holds true and we can write $\mathbf{\Theta}_{\kappa\rho_0}(t_0,x_0)\le \frac{1}{2}\mathbf{\Theta}_{\rho_0}(t_0,x_0)+\epsilon$. We can iterate this process and we obtain for all $n>1$,
$$\mathbf{\Theta}_{\kappa^n\rho_0}(t_0,x_0)\le \frac{1}{2^n}\mathbf{\Theta}_{\rho_0}(t_0,x_0)+\epsilon\sum_{j=0}^{n-1}2^{-j},$$
and therefore there exists $N\ge 1$ such that for all $n\ge N$ we have
$\mathbf{\Theta}_{\kappa^n\rho_0}(t_0,x_0)\le 4\epsilon$ from which we obtain (using the definition of $\mathbf{\Theta}_\rho$ given in (\ref{theta})) that
$$\mathbf{A}_{\kappa^N\rho_0}(t_0,x_0)\le \frac{1}{8}C,\quad \mathbf{H}_{\kappa^N\rho_0}(t_0,x_0)\le \frac{1}{8}C\quad \text{and}\quad
\mathbf{P}_{\kappa^N\rho_0}(t_0,x_0)\le \frac{1}{32}C.$$
This information is centered at the point $(t_0,x_0)$, in order to treat
the uncentered bound, we can let $\frac{1}{2}\kappa^N\rho_0$ to be the radius $R_1$
we want to find, thus for all points $(t,x)\in Q_{R_1}$ we have that
$Q_{R_1}\subset Q_{2R_1}(t_0,x_0)$, which implies 
\begin{equation*}
     \mathbf{A}_{R_1}(t,x)\le 2^{3-\frac{10}{\tau_0}} \mathbf{A}_{2R_1}(t_0,x_0)
    \le 8 \mathbf{A}_{2R_1}(t_0,x_0)\le 8 \mathbf{A}_{\kappa^N \rho }(t_0,x_0)<C,
\end{equation*}
\begin{equation*}
\mathbf{H}_{R_1}(t,x)\le 2^{-\frac{10}{\tau_0}}\mathbf{H}_{2R_1}(t_0,x_0)
    \le  \mathbf{H}_{2R_1}(t_0,x_0)\le  \mathbf{H}_{\kappa^N \rho }(t_0,x_0)<C,
\end{equation*}
and $\mathbf{P}_{R_1}(t,x)\le 2^{5-\frac{3}{2}(1+\frac{5}{\tau_0})}\mathbf{P}_{2R_1}(t_0,x_0)\le 32\mathbf{P}_{2R_1}(t_0,x_0)\le 8 \mathbf{P}_{\kappa^N \rho }(t_0,x_0)<C$. Having obtained these bounds, by the definition of $\mathbf{\Theta}_{R_1}$, we thus get $\mathbf{\Theta}_{R_1}(t,x)\le C$. Applying Lemma \ref{thetaire} and iterating once more, we find that the same will be true for $\kappa R_1$ and then, for all $\kappa^n R_1$, $n\in \mathbb{N}$. Since by definition we have  $\mathbf{A}_{\kappa^n R_1}(t,x)\le  \mathbf{\Theta}_{\kappa^n R_1}(t,x)$ we have finally obtained the estimate $\mathbf{A}_{\kappa^n R_1}(t,x)\le C$ and the inequality (\ref{iterative}) is proven which implies Proposition \ref{PropositionFirstMorreySpace}.\hfill$\blacksquare$
\begin{coro}\label{corolarioMorrey}
Under the hypotheses of the Proposition \ref{PropositionFirstMorreySpace}, we have 
\begin{equation}\label{ConlusioncorolarioMorrey}
\mathds{1}_{Q_{R_1}}\grad\otimes \vu\in \mathcal{M}^{2,\tau_1}_{t,x},\quad\mathds{1}_{Q_{R_1}}\grad\otimes \vb\in \mathcal{M}^{2,\tau_1}_{t,x}, \quad \mathds{1}_{Q_{R_1}}\grad\otimes \vw\in \mathcal{M}^{2,\tau_1}_{t,x}\quad\text{and}\quad\mathds{1}_{Q_{R_1}}div(\vw)\in \mathcal{M}^{2,\tau_1}_{t,x},
\end{equation}
where $\frac{1}{\tau_1}=\frac{1}{\tau_0}+\frac{1}{5}$.
\end{coro}
\textbf{Proof.} Indeed, from the general notation given in Definition \ref{Def_Invariants}, we have
$$\frac{1}{r}\bigg(\iint _{Q_r}|\grad \otimes \vu|^2dyds+\iint _{Q_r}|\grad \otimes \vb|^2dyds+\iint _{Q_r}|\grad \otimes \vw|^2dyds+\iint _{Q_r}|div(\vw)|^2dyds\bigg)=(\alpha_r+\beta_r+\gamma_r+\mathcal{W}_r),$$
and by the definition of $\mathbf{A}_r$ given in (\ref{defA}) we obtain
$$\frac{1}{r}\bigg(\iint _{Q_r}|\grad \otimes \vu|^2dyds+\iint _{Q_r}|\grad \otimes \vb|^2dyds+\iint _{Q_r}|\grad \otimes \vw|^2dyds+\iint _{Q_r}|div(\vw)|^2dyds\bigg)\leq r^{2(1-\frac{5}{\tau_0})}\mathbf{A}_r.$$
But since the quantity $\mathbf{A}_r$ is bounded for $0<r<1$ small enough (by the estimate (\ref{iterative})), we can thus write:
$$\iint _{Q_r}|\grad \otimes \vu|^2dyds+\iint _{Q_r}|\grad \otimes \vb|^2dyds+\iint _{Q_r}|\grad \otimes \vw|^2dyds+\iint _{Q_r}|div(\vw)|^2dyds\leq Cr^{3-\frac{10}{\tau_0}}=Cr^{5(1-\frac{2}{\tau_1})},$$
since we have $\frac{1}{\tau_1}=\frac{1}{\tau_0}+\frac{1}{5}$ and by the definition of Morrey spaces fiven in (\ref{DefMorreyparabolico}) this condition expresses the fact that each term of the left-hand side above belongs locally to the Morrey space $\mathcal{M}^{2,\tau_1}_{t,x}$.\hfill$\blacksquare$
\begin{rema}
From the Corollary \ref{corolarioMorrey}, we can easily deduce that
\begin{equation*}
\mathds{1}_{Q_{R_1}}\rot \vu\in \mathcal{M}^{2,\tau_1}_{t,x},\quad\mathds{1}_{Q_{R_1}}\rot \vb\in \mathcal{M}^{2,\tau_1}_{t,x} \quad \text{and}\quad    \mathds{1}_{Q_{R_1}}\rot\vw\in \mathcal{M}^{2,\tau_1}_{t,x}.
\end{equation*}
\end{rema}
We have proven the points 1), 2)  of the hypotheses of Proposition \ref{HolderRegularityproposition} (recall that the point 6) is given for free, due to the hypotheses on the external forces) and we still need to prove the points 3), 4) and 5). In order to achieve this task, we will need different arguments that are displayed in the next section.
%%%%%%%%%%%%%%%%%%%%%%%%%%%%%%%%%%%%%%
\section{More estimates}\label{Secc_MasEstimates}
Let $0<\mathfrak{a}<5$ be a parameter, we define the parabolic Riesz potential $\mathcal{L}_{\mathfrak{a}}$ of a locally integrable function
$\vec f: \mathbb{R}\times\mathbb{R}^3\longrightarrow \mathbb{R}^3$ by
\begin{equation}\label{DefinitionPotentielRiesz6}
\mathcal{L}_{\mathfrak{a}}(\vf)(t,x)=\int_{\mathbb{R}}\int_{\mathbb{R}^3}
\frac{1}{(|t-s|^{\frac{1}{2}}+|x-y|)^{5-\mathfrak{a}}}\vec{f}(s,y)dyds.
\end{equation}
Then, we have the following property
%%%%%%%%%%%%%%%%%%%%%%%%%%%%%%%%%%%%%%
\begin{lem}[Adams-Hedberg inequality]\label{Lemme_Hed}
If $0<\mathfrak{a}<\frac{5}{q}$, $1<p\le q<+\infty$ and $\vf\in \M^{p,q}$,
then for $\lambda=1-\frac{\mathfrak{a}q}{5}$ we have the following boundedness property in Morrey spaces:
\begin{equation*}
\|\mathcal{L}_{\mathfrak{a}}(\vf)\|_{\M^{\frac{p}{\lambda},\frac{q}{\lambda}}}\le\|\vf\|_{\M^{p,q}}. 
\end{equation*}
\end{lem}
See a proof of this fact in the book \cite[Corollary 5.1]{PGLR1}. We will use this result in the next result to obtain the hypothesis 4) of the Proposition \ref{HolderRegularityproposition}. 
%%%%%%%%%%%%%%%%%%%%%%%%%%%%%%%%%%%%%%
\begin{propo}\label{MorreyPropo2}
Let $(\vu,\vb,\vw, p)$ be a suitable solution of the system (\ref{EquationMMP}) over $\Omega$. Then for some radius $R_2>0$ such that $R_2<R_1$, we have (with $\frac{1}{\delta}+\frac{1}{\tau_0}<\frac{1-\alpha}{5}$):
\begin{equation*}
\mathds{1}_{Q_{R_2}}\vu\in \M^{3,\delta},\quad\mathds{1}_{Q_{R_2}}\vb\in \M^{3,\delta}, \quad \text{and}\quad\mathds{1}_{Q_{R_2}}\vw\in \M^{3,\delta},
\end{equation*}
\end{propo}
%%%%%%%%%%%%%%%%%%%%%%%%%%%%%%%%%%%%%%
\textbf{Proof.} For a point $(t_0,x_0)\subset \Omega$ we consider the radii
$$0<R_2<\Bar{R}<\tilde R<R_1<4 {\bf R},$$  
(recall (\ref{DefinitionRayons})) and the associated parabolic balls $Q_{R_2}\subset Q_{\bar R}\subset Q_{\tilde R}\subset Q_{R_1}\subset Q_{4 {\bf R}}$. Consider now $\bar \phi, \bar \psi: \mathbb{R}\times\mathbb{R}^3\longrightarrow \mathbb{R}$ two  non-negative functions such that $\phi, \psi\in \mathcal{C}_0^{\infty}
(\mathbb{R}\times \mathbb{R}^3)$ and such that
\begin{equation}\label{ProprieteLocalisation6}
\bar \phi \equiv 1\;\; \text{over}\; \; Q_{R_2},\; \; supp(\bar \phi)\subset Q_{\bar R} \quad \mbox{and}\quad\bar \psi \equiv 1\;\; \text{over}\; \; Q_{\tilde R},\; \; supp(\bar \psi)\subset Q_{ R_1}.
\end{equation}
Since $R_1<{\bf R}<t_0$, we have $\bar \phi (0,\cdot)=\bar \psi (0,\cdot)=0$ and we also have the identity $\bar \phi \bar \psi=\bar \phi$.\\

With the help of these localizing functions we will study the evolution of the variable $   \vec{\mathcal{V}}=\bar \phi(\vu+\vb+\vw)$ and we obtain the system
\begin{equation}\label{equationV}
\begin{cases}
\partial_t \vec {\mathcal{V}}=\Delta \vec {\mathcal{V}}+\vec {\mathcal{N}},\\[3mm]
\vec {\mathcal{V}}(0,x)=0,
\end{cases}
\end{equation}
where, following the same computations of (\ref{SystemeLocalise}) we have
\begin{eqnarray}
\vec {\mathcal{N}}&=&(\partial_t\bar\phi- \Delta \bar\phi)(\vu+\vb+\vw)-2\sum_{i=1}^3 (\partial_i \bar\phi)(\partial_i (\vu+\vb+\vw))-\bar\phi\bigg((\vb \cdot \grad)\vu+(\vu \cdot \grad)\vb\bigg)\notag\\
&&-2\bar\phi\grad p+\bar\phi(\rot\vw) +\bar\phi(\vf+\vg)- \bar\phi\left(div((\vu+\vb) \otimes \va+\va \otimes (\vu+\vb))\right)\label{DefinitionVariableN}\\
&&+ \bar\phi \left(\grad div(\vw)-\vw-\frac{1}{2}((\vu+\vb)\cdot \grad)\vw+\frac{1}{4}\rot(\vu+\vb)\right).\notag
\end{eqnarray}
Now we will perform some computations over the term $\bar\phi\grad p$ that contains the pressure. Indeed, as we have the identity $p=\bar \psi p$ over $Q_{\tilde R}$, then over the smaller ball $Q_{R_2}$  (recalling that $\bar \psi=1$ over $Q_{R_2}$ by (\ref{ProprieteLocalisation6}) since $Q_{R_2}\subset Q_{\tilde R}$), we can write $-\Delta(\bar \psi p)=-\bar \psi \Delta p+(\Delta \bar \psi)p-2\displaystyle{\sum_{i=1}^3}\partial_i((\partial_i\bar \psi)p)$
from which we deduce the identity 
\begin{equation}\label{ExpressionPression6}
\bar \phi \grad p=\bigg(\bar \phi\frac{\grad (-\bar \psi \Delta p)}{(-\Delta)}\bigg)+\bar \phi\frac{\grad ((\Delta \bar \psi)p)}{(-\Delta)}-2\sum_{i=1}^3\bar \phi\frac{\grad (\partial_i((\partial_i\bar \psi)p))}{(-\Delta)}.
\end{equation}
We recall now that by (\ref{FormulePressionIntro}) we have $2\Delta p = -\displaystyle{\sum_{i,j=1}^3}\partial_i\partial_j\left(u_i b_j+(u_i+b_i)a_j+a_i(u_j+b_j)\right)$ and thus, the first term of the right-hand side of the previous formula can be written in the following manner:
\begin{equation}\label{FormulaHipotesis5}
\bar \phi\frac{\grad (-\bar \psi \Delta p)}{(-\Delta)}= \bar \phi\frac{\grad}{2(-\Delta)}\left( \sum_{i,j=1}^3\bar \psi\left[\partial_i\partial_j \big(u_i b_j+(u_i+b_i)a_j+a_i(u_j+b_j)\big)\right]\right),
\end{equation}
and introducing the function $\bar \psi$ inside the derivatives we obtain
\begin{align}
\bar \phi\frac{\grad (-\bar \psi \Delta p)}{(-\Delta)}&= \sum_{i,j=1}^3 \bar \phi\frac{\grad}{2(-\Delta)}\biggl(\partial_i\partial_j(\bar \psi u_i b_j)-\partial_i((\partial_j \bar \psi)u_i b_j)-\partial_j((\partial_i \bar \psi)u_i b_j)+(\partial_i \partial_j \bar \psi)(u_i b_j)\label{FormulaHipotesis51}\\
&\;\;+ \partial_i\partial_j(\bar \psi (u_i+b_i) a_j)-\partial_i((\partial_j \bar \psi)(u_i +b_i)a_j)-\partial_j((\partial_i \bar \psi)(u_i+b_i) a_j)+(\partial_i \partial_j \bar \psi)((u_i+b_i) a_j)\notag\\
&\;\;+ \partial_i\partial_j(\bar \psi  a_i(u_j+b_j))-\partial_i((\partial_j \bar \psi)a_i(u_j +b_j))-\partial_j((\partial_i \bar \psi) a_i(u_j+b_j))+(\partial_i \partial_j \bar \psi)( a_i(u_j+b_j))\biggr).\notag
\end{align}
Now for the first terms of each line above we use the identities (recall that $\bar \phi \bar \psi=\bar \phi$):
\begin{align*}
\bar \phi\frac{\grad}{(-\Delta)}\partial_i\partial_j(\bar \psi u_i b_j)&=\left[\bar \phi, \frac{\grad\partial_i\partial_j}{(-\Delta)}\right](\bar \psi u_i b_j)+ \frac{\grad\partial_i\partial_j}{(-\Delta)}(\bar \phi u_i b_j)\\
\bar \phi\frac{\grad}{(-\Delta)}\partial_i\partial_j(\bar \psi(u_i+b_i) a_j)&=\left[\bar \phi, \frac{\grad\partial_i\partial_j}{(-\Delta)}\right](\bar \psi(u_i+b_i) a_j)+ \frac{\grad\partial_i\partial_j}{(-\Delta)}(\bar \phi(u_i+b_i) a_j)\\
\bar \phi\frac{\grad}{(-\Delta)}\partial_i\partial_j(\bar \psi a_i(u_j+b_j))&=\left[\bar \phi, \frac{\grad\partial_i\partial_j}{(-\Delta)}\right](\bar \psi a_i (u_j+b_j))+\frac{\grad\partial_i\partial_j}{(-\Delta)}(\bar \phi  (u_j+b_j)a_i),
\end{align*}
and with this lengthy and tedious formulation for the first term of (\ref{ExpressionPression6}), we come back to the term $\vec{\mathcal{N}}$ given in (\ref{DefinitionVariableN}) to obtain
\begin{align}
&\vec {\mathcal{N}}= \underbrace{(\partial_t \bar \phi - \Delta \bar \phi)(\vu+\vb+\vw)}_{(1)}-2\sum_{i=1}^3 \underbrace{(\partial_i \bar \phi)(\partial_i (\vu+\vb+\vw))}_{(2)}-\underbrace{\bar \phi\bigg((\vb\cdot \grad)\vu+(\vu\cdot \grad)\vb\bigg)}_{(3)}\label{FormulePourN}\\
&-\bigg\lbrace\bigg(\sum_{i,j=1}^3\underbrace{\left[\bar \phi,\frac{\grad\partial_i\partial_j}{(-\Delta)}\right](\bar \psi u_i b_j)}_{(4)}+\sum_{i,j=1}^3 \underbrace{\frac{\grad\partial_i\partial_j}{(-\Delta)}(\bar \phi u_i b_j)}_{(5)}- \sum_{i,j=1}^3 \frac{\bar \phi\grad}{(-\Delta)}\big[\underbrace{\partial_i((\partial_j \bar \psi)u_i b_j)}_{(6)}+\underbrace{\partial_j((\partial_i \bar \psi)u_i b_j)}_{(7)}\notag\\
&-\underbrace{(\partial_i \partial_j \bar \psi)(u_i b_j)}_{(8)}\big]\bigg)+\bigg(\underbrace{\left[\bar \phi, \frac{\grad\partial_i\partial_j}{(-\Delta)}\right](\bar \psi(u_i+b_i) a_j)}_{(9)}+\underbrace{\frac{\grad\partial_i\partial_j}{(-\Delta)}(\bar \phi(u_i+b_i) a_j)}_{(10)}- \sum_{i,j=1}^3 \frac{\bar \phi\grad}{(-\Delta)}\big[\underbrace{\partial_i((\partial_j \bar \psi)(u_i +b_i)a_j)}_{(11)}\notag\\
&+\underbrace{\partial_j((\partial_i \bar \psi)(u_i+b_i) a_j)}_{(12)}-\underbrace{(\partial_i \partial_j \bar \psi)((u_i+b_i) a_j)}_{(13)}\big]\bigg)+\bigg(\underbrace{\left[\bar \phi, \frac{\grad\partial_i\partial_j}{(-\Delta)}\right](\bar \psi a_i (u_j+b_j))}_{(14)}+\underbrace{\frac{\grad\partial_i\partial_j}{(-\Delta)}(\bar \phi  (u_j+b_j)a_i)}_{(15)}\notag\\
&- \sum_{i,j=1}^3 \frac{\bar \phi\grad}{(-\Delta)}\big[\underbrace{\partial_i((\partial_j \bar \psi)a_i(u_j +b_j))}_{(16)}+\underbrace{\partial_j((\partial_i \bar \psi) a_i(u_j+b_j))}_{(17)}-\underbrace{(\partial_i \partial_j \bar \psi)( a_i(u_j+b_j))}_{(18)}\big]\bigg)+2\underbrace{\bar \phi\frac{\grad ((\Delta \bar \psi)p)}{(-\Delta)}}_{(19)}\notag\\
&-4\sum_{i=1}^3\underbrace{\bar \phi\frac{\grad (\partial_i((\partial_i\bar \psi)p))}{(-\Delta)}}_{(20)}\bigg\rbrace+\underbrace{\bar\phi(\rot\vw)}_{(21)}+\underbrace{\bar\phi(\vf+\vg)}_{(22)}-\underbrace{\bar\phi\left(div((\vu+\vb) \otimes \va+\va \otimes (\vu+\vb))\right)}_{(23)}\notag\\
&+ \underbrace{\bar\phi \left(\grad div(\vw)-\vw-\frac{1}{2}((\vu+\vb)\cdot \grad)\vw+\frac{1}{4}\rot(\vu+\vb)\right)}_{(24)}:=\sum_{i=1}^{24}\vec{\mathcal{N}}_k.\notag
\end{align}
With this expression of $\vec{\mathcal{N}}$, we obtain that the solution of the equation (\ref{equationV}) is given by
$$\vec{\mathcal{V}}=\int_0^{t}e^{(t-s)\Delta}\vec{\mathcal{N}}(s,\cdot)ds:=\sum_{k=1}^{24}\vec{\mathcal{V}}_k=\sum_{k=1}^{24}\int_0^{t}e^{(t-s)\Delta}\vec{\mathcal{N}}_k(s, \cdot)ds, $$
and we will study each term of this expression with the following lemma:
\begin{lem}\label{LemmaGainenMorreySigma}
Under the general hypotheses of the Theorem \ref{HolderRegularity_theorem}, for all $k=1,\dots, 24$ we have 
\begin{equation*}
\mathds{1}_{Q_{R_2}}\vec{\mathcal{V}}_k \in \M^{3,\sigma}
\end{equation*}
for some $\sigma$ close to $\tau_0$ such that $\tau_0<\sigma$.
\end{lem}
\textbf{Proof.} Fortunately many of the terms above share a similar structure as we have essentially the same information over the variables $\vu, \vb$ and $\vw$. Recall that we have proven so far the estimates (\ref{ConclusionPropositionFirstMorreySpace}) and (\ref{ConlusioncorolarioMorrey}). 
\begin{itemize}
\item For $\vec{\mathcal{V}}_1$, recalling that  $e^{(t-s)\Delta}\vec{\mathcal{N}}_1=\mathfrak{g}_{t-s}\ast \vec{\mathcal{N}}_1$ where $\mathfrak{g}_{t}$ is the usual 3D heat kernel, we have
$$|\mathds{1}_{Q_{R_2}}\vec{\mathcal{V}}_1(t,x)|=\left|\mathds{1}_{Q_{R_2}}\int_{0}^{t} \int_{\mathbb{R}^3}  \mathfrak{g}_{t-s}(x-y)[(\partial_t \bar\phi - \Delta \bar\phi)(\vu+\vb+\vw)](s,y)dyds\right|.$$
Thus, by the decay properties of the heat kernel as well as the properties of the test function $\bar\phi$ (see (\ref{ProprieteLocalisation6})), we have
$$|\mathds{1}_{Q_{R_2}}\vec{\mathcal{V}}_1(t,x)|\leq C\mathds{1}_{Q_{R_2}}\int_{\mathbb{R}} \int_{\mathbb{R}^3} \frac{1}{(|t-s|^{\frac{1}{2}}+|x-y|)^3}\left| \mathds{1}_{Q_{\bar R}}(\vu+\vb+\vw)(s,y)
\right| \,dy \,ds.$$
Now, recalling the definition of the Riesz potential given in (\ref{DefinitionPotentielRiesz6}) and since $Q_{R_2}\subset Q_{\bar R}$ we obtain the pointwise estimate$|\mathds{1}_{Q_{R_2}}\vec{\mathcal{V}}_1(t,x)|\leq C\mathds{1}_{Q_{\bar R}}\mathcal{L}_{2}(  |\mathds{1}_{Q_{\bar R}}(\vu+\vb+\vw)|)(t,x)$ and taking Morrey $\mathcal{M}_{t,x}^{3, \sigma}$ norm we obtain
$$\|\mathds{1}_{Q_{R_2}}\vec{\mathcal{V}}_1(t,x)\|_{\mathcal{M}_{t,x}^{3, \sigma}}\leq C\|\mathds{1}_{Q_{\bar R}}\mathcal{L}_{2}(  |\mathds{1}_{Q_{\bar R}}  (\vu + \vb+\vw)|)\|_{\mathcal{M}_{t,x}^{3, \sigma}}.$$
Now, for some $2<q<\frac{5}{2}$ we set $\lambda=1-\frac{2q}{5}$ and we define $3=\frac{a}{\lambda}$ and $\sigma< 10 <\frac{q}{\lambda}$ (remark that $a\leq q$). Thus, by Lemma \ref{lemma_locindi} and by Lemma \ref{Lemme_Hed} we can write:
\begin{eqnarray*}
\|\mathds{1}_{Q_{\bar R}}\mathcal{L}_{2}(  |\mathds{1}_{Q_{\bar R}}  (\vu + \vb+\vw)|)\|_{\mathcal{M}_{t,x}^{3, \sigma}}\leq C\|\mathcal{L}_{2}(  |\mathds{1}_{Q_{\bar R}}  (\vu + \vb+\vw)|)\|_{\mathcal{M}_{t,x}^{\frac{a}{\lambda},\frac{q}{\lambda}}}\notag\\
\leq C\|\mathds{1}_{Q_{\bar R}}  (\vu + \vb+\vw)\|_{\mathcal{M}_{t,x}^{a, q}}\leq C\|\mathds{1}_{Q_{R_3}}  (\vu + \vb+\vw)\|_{\mathcal{M}_{t,x}^{3, \tau_0}}<+\infty,
\end{eqnarray*}
where in the last estimate we applied again Lemma \ref{lemma_locindi}  (noting that $a\leq 3$ and $q<\tau_0$) and we used the estimates over $\vu, \vb$ and $\vw$ available in (\ref{ConclusionPropositionFirstMorreySpace}).
\item For $\vec{\mathcal{V}}_2$ we write $(\partial_{i}\bar\phi) (\partial_i(\vu+\vb+\vw))=\partial_i((\partial_i\bar\phi)(\vu+\vb+\vw))-(\partial_i^2\bar\phi)(\vu+\vb+\vw)$ and we have
\begin{equation}\label{EstimationPonctuelleVV2}
\begin{split}
|\mathds{1}_{Q_{R_2}}\vec{\mathcal{V}}_2(t,x)|\leq \sum_{i=1}^{3}\left|\mathds{1}_{Q_{R_2}}\int_{0}^{t} e^{ (t-s)\Delta} \partial_i\big((\partial_{i}\bar\phi) (\vu+\vb+\vw)\big)ds\right|\\
+\left|\mathds{1}_{Q_{R_2}}\int_{0}^{t} e^{ (t-s)\Delta} (\partial_{i}^2\bar\phi)(\vu+\vb+\vw)ds\right|.
\end{split}
\end{equation}
Remark that the second term of the right-hand side of  (\ref{EstimationPonctuelleVV2}) can be treated in the same manner as the term $\vec{\mathcal{V}}_1$ so we will only study the first term: by the properties of the heat kernel and by the definition of the Riesz potential $\mathcal{L}_{1}$ (see (\ref{DefinitionPotentielRiesz6})), we obtain
\begin{eqnarray*}
A_2:=\left|\mathds{1}_{Q_{R_2}}\int_{0}^{t} e^{ (t-s)\Delta} \partial_i\big((\partial_{i}\bar\phi) (\vu+\vb+\vw)\big)ds\right|=\left|\mathds{1}_{Q_{R_2}}\int_{0}^{t} \int_{\mathbb{R}^3}\partial_i\mathfrak{g}_{t-s}(x-y)(\partial_{i}\bar\phi) (\vu+\vb+\vw)(s,y)dyds\right|\\
\leq C\mathds{1}_{Q_{R_2}}\int_{\mathbb{R}}\int_{\mathbb{R}^3} \frac{|\mathds{1}_{Q_{\bar R}} (\vu+\vb+\vw)(s,y)|}{(|t-s|^{\frac{1}{2}}+|x-y|)^4}dyds
\leq C\mathds{1}_{Q_{R_2}}(\mathcal{L}_1(|\mathds{1}_{Q_{\bar R}} (\vu+\vb+\vw)|))(t,x).
\end{eqnarray*}
Taking the Morrey $\mathcal{M}_{t,x}^{3, \sigma}$ norm we obtain $\|A_2\|_{\mathcal{M}_{t,x}^{3, \sigma}}\leq C\|\mathds{1}_{Q_{R_2}}(\mathcal{L}_1(|\mathds{1}_{Q_{\bar R}} (\vu+\vb+\vw)|))\|_{\mathcal{M}_{t,x}^{3, \sigma}}$.
Now, for some $4\leq q <5$ we define $\lambda=1-\frac{q}{5}$, noting that $3\leq \frac{3}{\lambda}$ and $\sigma <10<\frac{q}{\lambda}$, by Lemma \ref{Lemme_Hed}, we can write
\begin{eqnarray*}
\|\mathds{1}_{Q_{R_2}}(\mathcal{L}_1(|\mathds{1}_{Q_{\bar R}} (\vu+\vb+\vw)|))\|_{\mathcal{M}_{t,x}^{3, \sigma}}&\leq &C\|\mathcal{L}_1(|\mathds{1}_{Q_{\bar R}} (\vu+\vb+\vw)|)\|_{\mathcal{M}_{t,x}^{\frac{3}{\lambda}, \frac{q}{\lambda}}}\leq C\|\mathds{1}_{Q_{\bar R}} (\vu+\vb+\vw)\|_{\mathcal{M}_{t,x}^{3, q}}\\
&\leq & C\|\mathds{1}_{Q_{R_1}} (\vu+\vb+\vw)\|_{\mathcal{M}_{t,x}^{3, \tau_0}}<+\infty,
\end{eqnarray*}
from which we deduce that $\|\mathds{1}_{Q_{R_2}}\vec{\mathcal{V}}_2\|_{\mathcal{M}_{t,x}^{3, \sigma}}<+\infty$.
\item For the term $\vec{\mathcal{V}}_3$ we have
\begin{eqnarray*}
|\mathds{1}_{Q_{R_2}}\vec{\mathcal{V}}_3(t,x)|&=&\left|\mathds{1}_{Q_{R_2}}\int_{0}^{t} \int_{\mathbb{R}^3}  \mathfrak{g}_{t-s}(x-y)\left[\bar\phi \left( (\vb\cdot\vn)\vu+(\vu\cdot\vn)\vb\right)\right](s,y)dyds\right|\\
&\leq & C\mathds{1}_{Q_{R_2}}\mathcal{L}_2\left(\left|\mathds{1}_{Q_{\bar R}}\left( (\vb\cdot\vn)\vu+(\vu\cdot\vn)\vb\right)\right|\right)(t,x),
\end{eqnarray*}
from which we deduce 
\begin{equation}\label{Decomposition2termesVV3}
\|\mathds{1}_{Q_{R_2}}\vec{\mathcal{V}}_3\|_{\mathcal{M}_{t,x}^{3, \sigma}}\leq C\left\|\mathds{1}_{Q_{R_2}}\mathcal{L}_2\left(|\mathds{1}_{Q_{\bar R}} (\vb\cdot\vn)\vu|\right)\right\|_{\mathcal{M}_{t,x}^{3, \sigma}}+C\left\|\mathds{1}_{Q_{R_2}}\mathcal{L}_2\left(|\mathds{1}_{Q_{\bar R}} (\vu\cdot\vn)\vb|\right)\right\|_{\mathcal{M}_{t,x}^{3, \sigma}}.
\end{equation}
As we have completely symmetric information on $\vu$ and $\vb$ it is enough the study one of these terms and we will treat the first one. We set now $\frac{5}{3-\alpha}<q<\frac{5}{2}$ and $\lambda=1-\frac{2q}{5}$. Since $3\leq \frac{6}{5\lambda}$ and $\tau_0<\sigma<\frac{q}{\lambda}$, applying Lemma \ref{lemma_locindi} and Lemma \ref{Lemme_Hed} we have
$$\left\|\mathds{1}_{Q_{R_2}}\mathcal{L}_2\left(|\mathds{1}_{Q_{\bar R}} (\vb\cdot\vn)\vu|\right)\right\|_{\mathcal{M}_{t,x}^{3, \sigma}}\leq C\left\|\mathds{1}_{Q_{R_2}}\mathcal{L}_2\left(|\mathds{1}_{Q_{\bar R}} (\vb\cdot\vn)\vu|\right)\right\|_{\mathcal{M}_{t,x}^{\frac{6}{5\lambda},\frac{q}{\lambda}}}\leq C\left\|\mathds{1}_{Q_{\bar R}} (\vb\cdot\vn)\vu\right\|_{\mathcal{M}_{t,x}^{\frac{6}{5}, q}}.$$
Recall that we have $\frac{5}{1-\alpha}<\tau_0<\sigma< 10$ and by the H\"older inequality in Morrey spaces (see Lemma \ref{lemma_Product}) we obtain
$$\left\|\mathds{1}_{Q_{\bar R}} (\vb\cdot\vn)\vu\right\|_{\mathcal{M}_{t,x}^{\frac{6}{5}, q}}\leq \left\|\mathds{1}_{Q_{R_3}}\vb\right\|_{\mathcal{M}_{t,x}^{3, \tau_0}}\left\|\mathds{1}_{Q_{R_3}}\vn \otimes \vu\right\|_{\mathcal{M}_{t,x}^{2, \tau_1}}<+\infty,$$
where $\frac{1}{q}=\frac{1}{\tau_0}+\frac{1}{\tau_1}=\frac{2}{\tau_0}+\frac{1}{5}$. Note that the condition $\frac{5}{1-\alpha}<\tau_0<\sigma< 10$ and the relationship $\frac{1}{q}=\frac{2}{\tau_0}+\frac{1}{5}$ are compatible with the fact that $\frac{5}{3-\alpha}<q < \frac{5}{2}$. Applying exactly the same ideas in the second term of (\ref{Decomposition2termesVV3}) we obtain $\|\mathds{1}_{Q_{R_2}}\vec{\mathcal{V}}_3\|_{\mathcal{M}_{t,x}^{3,\sigma}}<+\infty$.
\item The term $\vec{\mathcal{V}}_4$ is the most technical one. Indeed, we write
\begin{eqnarray*}
|\mathds{1}_{Q_{R_2}}\vec{\mathcal{V}}_4|\leq \sum^3_{i,j= 1} \mathds{1}_{Q_{R_2}}\int_{\mathbb{R}} \int_{\mathbb{R}^3}\frac{\left|\left[\bar\phi, \, \frac{ \vn \partial_i\partial_j}{(-\Delta)}\right] (\bar\psi u_i b_j)(s,y)\right|}{(|t-s|^{\frac{1}{2}}+|x-y|)^3}dyds\leq \sum^3_{i,j= 1}\mathds{1}_{Q_{R_2}}\mathcal{L}_{2}\left(\left|\left[\bar\phi, \, \frac{ \vn \partial_i\partial_j}{(-\Delta)}\right] (\bar\psi u_i b_j)\right|\right),
\end{eqnarray*}
and taking the $\mathcal{M}_{t,x}^{3, \sigma}$-norm we have
$\|\mathds{1}_{Q_{R_2}}\vec{\mathcal{V}}_4\|_{\mathcal{M}_{t,x}^{3, \sigma}}\leq \sum^3_{i,j= 1}\left\|\mathds{1}_{Q_{R_2}}\mathcal{L}_{2}\left(\left|\left[\bar\phi, \, \frac{ \vn \partial_i\partial_j}{(-\Delta)}\right] (\bar\psi u_i b_j)\right|\right)\right\|_{\mathcal{M}_{t,x}^{3, \sigma}}$. 
If we set $\frac{1}{q}=  \frac{2}{\tau_0} + \frac{1}{5}$ and $\lambda=1-\frac{2q}{5}$ then we have $3\leq \frac{3}{2\lambda}$ and $\sigma \leq  \frac{q}{\lambda} = \frac{5 \tau_0}{10- \tau_0}$ and by Lemma \ref{lemma_locindi} and Lemma \ref{Lemme_Hed} we obtain:
\begin{eqnarray*}
\left\|\mathds{1}_{Q_{R_2}}\mathcal{L}_{2}\left(\left|\left[\bar\phi, \, \frac{ \vn \partial_i\partial_j}{(-\Delta)}\right] (\bar\psi u_i b_j)\right|\right)\right\|_{\mathcal{M}_{t,x}^{3, \sigma}}&\leq &C\left\|\mathds{1}_{Q_{R_2}}\mathcal{L}_{2}\left(\left|\left[\bar\phi, \, \frac{ \vn \partial_i\partial_j}{(-\Delta)}\right] (\bar\psi u_i b_j)\right|\right)\right\|_{\mathcal{M}_{t,x}^{\frac{3}{2\lambda}, \frac{q}{\lambda}}}\\
&\leq &C\left\|\left[\bar\phi, \, \frac{ \vn \partial_i\partial_j}{(-\Delta)}\right] (\bar\psi u_i b_j)\right\|_{\mathcal{M}_{t,x}^{\frac{3}{2}, q}},
\end{eqnarray*}
We will study this norm and by the definition of Morrey spaces (\ref{DefMorreyparabolico}), if we introduce a threshold $\mathfrak{r}=\frac{\bar R-R_2}{2}$, we have
\end{itemize}
\begin{equation}\label{CommutatorEstimatevVV4}
\begin{split}
\left\|\left[\bar\phi, \, \frac{ \vn \partial_i\partial_j}{(-\Delta)}\right] (\bar\psi u_i b_j)\right\|_{\mathcal{M}_{t,x}^{\frac{3}{2}, q}}^{\frac{3}{2}}&\leq\underset{\underset{0<r<\mathfrak{r}}{(\mathfrak{t},\bar{x})}}{\sup}\;\frac{1}{r^{5(1-\frac{3}{2q})}}\int_{Q_r(\mathfrak{t},\bar x)}\left|\left[\bar\phi, \, \frac{ \vn \partial_i\partial_j}{(-\Delta)}\right] (\bar\psi u_i b_j)\right|^{\frac{3}{2}}dxdt\qquad\\
&+\underset{\underset{\mathfrak{r}< r}{(\mathfrak{t},\bar{x})}}{\sup}\;\frac{1}{r^{5(1-\frac{3}{2q})}}\int_{Q_r(\mathfrak{t},\bar{x})}\left|\left[\bar\phi, \, \frac{ \vn \partial_i\partial_j}{(-\Delta)}\right] (\bar\psi u_i b_j)\right|^{\frac{3}{2}}dxdt.\qquad
\end{split}
\end{equation}
\begin{itemize}
\item[]
Now, we study the second term of the right-hand side above, which is easy to handle as we have $\mathfrak{r}<r$ and we can write
$$\underset{\underset{\mathfrak{r}< r}{(\mathfrak{t},\bar{x}) \in \mathbb{R}\times \R}}{\sup}\;\frac{1}{r^{5(1-\frac{3}{2q})}}\int_{Q_r(\mathfrak{t},\bar{x})}\left|\left[\bar\phi, \, \frac{ \vn \partial_i\partial_j}{(-\Delta)}\right] (\bar\psi u_i b_j)\right|^{\frac{3}{2}}dxdt\leq C_{\mathfrak{r}}\left\|\left[\bar\phi, \, \frac{ \vn \partial_i\partial_j}{(-\Delta)}\right] (\bar\psi u_i b_j)\right\|_{L^{\frac{3}{2}}_{t,x}}^{\frac{3}{2}},$$
and since $\bar\phi$ is a regular function and $\frac{ \vn \partial_i\partial_j}{(-\Delta)}$ is a Calder\'on-Zydmund operator, by the Calder\'on commutator theorem (see the book \cite{PGLR0}), we have that the operator $\left[\bar\phi, \, \frac{ \vn \partial_i\partial_j}{(-\Delta)}\right]$ is bounded in the space $L^{\frac{3}{2}}_{t,x}$ and we can write (using the support properties of $\bar \psi$ given in (\ref{ProprieteLocalisation6}) and the information given in (\ref{ConclusionPropositionFirstMorreySpace})):
\begin{eqnarray*}
\left\|\left[\bar\phi, \, \frac{ \vn \partial_i\partial_j}{(-\Delta)}\right] (\bar\psi u_i b_j)\right\|_{L^{\frac{3}{2}}_{t,x}}&\leq &C\left\|\bar\psi u_i b_j\right\|_{L^{\frac{3}{2}}_{t,x}}\leq C\|\mathds{1}_{Q_{R_1}} u_i b_j\|_{\mathcal{M}^{\frac{3}{2}, \frac{3}{2}}_{t,x}}\\
&\leq & C\|\mathds{1}_{Q_{R_1}} \vu\|_{\mathcal{M}^{3,3}_{t,x}}\|\mathds{1}_{Q_{R_1}} \vb\|_{\mathcal{M}^{3, 3}_{t,x}}\leq C\|\mathds{1}_{Q_{R_1}} \vu\|_{\mathcal{M}^{3,\tau_0}_{t,x}}\|\mathds{1}_{Q_{R_1}} \vb\|_{\mathcal{M}^{3, \tau_0}_{t,x}}<+\infty,
\end{eqnarray*}
where in the last line we used H\"older inequalities in Morrey spaces and we applied Lemma \ref{lemma_locindi}.\\

The first term of the right-hand side of  (\ref{CommutatorEstimatevVV4}) requires some extra computations: indeed, as we are interested to obtain information over the parabolic ball $Q_{r}(\mathfrak{t}, \bar{x})$ we can write
for some  $0<r<\mathfrak{r}$:
\begin{equation}\label{CommutatorEstimatevVV401}
\mathds{1}_{Q_{r}}\left[\bar\phi, \, \frac{ \vn \partial_i\partial_j}{(-\Delta)}\right] (\bar\psi u_i b_j))=\mathds{1}_{Q_{r}}\left[\bar\phi, \, \frac{ \vn \partial_i\partial_j}{(-\Delta)}\right] (\mathds{1}_{Q_{2r}}\bar\psi u_i b_j)+\mathds{1}_{Q_{r}}\left[\bar\phi, \, \frac{ \vn \partial_i\partial_j}{(-\Delta)}\right] ((\mathbb{I}-\mathds{1}_{Q_{2r}})\bar\psi u_i b_j),
\end{equation}
and as before we will study the $L^{\frac{3}{2}}_{t,x}$ norm of these two terms. For the first quantity in the right-hand side of (\ref{CommutatorEstimatevVV401}), by the Calder\'on commutator theorem, by the definition of Morrey spaces and by the H\"older inequalities we have 
\begin{eqnarray*}
\left\|\mathds{1}_{Q_{r}}\left[\bar\phi, \, \frac{ \vn \partial_i\partial_j}{(-\Delta)}\right] (\mathds{1}_{Q_{2r}}\bar\psi u_i b_j)\right\|_{L^{\frac{3}{2}}_{t,x}}^{\frac{3}{2}}&\leq& C\|\mathds{1}_{Q_{2r}}\bar\psi u_i b_j\|_{L^{\frac{3}{2}}_{t,x}}^{\frac{3}{2}}\leq Cr^{5 (1-\frac{3}{\tau_0})} \|\mathds{1}_{Q_{R_1}}u_i b_j\|_{\mathcal{M}^{\frac{3}{2}, \frac{\tau_0}{2}}_{t,x}}^{\frac{3}{2}}\\
&\leq & Cr^{5 (1-\frac{3}{\tau_0})} \|\mathds{1}_{Q_{R_1}}\vu\|_{\mathcal{M}^{3, \tau_0}_{t,x}}^{\frac{3}{2}}\|\mathds{1}_{Q_{R_1}}\vb\|_{\mathcal{M}^{3, \tau_0}_{t,x}}^{\frac{3}{2}},
\end{eqnarray*}
for all $0<r<\mathfrak{r}$, from which we deduce that 
$$\underset{\underset{0<r<\mathfrak{r}}{(\mathfrak{t},\bar{x}) }}{\sup}\;\frac{1}{r^{5(1-\frac{3}{2q})}}\int_{Q_r(\mathfrak{t},\bar{x})}\left|\mathds{1}_{Q_{r}}\left[\bar\phi, \, \frac{ \vn \partial_i\partial_j}{(-\Delta)}\right] (\mathds{1}_{Q_{2r}}\bar\psi u_i b_j)\right|^{\frac{3}{2}}dxdt\leq C \|\mathds{1}_{Q_{R_1}}\vu\|_{\mathcal{M}^{3, \tau_0}_{t,x}}^{\frac{3}{2}}\|\mathds{1}_{Q_{R_1}}\vb\|_{\mathcal{M}^{3, \tau_0}_{t,x}}^{\frac{3}{2}}<+\infty.$$
We study now the second term of the right-hand side of (\ref{CommutatorEstimatevVV401}) and for this we consider the following operator:
$$T: f \mapsto  \left(\mathds{1}_{Q_{r}}  \left[\bar\phi, \, \frac{ \vn \partial_i \partial_j}{- \Delta }\right] (\mathbb{I} - \mathds{1}_{Q_{2r}}) \bar \varphi
\right) f,$$
and by the properties of the convolution kernel of the operator $\frac{1}{(-\Delta)}$ we obtain
$$|T(f)(x)|\leq C\mathds{1}_{Q_{r}}(x)\int_{\mathbb{R}^3}\frac{(\mathbb{I} - \mathds{1}_{Q_{2r}})(y) \mathds{1}_{Q_{R_1}}(y) |f(y)| |\bar \phi(x)-\bar \phi(y)|}{|x-y|^4} dy.$$
Recalling that $0<r<\mathfrak{r}=\frac{\bar R-R_2}{2}$, by the support properties of the test function $\bar\phi$ (see (\ref{ProprieteLocalisation6})), the integral above is meaningful if $|x-y|>r$ and thus we can write
\begin{eqnarray*}
\left\|\mathds{1}_{Q_{r}}\left[\bar\phi, \, \frac{ \vn \partial_i\partial_j}{(-\Delta)}\right] ((\mathbb{I}-\mathds{1}_{Q_{2r}})\bar\psi u_i b_j)\right\|_{L^\frac{3}{2}_{t,x}}^\frac{3}{2}\leq C\left\|\mathds{1}_{Q_{r}} \int_{\mathbb{R}^3} \frac{\mathds{1}_{|x-y| > r}}{|x-y|^4}(\mathbb{I} - \mathds{1}_{Q_{2r}})(y) \mathds{1}_{Q_{R_1}}(y)|u_i b_j |dy\right\|_{L^\frac{3}{2}_{t,x}}^\frac{3}{2}\\
\leq C\left(\int_{|y|>r}\frac{1}{|y|^4}\| \mathds{1}_{Q_{R_1}}
|u_i b_j |(\cdot-y)\|_{L^\frac{3}{2}_{t,x}(Q_r)}dy\right)^\frac{3}{2}\leq Cr^{-\frac{3}{2}}\| \mathds{1}_{Q_{R_1}}u_i b_j \|_{L^\frac{3}{2}_{t,x}(Q_{r})}^\frac{3}{2},
\end{eqnarray*}
with this estimate at hand and using the definition of Morrey spaces, we can write
\begin{eqnarray*}
\int_{Q_r(\mathfrak{t},\bar{x})}\left|\mathds{1}_{Q_{r}}\left[\bar\phi, \, \frac{ \vn \partial_i\partial_j}{(-\Delta)}\right] ((\mathbb{I}-\mathds{1}_{Q_{2r}})\bar\psi u_i b_j)\right|^{\frac{3}{2}}dxdt&\leq &Cr^{-\frac{3}{2}}r^{5 (1-\frac{3}{\tau_0})}\| \mathds{1}_{Q_{R_1}}u_i b_j \|_{\mathcal{M}^{\frac{3}{2}, \frac{\tau_0}{2}}_{t,x}}^\frac{3}{2}\\
&\leq &Cr^{5(1-\frac{3}{2q})}\| \mathds{1}_{Q_{R_1}}u_i b_j \|_{\mathcal{M}^{\frac{3}{2}, \frac{\tau_0}{2}}_{t,x}}^\frac{3}{2},
\end{eqnarray*}
where in the last inequality we used the fact that $\frac{1}{q}=  \frac{2}{\tau_0} + \frac{1}{5}$, which implies  $r^{-\frac{3}{2}}r^{5 (1-\frac{3}{\tau_0})}= r^{5(1-\frac{3}{2q})}$. Thus we finally obtain
$$\underset{\underset{0<r<\mathfrak{r}}{(\mathfrak{t},\bar{x}) }}{\sup}\;\frac{1}{r^{5(1-\frac{3}{2q})}}\int_{Q_r(\mathfrak{t},\bar{x})}\left|\mathds{1}_{Q_{r}}\left[\bar\phi, \, \frac{ \vn \partial_i\partial_j}{(-\Delta)}\right] ((\mathbb{I}-\mathds{1}_{Q_{2r}})\bar\psi u_i b_j)\right|^{\frac{3}{2}}dxdt\leq C \|\mathds{1}_{Q_{R_1}}\vu\|_{\mathcal{M}^{3, \tau_0}_{t,x}}^{\frac{3}{2}}\|\mathds{1}_{Q_{R_1}}\vb\|_{\mathcal{M}^{3, \tau_0}_{t,x}}^{\frac{3}{2}}<+\infty.$$
We have proven that all the term in (\ref{CommutatorEstimatevVV4}) are bounded and we can conclude that $\|\mathds{1}_{Q_{R_2}}\vec{\mathcal{V}}_4\|_{\mathcal{M}_{t,x}^{3, \sigma}}<+\infty$.

\item For the quantity $\vec{\mathcal{V}}_5$, based in the expression (\ref{FormulePourN}) we write
\begin{eqnarray*}
|\mathds{1}_{Q_{R_2}}\vec{\mathcal{V}}_5(t,x)|&\leq &C\sum^3_{i,j= 1} \mathds{1}_{Q_{R_2}}\int_{\mathbb{R}}\int_{\mathbb{R}^3} \frac{|\mathcal{R}_i\mathcal{R}_j( \bar\phi u_i b_j )(s,y)|}{(|t-s|^{\frac{1}{2}}+|x-y|)^4}dyds\leq C\sum^3_{i,j= 1} \mathds{1}_{Q_{R_2}} \mathcal{L}_{1}\left(|\mathcal{R}_i\mathcal{R}_j( \bar\phi u_i b_j )|\right)(t,x),
\end{eqnarray*}
where we used the decaying properties of the heat kernel (recall that $\mathcal{R}_i=\frac{\partial_i}{\sqrt{- \Delta}}$ are the Riesz transforms). 
Now taking the Morrey $\mathcal{M}^{3, \sigma}_{t,x}$ norm and by Lemma \ref{lemma_locindi} (with $\nu=\frac{4\tau_0+5}{5\tau_0}$, $p=3$, $q=\tau_0$ such that $\frac{p}{\nu}>3$ and $\frac{q}{\nu}>\sigma$ which is compatible with the condition $\tau_0<\sigma$) we have
\begin{eqnarray*}
\|\mathds{1}_{Q_{R_2}}\vec{\mathcal{V}}_5\|_{\mathcal{M}^{3, \sigma}_{t,x}}&\leq &C \sum^3_{i,j= 1}\|  \mathds{1}_{Q_{R_2}}\mathcal{L}_{1}\left(|\mathcal{R}_i\mathcal{R}_j( \bar\phi u_i b_j )|\right)\|_{\mathcal{M}^{\frac{p}{\nu}, \frac{q}{\nu}}_{t,x}}
\end{eqnarray*}
Then by Lemma \ref{Lemme_Hed} with $\lambda= 1 - \tfrac{\tau_0 /2}{5}$ (recall $\frac{5}{1-\alpha}<\tau_0 <10$ so that $\nu > 2 \lambda$) and by the boundedness of Riesz transforms in Morrey spaces we obtain:
\begin{eqnarray*}
\|\mathds{1}_{Q_{R_2}} \mathcal{L}_{1}\left(|\mathcal{R}_i\mathcal{R}_j( \bar\phi u_i b_j )|\right)\|_{\mathcal{M}^{\frac{p}{\nu}, \frac{q}{\nu}}_{t,x}}
&\leq&C\|\mathcal{L}_{1}\left(|\mathcal{R}_i\mathcal{R}_j( \bar\phi u_i b_j )|\right)\|_{\mathcal{M}^{\frac{p}{2\lambda}, \frac{q}{2\lambda}}_{t,x}}
\leq C \|\mathcal{R}_i\mathcal{R}_j( \bar\phi u_i b_j )\|_{\mathcal{M}^{\frac{3}{2},\frac{\tau_0}{2}}_{t,x}}\\
&\leq &\| \mathds{1}_{Q_{R_1}} u_i b_j \|_{\mathcal{M}^{\frac{3}{2},\frac{\tau_0}{2}}_{t,x}}
\leq  C\|\mathds{1}_{Q_{R_1}} \vu\|_{\mathcal{M}^{3,\tau_0}_{t,x}}\|  \mathds{1}_{Q_{R_1}} \vb\|_{\mathcal{M}^{3,\tau_0}_{t,x}}<+\infty.
\end{eqnarray*}
\item The quantities $\vec{\mathcal{V}}_6$ and $\vec{\mathcal{V}}_7$ based in the corresponding terms of (\ref{FormulePourN}) can be treated in a very similar fashion since their inner structure is essentially the same. We thus only treat here the term $\vec{\mathcal{V}}_6$ and following the same ideas we have
$$|\mathds{1}_{Q_{R_2}}\vec{\mathcal{V}}_6|\leq C\sum^3_{i,j= 1}  \mathds{1}_{Q_{R_2}}\int_{\mathbb{R}}  \int_{\mathbb{R}^3}\frac{\left|\frac{\bar\phi \vn\partial_i}{(- \Delta )} (\partial_j \bar\psi) u_i b_j(s,y)\right|}{(|t-s|^{\frac{1}{2}}+|x-y|)^3} dyds=C\sum^3_{i,j= 1}  \mathds{1}_{Q_{R_2}}\mathcal{L}_{2}\left(\left|\frac{\bar\phi \vn\partial_i}{(- \Delta )} (\partial_j \bar\psi) u_i b_j\right|\right).$$
For $2<q<\frac{5}{2}$, define $\lambda=1-\frac{2q}{5}$, we thus have $3\leq \frac{3}{2\lambda}$ and $\sigma<10\leq \frac{q}{\lambda}$. Then, by Lemma \ref{lemma_locindi} and Lemma \ref{Lemme_Hed} we can write
\begin{eqnarray*}
\left\|\mathds{1}_{Q_{R_2}}\mathcal{L}_{2}\left|\frac{\bar\phi \vn\partial_i}{(- \Delta )} (\partial_j \bar\psi) u_i b_j\right|\right\|_{\mathcal{M}^{3, \sigma}_{t,x}}\leq C\left\|\mathds{1}_{Q_{R_2}}\mathcal{L}_{2}\left|\frac{\bar\phi \vn\partial_i}{(- \Delta )} (\partial_j \bar\psi) u_i b_j\right|\right\|_{\mathcal{M}^{\frac{3}{2\lambda}, \frac{q}{\lambda}}_{t,x}}\leq  C\left\|\frac{\bar\phi \vn\partial_i}{(- \Delta )} (\partial_j \bar\psi) u_i b_j\right\|_{\mathcal{M}^{\frac{3}{2}, q}_{t,x}},
\end{eqnarray*}
but since the operator $\frac{\bar\phi \vn\partial_i}{(- \Delta )}$ is bounded in Morrey spaces and since $2<q<\frac{5}{2}< \tfrac{\tau_0}{2}$, one has by Lemma \ref{lemma_locindi} and by the H\"older inequalities
\begin{eqnarray*}
\left\|\frac{\bar\phi \vn\partial_i}{(- \Delta )} (\partial_j \bar\psi) u_i b_j\right\|_{\mathcal{M}^{\frac{3}{2}, q}_{t,x}}\leq C\left\| \mathds{1}_{Q_{R_1}}u_i b_j\right\|_{\mathcal{M}^{\frac{3}{2}, q}_{t,x}}\leq C\| \mathds{1}_{Q_{R_1}}u_i b_j\|_{\mathcal{M}^{\frac{3}{2}, \frac{\tau_0}{2}}_{t,x}}\leq C\|\mathds{1}_{Q_{R_1}} \vu\|_{\mathcal{M}^{3,\tau_0}_{t,x}}\|  \mathds{1}_{Q_{R_1}} \vb\|_{\mathcal{M}^{3,\tau_0}_{t,x}},
\end{eqnarray*}
from which we deduce $\|\mathds{1}_{Q_{R_2}}\vec{\mathcal{V}}_6\|_{\mathcal{M}^{3, \sigma}_{t,x}}<+\infty$. The same computations can be performed to obtain that $\|\mathds{1}_{Q_{R_2}}\vec{\mathcal{V}}_7\|_{\mathcal{M}^{3, \sigma}_{t,x}}<+\infty$.

\item The quantity $\vec{\mathcal{V}}_8$ based in (\ref{FormulePourN}) is treated in the following manner: we first write
$$\|\mathds{1}_{Q_{R_2}}\vec{\mathcal{V}}_8\|_{\mathcal{M}^{3, \sigma}_{t,x}}\leq C \sum^3_{i,j= 1}\left\| \mathds{1}_{Q_{R_2}}\left(\mathcal{L}_2\left|\bar\phi \frac{\vn}{(-\Delta )}(\partial_i \partial_j\bar \psi) (u_i b_j)\right|\right)\right\|_{\mathcal{M}^{3,\sigma}_{t,x}}.$$
We set $1<\nu<\frac{3}{2}$, $2\nu <q<\frac{5\nu}{2}$ and $\lambda=1-\frac{2q}{5\nu}$, thus we have $3\leq \frac{\nu}{\lambda}$ and $\sigma <10< \frac{q}{\lambda}$, then, by Lemma \ref{lemma_locindi} and by Lemma \ref{Lemme_Hed} we can write
\begin{eqnarray}
\left\|\mathds{1}_{Q_{R_2}}\left(\mathcal{L}_2\left|\bar\phi \frac{\vn}{(-\Delta )}(\partial_i \partial_j\bar \psi) (u_i b_j)\right|\right)\right\|_{\mathcal{M}^{3,\sigma}_{t,x}}\leq C\left\| \mathds{1}_{Q_{R_2}}\left(\mathcal{L}_2\left|\bar\phi \frac{\vn}{(-\Delta )}(\partial_i \partial_j\bar \psi) (u_i b_j)\right|\right)\right\|_{\mathcal{M}^{\frac{\nu}{\lambda},\frac{q}{\lambda}}_{t,x}}&&\notag\\
\leq C\left\| \bar\phi \frac{\vn}{(-\Delta )}(\partial_i \partial_j\bar \psi) (u_i b_j)\right\|_{\mathcal{M}^{\nu,q}_{t,x}}\leq C\left\| \bar\phi \frac{\vn}{(-\Delta )}(\partial_i \partial_j\bar \psi) (u_i b_j)\right\|_{\mathcal{M}^{\nu,\frac{5\nu}{2}}_{t,x}}\leq C\left\| \bar\phi \frac{\vn}{(-\Delta )}(\partial_i \partial_j\bar \psi) (u_i b_j)\right\|_{L^{\nu}_tL^{\infty}_x}\label{Formula_intermediairevVV80}
\end{eqnarray}
where in the last estimate we used the space inclusion $L^{\nu}_tL^{\infty}_x\subset \mathcal{M}^{\nu,\frac{5\nu}{2}}_{t,x}$. Let us focus now in the $L^\infty$ norm above (\emph{i.e.} without considering the time variable). Remark that due to the support properties of the auxiliary function $\bar\psi$ given in (\ref{ProprieteLocalisation6}) we have $supp(\partial_i \partial_j\bar \psi) =Q_{R_1}\setminus Q_{\tilde{R}}$ and recall by (\ref{ProprieteLocalisation6}) we have $supp\; \bar \phi = Q_{\bar R}$ where $\bar{R}<\tilde{R}<R_1$, thus by the properties of the kernel of the operator $\frac{\vn}{(-\Delta)}$ we can write
\begin{eqnarray}
\left| \bar\phi \frac{\vn}{(-\Delta )}(\partial_i \partial_j\bar \psi) (u_i b_j)\right|&\leq& C\left|\int_{\mathbb{R}^3} \frac{1}{|x-y|^2}\mathds{1}_{Q_{\bar R}}(x)\mathds{1}_{Q_{R_1}\setminus Q_{\tilde{R}}}(y)(\partial_i \partial_j\bar \psi) (u_i b_j)(\cdot,y)dy\right|\notag\\
&\leq & C\left|\int_{\mathbb{R}^3} \frac{\mathds{1}_{|x-y|>(\tilde{R}-\bar R)}}{|x-y|^2}\mathds{1}_{Q_{\bar R}}(x)\mathds{1}_{Q_{R_1}\setminus Q_{\tilde{R}}}(y)(\partial_i \partial_j\bar \psi) (u_i b_j)(\cdot,y)dy\right|,\label{Formula_intermediairevVV801}
\end{eqnarray}
and the previous expression is nothing but the convolution between the function $(\partial_i \partial_j\bar \psi) (u_i b_j)$ and a $L^\infty$-function, thus we have 
\begin{equation}\label{Formula_intermediairevVV81}
\left\| \bar\phi \frac{\vn}{(-\Delta )}(\partial_i \partial_j\bar \psi) (u_i b_j) (t,\cdot)\right\|_{L^\infty}\leq C\|(\partial_i \partial_j\bar \psi) (u_i b_j)(t,\cdot)\|_{L^1}\leq C\|\mathds{1}_{Q_{R_1}}(u_i b_j)(t,\cdot)\|_{L^{\nu}},
\end{equation}
and taking the $L^\nu$-norm in the time variable we obtain
\begin{eqnarray*}
\left\| \bar\phi \frac{\vn}{(-\Delta )}(\partial_i \partial_j\bar \psi) (u_i b_j)\right\|_{L^{\nu}_t L^{\infty}_{x}}
&\leq &C\|\mathds{1}_{Q_{R_1}}u_i b_j\|_{L^{\nu}_{t,x}}\leq C\|\mathds{1}_{Q_{R_1}}\vu\|_{\mathcal{M}^{3,\tau_0}_{t,x}}\|\mathds{1}_{Q_{R_1}}\vb\|_{\mathcal{M}^{3,\tau_0}_{t,x}}<+\infty,
\end{eqnarray*}
where we used the fact that $1<\nu<\frac{3}{2}<\frac{\tau_0}{2}$ and we applied Hölder's inequality. Gathering together all these estimates we obtain 
$\|\mathds{1}_{Q_{R_2}}\vec{\mathcal{V}}_8\|_{\mathcal{M}^{3, \sigma}_{t,x}}<+\infty$.\\

The terms $\vec{\mathcal{V}}_9, \cdots,\vec{\mathcal{V}}_{18}$ are studied in the following lemma.
\begin{lem}
\begin{itemize}
\item[]
\item[1)]The quantities $\vec{\mathcal{V}}_9$ and $\vec{\mathcal{V}}_{14}$ based in the corresponding terms of (\ref{FormulePourN}) can be treated in the same way as the term $\vec{\mathcal{V}}_4$.
\item[2)] The terms $\vec{\mathcal{V}}_{10}$ and $\vec{\mathcal{V}}_{15}$ are controlled as $\vec{\mathcal{V}}_5$.
\item[3)] The terms $\vec{\mathcal{V}}_{11}$, $\vec{\mathcal{V}}_{12}$, $\vec{\mathcal{V}}_{16}$ and $\vec{\mathcal{V}}_{17}$ are controlled as $\vec{\mathcal{V}}_6$.
\item[4)] The terms $\vec{\mathcal{V}}_{13}$ and $\vec{\mathcal{V}}_{18}$  are controlled as $\vec{\mathcal{V}}_8$.
\end{itemize}
\end{lem}
{\bf Proof.} Following the estimates given previously for the terms $\vec{\mathcal{V}}_4$, $\vec{\mathcal{V}}_5$, $\vec{\mathcal{V}}_6$ and $\vec{\mathcal{V}}_8$, all the terms $\vec{\mathcal{V}}_9, \cdots,\vec{\mathcal{V}}_{18}$ can be controlled by the quantities  $\|\mathds{1}_{Q_{R_1}}\vu\|_{\mathcal{M}^{3,\tau_0}_{t,x}}$, $\|\mathds{1}_{Q_{R_1}}\vb\|_{\mathcal{M}^{3,\tau_0}_{t,x}}$ or $\|\mathds{1}_{Q_{R_1}}\va\|_{\mathcal{M}^{3,\tau_0}_{t,x}}$. It is enough to observe that we have
$\|\mathds{1}_{Q_{R_1}}\va\|_{\mathcal{M}^{3,\tau_0}_{t,x}}\leq C\|\mathds{1}_{Q_{R_1}}\va\|_{\mathcal{M}^{6,6}_{t,x}}=\|\va\|_{L^{6,6}_{t,x}(\Omega)}<+\infty$ since $\frac{5}{1-\alpha}<\tau_0<6$, which is possible if $0<\alpha<\frac{1}{12}$. \hfill $\blacksquare$

\item The quantity $\vec{\mathcal{V}}_{19}$ based in (\ref{FormulePourN})  can be treated in the same way as the term $\vec{\mathcal{V}}_8$. Indeed, by the same arguments displayed to deduce (\ref{Formula_intermediairevVV80}), we can write (recall that $1<\nu<\frac{3}{2}$): $\displaystyle{
\|\mathds{1}_{Q_{R_2}}\vec{\mathcal{V}}_{19}\|_{\mathcal{M}^{3, \sigma}_{t,x}}\leq C\left\| \bar\phi \frac{\vn}{(-\Delta )}( (\Delta \bar\psi) p)\right\|_{L^{\nu}_t L^{\infty}_{x}}}$ and if we study the $L^\infty$-norm in the space variable of this term, by the same ideas used in (\ref{Formula_intermediairevVV801})-(\ref{Formula_intermediairevVV81}) we obtain $\left\| \bar\phi \frac{\vn}{(-\Delta )}( (\Delta \bar\psi) p) (t,\cdot)\right\|_{L^{\infty}}\leq C\|(\Delta \bar\psi) p (t,\cdot)\|_{L^1}\leq C\|\mathds{1}_{Q_{R_1}}p (t,\cdot)\|_{L^{\nu}}$. Thus, taking the $L^{\nu}$-norm in the time variable we have
$$\|\mathds{1}_{Q_{R_2}}\vec{\mathcal{V}}_{19}\|_{\mathcal{M}^{3, \sigma}_{t,x}}
\leq C\left\| \bar\phi \frac{\vn}{(-\Delta )}( (\Delta \bar\psi) p)\right\|_{L^{\nu}_t L^{\infty}_{x}}\leq C \|\mathds{1}_{Q_{R_1}}p\|_{L^{\nu}_{t,x}}\leq C \|\mathds{1}_{Q_{R_1}}p\|_{L^{\frac32}_{t,x}}<+\infty.$$

\item The study of the quantity $\vec{\mathcal{V}}_{20}$ follows almost the same lines as the terms $\vec{\mathcal{V}}_{8}$ and $\vec{\mathcal{V}}_{9}$. However instead of (\ref{Formula_intermediairevVV801}) we have
$$\left|\bar\phi \frac{\vn\partial_i }{(-\Delta)} ( (\partial_i \bar\psi)p)\right|
\leq  C\left|\int_{\mathbb{R}^3} \frac{\mathds{1}_{|x-y|>(\tilde{R}-\bar R)}}{|x-y|^3}\mathds{1}_{Q_{\bar R}}(x)\mathds{1}_{Q_{R_1}\setminus Q_{\tilde{R}}}(y)(\partial_i \bar \psi) p(t,y)dy\right|,$$
and thus we can write:
$$\|\mathds{1}_{Q_{R_2}}\vec{\mathcal{V}}_{20}\|_{\mathcal{M}^{3, \sigma}_{t,x}}\leq  \left\|\bar\phi \frac{\vn\partial_i }{(-\Delta)} ( (\partial_i \bar\psi) p)\right\|_{L^{\nu}_t L^{\infty}_{x}}\leq C \|\mathds{1}_{Q_{R_1}}p\|_{L^{\nu}_{t,x}}\leq C \|\mathds{1}_{Q_{R_1}}p\|_{L^{\frac32}_{t,x}}<+\infty.$$

\item For the term $\vec{\mathcal{V}}_{21}$ based in (\ref{FormulePourN}) can be treated in the same manner as $\vec{\mathcal{V}}_{2}$ and we easily obtain $\|\mathds{1}_{Q_{R_2}}\vec{\mathcal{V}}_{21}\|_{\mathcal{M}_{t,x}^{3, \sigma}}<+\infty$

\item The study of the quantity $\vec{\mathcal{V}}_{22}$ is easy to handle, indeed, we have 
\begin{eqnarray*}
|\mathds{1}_{Q_{R_2}}\vec{\mathcal{V}}_{22}|&\leq &\left|\mathds{1}_{Q_{R_2}}\int_{0}^{t} e^{ (t-s)\Delta} \bar\phi(\vf+\vg)ds\right|\leq C\mathds{1}_{Q_{R_2}}\int_{\mathbb{R}}\int_{\mathbb{R}^3}\frac{|\bar\phi(\vf+\vg)(s,y)|}{(|t-s|^{\frac{1}{2}}+|x-y|)^3}dyds\\
&\leq& C\mathds{1}_{Q_{R_2}}\mathcal{L}_{2}(\mathds{1}_{Q_{R_3}}|\vf+\vg|)(t,x),
\end{eqnarray*} 
and taking the Morrey $\mathcal{M}_{t,x}^{3, \sigma}$ norm we obtain
$\|\mathds{1}_{Q_{R_2}}\vec{\mathcal{V}}_{22}\|_{\mathcal{M}_{t,x}^{3, \sigma}}\leq C\|\mathds{1}_{Q_{R_2}}\mathcal{L}_{2}(\mathds{1}_{Q_{R_1}}|\vf+\vg|)\|_{\mathcal{M}_{t,x}^{3, \sigma}}$, then if we set $\frac{11}{5}<q<\frac{5}{2}$ and $\lambda=1-\frac{2q}{5}$ we thus have $3\leq \frac{10}{7\lambda}$ and $\sigma<10<\frac{q}{\lambda}$. Now by Lemma \ref{lemma_locindi} and Lemma \ref{Lemme_Hed} we have
$\|\mathds{1}_{Q_{R_2}}\mathcal{L}_{2}(\mathds{1}_{Q_{R_1}}|\vf+\vg|)\|_{\mathcal{M}_{t,x}^{3, \sigma}}\leq C\|\mathcal{L}_{2}(\mathds{1}_{Q_{R_1}}|\vf+\vg|)\|_{\mathcal{M}_{t,x}^{\frac{10}{7\lambda}, \frac{q}{\lambda}}}\leq C\|\mathds{1}_{Q_{R_1}}|\vf+\vg|\|_{\mathcal{M}_{t,x}^{\frac{10}{7}, q}}$ but since $q<\frac{5}{2}<\frac{5}{2-\alpha}<\tau_a,\tau_b$, by  Lemma \ref{lemma_locindi} we obtain 
$$\|\mathds{1}_{Q_{R_1}}|\vf+\vg|\|_{\mathcal{M}_{t,x}^{\frac{10}{7}, q}}\leq C\left(\|\mathds{1}_{Q_{R_1}}\vf\|_{\mathcal{M}_{t,x}^{\frac{10}{7}, \tau_a}}+\|\mathds{1}_{Q_{R_1}}\vg\|_{\mathcal{M}_{t,x}^{\frac{10}{7}, \tau_b}}\right)<+\infty,$$
thus, gathering all the estimates above we have $\|\mathds{1}_{Q_{R_2}}\vec{\mathcal{V}}_{22}\|_{\mathcal{M}_{t,x}^{3, \sigma}}<+\infty$.

\item For the quantity $\vec{\mathcal{V}}_{23}$ of (\ref{FormulePourN}) we first note that the quantity $\bar\phi div((\vu+\vb) \otimes \va+\va \otimes (\vu+\vb))$ can be decomposed as $\bar\phi \partial_i(u_ja_k)$ with $1\leq i,j,k\leq 3$ (and other similar terms with $b_j$ instead of $u_j$) and thus we have:
\begin{align*}
\left|\mathds{1}_{Q_{R_2}}\int_0^t\int_{\mathbb{R}^3}\mathfrak{g_{t-s}}(x-y)[\bar \phi\partial_i(u_ja_k)](s,y)dyds \right|&\leq\bigg|\mathds{1}_{Q_{R_2}}\int_0^t\int_{\mathbb{R}^3}\partial_i\mathfrak{g_{t-s}}(x-y)[ \bar  \phi u_ja_k](s,y)dyds\bigg| \\
&+\bigg|\mathds{1}_{Q_{R_2}}\int_0^t\int_{\mathbb{R}^3}\partial_i\mathfrak{g_{t-s}}(x-y)[\partial_i( \bar \phi)(u_ja_k)](s,y)dyds \biggr|,
\end{align*}
and by the same arguments as in the previous lines we obtain
\begin{equation}\label{Estimation23}
\begin{split}
\left\|\mathds{1}_{Q_{R_2}}\int_0^t\int_{\mathbb{R}^3}\mathfrak{g_{t-s}}(x-y)[\bar \phi \partial_j(u_ja_k)](s,y)dyds\right\|_{\M^{3,\sigma}}\le &C\Big(\|\mathds{1}_{Q_{R_2}}\mathcal{L}_{1}|\mathds{1}_{Q_{\bar R}}u_ja_k|\|_{\M^{3,\sigma}}\\
&+\|\mathds{1}_{Q_{R_2}}\mathcal{L}_{2}|\mathds{1}_{Q_{\bar R}}u_ja_k|\|_{\M^{3,\sigma}}\Big).
\end{split}
\end{equation}
For the first term of the right-hand side above we set $p=2$, $q=\frac{6\tau_0}{6+\tau_0}$ and $\lambda=\frac{30-\tau_0}{30+5\tau_0}$. Note that $\frac{p}{\lambda}\geq 3$ and $\frac{q}{\lambda}\geq \sigma$ (if $\sigma>\tau_0>5$ is close enough to $\tau_0$) and thus, by the Lemma \ref{lemma_locindi} and by Lemma \ref{Lemme_Hed}, we have $\|\mathds{1}_{Q_{R_2}}\mathcal{L}_{1}|\mathds{1}_{Q_{\bar R}}u_ja_k|\|_{\M^{3,\sigma}}\leq C\|\mathcal{L}_{1}|\mathds{1}_{Q_{\bar R}}u_ja_k|\|_{\M^{\frac{p}{\lambda},\frac{q}{\lambda}}}\leq C\|\mathds{1}_{Q_{\bar R}}u_ja_k\|_{\M^{p,q}}=C\|\mathds{1}_{Q_{R_1}}u_ja_k\|_{\M^{2,\frac{6\tau_0}{6+\tau_0}}}$ and by the H\"older inequalities in the Morrey spaces we obtain
$\|\mathds{1}_{Q_{R_1}}u_ja_k\|_{\M^{2,\frac{6\tau_0}{6+\tau_0}}}\leq \|\mathds{1}_{Q_{R_1}}u_j\|_{\M^{3,\tau_0}}\|\mathds{1}_{Q_{R_1}}a_k\|_{\M^{6,6}}= \|\mathds{1}_{Q_{R_1}}u_j\|_{\M^{3,\tau_0}}\|a_k\|_{L^{6}_{t,x}(\Omega)}<+\infty$.\\

For the second term of the right-hand side of (\ref{Estimation23}), we fix $p,q=2$ and $\lambda=\frac15$ and we have $\frac{p}{\lambda}\geq 3$ and $\frac{q}{\lambda}\geq \sigma$. Thus, by the same arguments as above we can write
\begin{eqnarray*}
\|\mathds{1}_{Q_{R_2}}\mathcal{L}_{2}|\mathds{1}_{Q_{\bar R}}u_ja_k|\|_{\M^{3,\sigma}}&\leq &C\|\mathcal{L}_{2}|\mathds{1}_{Q_{\bar R}}u_ja_k|\|_{\M^{\frac{p}{\lambda},\frac{q}{\lambda}}}\leq C\|\mathds{1}_{Q_{\bar R}}u_ja_k\|_{\M^{p,q}}=C\|\mathds{1}_{Q_{R_1}}u_ja_k\|_{\M^{2,2}}\\
&\leq &C\|\mathds{1}_{Q_{R_1}}u_ja_k\|_{\M^{2,\frac{6\tau_0}{6+\tau_0}}}\leq \|\mathds{1}_{Q_{R_1}}u_j\|_{\M^{3,\tau_0}}\|a_k\|_{L^{6}_{t,x}(\Omega)}<+\infty.
\end{eqnarray*}
Applying these estimates to all the terms of the form $\bar\phi \partial_i(u_ja_k)$ and $\bar\phi \partial_i(b_ja_k)$ we finally obtain that $\|\mathds{1}_{Q_{R_2}}\vec{\mathcal{V}}_{23}\|_{\M^{3,\sigma}}<+\infty$.

\item For the last term $\vec{\mathcal{V}}_{24}$ given by the corresponding quantity in (\ref{FormulePourN}), we have
\begin{equation}\label{DernierTermeN}
|\mathds{1}_{Q_{R_2}}\vec{\mathcal{V}}_{24}|=\bigg|\mathds{1}_{Q_{R_2}}\int_0^t\int_{\mathbb{R}^3}\mathfrak{g_{t-s}}(x-y)\bar\phi \bigg(\underbrace{\grad div(\vw)}_{(a)}-\underbrace{\vw}_{(b)}-\frac{1}{2}\underbrace{((\vu+\vb)\cdot \grad)\vw}_{(c)}+\frac{1}{4}\underbrace{\rot(\vu+\vb)}_{(d)}\bigg)\bigg|,
\end{equation}
and we will study each of the previous term separately. Indeed, for the term (a) above, proceeding in a similar fashion as in (\ref{EstimationPonctuelleVV2}), we have (for $1\leq i\leq 3$):
\begin{align*}
\left|\mathds{1}_{Q_{R_2}}\int_0^t\int_{\mathbb{R}^3}\mathfrak{g_{t-s}}(x-y)[\bar\phi \partial _i div(\vw)](s,y)dyds \right|&\le C\mathds{1}_{Q_{R_2}}\left(\mathcal{L}_1(|\mathds{1}_{Q_{\bar R}}div(\vw)|)(t,x)+\mathcal{L}_2(|\mathds{1}_{Q_{\bar R}}div(\vw)|)(t,x)\right).
\end{align*}
Then, we only weed to study the quantities in the right-hand side above: $\|\mathds{1}_{Q_{R_2}}\mathcal{L}_{1}(|\mathds{1}_{Q_{\bar R}}div(\vw)|)\|_{\M^{3,\sigma}}$ and $\|\mathds{1}_{Q_{R_2}}\mathcal{L}_{2}(|\mathds{1}_{Q_{\bar R}}div(\vw)| )\|_{\M^{3,\sigma}}$. For the first term we fix $p=2$, $q=\frac{10}3$ and $\lambda=\frac13$, we thus have $\frac p\lambda\geq3$ and $\frac{q}{\lambda}=10\geq \sigma$ and by Lemma \ref{lemma_locindi} and by Lemma \ref{Lemme_Hed} we have $\|\mathds{1}_{Q_{R_2}}\mathcal{L}_{1}(\mathds{1}_{Q_{\bar R}}div(\vw))\|_{\M^{3,\sigma}}\leq C\|\mathcal{L}_{1}(|\mathds{1}_{Q_{\bar R}}div(\vw)|)\|_{\M^{\frac{p}{\lambda},\frac{q}{\lambda}}}\leq C\|\mathds{1}_{Q_{\bar R}}div(\vw)\|_{\M^{p,q}}\leq C\|\mathds{1}_{Q_{\bar R}}div(\vw)\|_{\M^{2,\tau_1}}<+\infty$, since $\tau_1>\tau_0>5$ (and by the Corollary \ref{corolarioMorrey} and its conclusion (\ref{ConlusioncorolarioMorrey})). For the second term we set $p,q=2$ and $\lambda=\frac15$ and by the same arguments we have $\|\mathds{1}_{Q_{R_2}}\mathcal{L}_{2}(|\mathds{1}_{Q_{\bar R}}div(\vw)|)\|_{\M^{3,\sigma}}\leq C\|\mathcal{L}_{2}(|\mathds{1}_{Q_{\bar R}}div(\vw)|)\|_{\M^{\frac{p}{\lambda},\frac{q}{\lambda}}}\leq C\|\mathds{1}_{Q_{\bar R}}div(\vw)\|_{\M^{2,2}}\leq C\|\mathds{1}_{Q_{\bar R}}div(\vw)\|_{\M^{2,\tau_1}}<+\infty$ and thus the term (a) is bounded in the Morrey space $\M^{3,\sigma}$.\\
For the term (b) we proceed just as for the term $\vec{\mathcal{V}}_{1}$ and we have
$$\left\|\mathds{1}_{Q_{R_2}}\int_0^t\int_{\mathbb{R}^3}\mathfrak{g_{t-s}}(\cdot-y)[\bar\phi \vw](s,y)dyds \right\|_{\M^{3,\sigma}}\leq C\|\mathds{1}_{Q_{R_2}}\mathcal{L}_2(|\mathds{1}_{Q_{\bar R}}\vw|)\|_{\M^{3,\sigma}}.$$
Setting $p,q=2$ and $\lambda=\frac15$, we have $\|\mathds{1}_{Q_{R_2}}\mathcal{L}_{2}(|\mathds{1}_{Q_{\bar R}}\vw|)\|_{\M^{3,\sigma}}\leq C\|\mathcal{L}_{2}(|\mathds{1}_{Q_{\bar R}}\vw|)\|_{\M^{\frac{p}{\lambda},\frac{q}{\lambda}}}\leq C\|\mathds{1}_{Q_{\bar R}}\vw\|_{\M^{2,2}}\leq C\|\mathds{1}_{Q_{\bar R}}\vw\|_{\M^{3,\tau_0}}<+\infty$
(since we have (\ref{ConclusionPropositionFirstMorreySpace})).\\
Due to the symmetric information available for the variables $\vu, \vb$ and $\vw$ it is easy to see that the term (c) of (\ref{DernierTermeN}) can be treated as the term $\vec{\mathcal{V}}_{3}$ while the term (d) of (\ref{DernierTermeN}) can be studied as $\vec{\mathcal{V}}_{2}$.\\
With all these remarks we finally obtain that $\|\mathds{1}_{Q_{R_2}}\vec{\mathcal{V}}_{24}\|_{\M^{3,\sigma}}<+\infty$.
\end{itemize}
With all these estimates Lemma \ref{LemmaGainenMorreySigma} is now proven. \hfill$\blacksquare$\\

\noindent{\bf End of the proof of Proposition \ref{MorreyPropo2}.} We have proven that $\mathds{1}_{Q_{R_2}} (\vu+\vb+\vw) \in \mathcal{M}_{t,x}^{3,\sigma}$ for $\tau_0<\sigma$ with $\sigma$ very close to $\tau_0$ (say $\sigma=\tau_0+\epsilon$). But this is not enough to ensure the condition $\frac{1}{\delta}+\frac{1}{\tau_{0}} < \frac{1-\alpha}{5}$ stated in Proposition \ref{MorreyPropo2}. In order to obtain this relationship, we will iterate the arguments above. Indeed, considering the information $\mathds{1}_{Q_{R_2}} (\vu+\vb+\vw) \in \mathcal{M}_{t,x}^{3, \tau_0+\epsilon}$ and reapplying Lemma \ref{LemmaGainenMorreySigma}, we will obtain $\mathds{1}_{Q_{\bar R_2}} (\vu+\vb+\vw) \in \mathcal{M}_{t,x}^{3, \sigma_1}$ where $\bar R_2<R_2$ and $\sigma_1=\sigma+\epsilon=\tau_0+2\epsilon$ and we can repeat these arguments until obtaining $\mathds{1}_{Q_{\bar {\bar{R_2}}}} (\vu+\vb+\vw) \in \mathcal{M}_{t,x}^{3, \sigma_n}$ where $\sigma_n=\tau_0+(n+1) \epsilon$ such that $\frac{1}{\sigma_n}+\frac{1}{\tau_0}<\frac{1-\alpha}{5}$ with $\bar {\bar{R_2}}<\bar{R_2}$. As we can see, at each iteration we have to consider smaller parabolic balls and without fear of confusion we can set $\delta=\sigma_n$ with the corresponding radius to be $R_2$. We thus have 
$\mathds{1}_{Q_{ R_2}} \vu \in \mathcal{M}_{t,x}^{3, \delta}$ and $\mathds{1}_{Q_{ R_2}} \vb \in \mathcal{M}_{t,x}^{3, \delta}$ with $\frac{1}{\delta}+\frac{1}{\tau_0}<\frac{1-\alpha}{5}$ and the proof of  Proposition \ref{MorreyPropo2} is finished. \hfill $\blacksquare$
\begin{rema}\label{RemarqueDelta}
Note that by iteration the value of $\delta$ can be made big enough.
\end{rema}
We have obtained the hypotheses 1), 2), 4) of the Proposition \ref{HolderRegularityproposition} and with these results at hand we will now study the hypothesis 5). 
%%%%%%%%%%%%%%%%%%%%%%%%%%%%%%%%%%%%%%
 \begin{coro}\label{CorolarioGainParabolique1}
Consider the general hypotheses of Theorem  (\ref{HypothesesLocal1}). Then, for $R_2$ such that $R_2<R_1<{\bf R}$ and for $1\leq i,j\leq 3$ we have 
$$\mathds{1}_{Q_{R_2}}\frac{\partial_i \partial_j}{(-\Delta)}(u_i b_j)\in \M^{\p,\q},\quad \mathds{1}_{Q_{R_2}}\frac{\partial_i \partial_j}{(-\Delta)}(u_i a_j)\in \M^{\p,\q}, \quad\mbox{ and }\quad  \mathds{1}_{Q_{R_2}}\frac{\partial_i \partial_j}{(-\Delta)}(b_i a_j)\in \M^{\p,\q},$$
with $\p_0\leq \p<+\infty$ and $\q_1\leq \q<+\infty$ where $1\le \mathfrak{p}_0\le \frac{6}{5}$ and $5<\q_1=\frac{5}{1-\alpha}<\frac{11}{2}$ with $0<\alpha<\frac{1}{12}$.
\end{coro}
{\bf Proof.} Recall that from (\ref{FormulePressionIntro}) we have the expression $p = \displaystyle{\sum_{i,j=1}^3}\frac{\partial_i\partial_j}{2(-\Delta)}\left(u_i b_j+(u_i+b_i)a_j+a_i(u_j+b_j)\right)$, which corresponds with the terms that we want to study and consequently we only need to prove that we have $\mathds{1}_{Q_{R_2}}p\in  \M^{\p,\q}$. Thus introducing suitable localizing functions $\bar\phi$ and $\bar\psi$ as in (\ref{ProprieteLocalisation6}) and following the computations made in (\ref{ExpressionPression6}), (\ref{FormulaHipotesis5}) and (\ref{FormulaHipotesis51}) we have
\begin{eqnarray}
\bar \phi p&=&\bigg(\bar \phi\frac{-\bar \psi\Delta p}{(-\Delta)}\bigg)+\bar \phi\frac{(\Delta \bar \psi)p}{(-\Delta)}-2\sum_{i=1}^3\bar \phi\frac{ \partial_i((\partial_i\bar \psi)p)}{(-\Delta)}\notag\\
&=&\sum_{i,j=1}^3 \frac{\bar \phi}{2(-\Delta)}\biggl(\underbrace{\partial_i\partial_j(\bar \psi u_i b_j)}_{(i)}+\underbrace{(\partial_i \partial_j \bar \psi)(u_i b_j)}_{(ii)}-[\underbrace{\partial_i((\partial_j \bar \psi)u_i b_j)+\partial_j((\partial_i \bar \psi)u_i b_j)}_{(iii)}]\label{EstimationsFinales}\\
&&+\underbrace{\partial_i\partial_j(\bar \psi (u_i+b_i) a_j)}_{(i)}+\underbrace{(\partial_i \partial_j \bar \psi)((u_i+b_i) a_j)}_{(ii)}-[\underbrace{\partial_i((\partial_j \bar \psi)(u_i +b_i)a_j)\partial_j((\partial_i \bar \psi)(u_i+b_i) a_j)}_{(iii)}]\notag\\
&&+\underbrace{\partial_i\partial_j(\bar \psi  a_i(u_j+b_j))}_{(i)}+\underbrace{(\partial_i \partial_j \bar \psi)( a_i(u_j+b_j))}_{(ii)}-[\underbrace{\partial_i((\partial_j \bar \psi)a_i(u_j +b_j))+\partial_j((\partial_i \bar \psi) a_i(u_j+b_j))}_{(iii)}]\biggr)\notag\\
&&+\underbrace{\bar \phi\frac{(\Delta \bar \psi)p}{(-\Delta)}}_{(iv)}-2\sum_{i=1}^3\underbrace{\bar \phi\frac{ \partial_i((\partial_i\bar \psi)p)}{(-\Delta)}}_{(v)},\notag
\end{eqnarray}
and we will prove that each one of these terms belong to the space $\M^{\frac65, \frac{11}{2}}$ (we are considering here $\p=\frac65$ and $\q= \frac{11}{2}$). Fortunately, many terms of (\ref{EstimationsFinales}) share a common structure.
\begin{itemize}
\item For the term of the form $(i)$ we write:
$$\bigg\|\frac{\bar \phi}{(-\Delta)}\partial_i\partial_j(\bar \psi u_i b_j)\bigg\|_{\M^{\frac65, \frac{11}{2}}}\leq C\|\mathds{1}_{Q_{R_2}} u_i b_j\|_{\M^{\frac65, \frac{11}{2}}}\leq C\|\mathds{1}_{Q_{R_2}} u_i\|_{\M^{3, \delta'}}\|\mathds{1}_{Q_{R_2}} b_j\|_{\M^{3, \delta'}}<+\infty,$$
where we used the boundedness of the Riesz transforms in Morrey spaces as well as the H\"older inequalities (and we considered $\delta'=11$ which is possible by Remark \ref{RemarqueDelta}). We consider now the terms of the form $ \frac{\bar \phi}{(-\Delta)}\partial_i\partial_j(\bar \psi u_i a_j)$ and we write by the same arguments as above
$$\bigg\|\frac{\bar \phi}{(-\Delta)}\partial_i\partial_j(\bar \psi u_i a_j)\bigg\|_{\M^{\frac65, \frac{11}{2}}}\leq C\|\mathds{1}_{Q_{R_2}} u_i a_j\|_{\M^{\frac65, \frac{11}{2}}}\leq C\|\mathds{1}_{Q_{R_2}} u_i\|_{\M^{3, \delta''}}\|\mathds{1}_{Q_{R_2}} a_j\|_{\M^{6,6}}<+\infty,$$
where $\delta''=66$.
\item For the terms of the form $(ii)$, we first have
$$\left\|\frac{\bar\phi }{(-\Delta )}(\partial_i \partial_j \bar\psi) (u_i b_j)\right\|_{\mathcal{M}^{\frac{6}{5}, \frac{11}{2}}_{t,x}}\leq C\left\|\frac{\bar\phi}{(-\Delta )}(\partial_i \partial_j \bar\psi) (u_i b_j)\right\|_{\mathcal{M}^{\frac{11}{5}, \frac{11}{2}}_{t,x}}\leq C\left\| \frac{\bar\phi}{(-\Delta )}(\partial_i \partial_j  \bar\psi) (u_i b_j)\right\|_{L^{\frac{11}{5}}_t L^{\infty}_{x}},$$
where we used the space inclusion $L^{\frac{11}{5}}_t L^{\infty}_{x}\subset \mathcal{M}^{\frac{11}{5}, \frac{11}{2}}_{t,x}$. Following the same ideas displayed in formulas (\ref{Formula_intermediairevVV80})-(\ref{Formula_intermediairevVV81}), due to the support properties of the auxiliary functions we obtain
$$\left\| \frac{\bar\phi}{(-\Delta )}(\partial_i \partial_j \bar\psi) (u_i b_j)\right\|_{L^{\frac{11}{5}}_t L^{\infty}_{x}}\leq \|\mathds{1}_{Q_{R_2}}u_i b_j\|_{L^{\frac{11}{5}}_{t}L^1_x}\leq C\|\mathds{1}_{Q_{R_2}}\vu\|_{L^{\frac{22}{5}}_tL^{\frac{66}{23}}_{x}}\|\mathds{1}_{Q_{R_2}}\vb\|_{L^{\frac{22}{5}}_tL^{\frac{66}{23}}_{x}}<+\infty,$$
as by interpolation we have $\|\vu\|_{L^{\frac{22}{5}}_tL^{\frac{66}{23}}_{x}(\Omega)}\leq \|\vu\|_{L^{\infty}_tL^{2}_{x}(\Omega)}^{\frac{5}{11}}\|\vu\|_{L^{2}_tL^{6}_{x}(\Omega)}^{\frac{6}{11}}$. The terms of the form $\frac{\bar\phi }{(-\Delta )}(\partial_i \partial_j \bar \psi) (u_i a_j)$ are treated in exactly the same fashion as we have $\|\va\|_{L^{\frac{22}{5}}_tL^{\frac{66}{23}}_{x}(\Omega)}\leq C\|\va\|_{L^6_{t,x}(\Omega)}$.
\item The term of the form $(iii)$ can be studied in exactly the same manner as the terms of the form $(ii)$. 
\item For the term $(iv)$, by the same arguments we obtain
$$\bigg\|\bar \phi\frac{(\Delta \bar \psi)p}{(-\Delta)}\bigg\|_{\mathcal{M}^{\frac{6}{5}, \frac{11}{2}}_{t,x}}\leq\left\| \frac{\bar\phi}{(-\Delta )}(\Delta \bar\psi) p\right\|_{L^{\frac{11}{5}}_t L^{\infty}_{x}}\leq C\|\mathds{1}_{Q_{R_2}}p\|_{L^{\frac{11}{5}}_{t}L^1_x}\leq C\|\mathds{1}_{Q_{R_2}}p\|_{L^{\frac{5}{2}}_{t}L^1_x}<+\infty.$$
\item The term $(v)$ can be treated in the same manner as the previous point. 
\begin{rema}
The condition $p\in L^{\frac{5}{2}}_{t}L^1_x(\Omega)$ is needed here in order to treat these two previous terms. If we have some additional information over the perturbation term (\emph{e.g.} $\va\in L^2_t\dot{H}^1_x(\Omega)$) then these terms can be controlled by the information $p\in L^{\frac{3}{2}}_{t,x}(\Omega)$.
\end{rema}

\end{itemize}
We have proven so far that all the terms of (\ref{EstimationsFinales}) can be controlled in the Morrey space $\M^{\frac65, \frac{11}{2}}$ and this ends the proof of the Corollary \ref{CorolarioGainParabolique1}. \hfill $\blacksquare$\\
%%%%%%%%%%%%%%%%%%%%%%%%%%%%%%%%%%%%%%
 
\noindent In order to obtain Proposition \ref{HolderRegularityproposition} (and thus Theorem \ref{HolderRegularity_theorem}) we only need to verify the hypothesis 3) \emph{i.e.} $\mathds{1}_{Q_{R_2}}div(\vw)\in \M^{\frac{6}{5}, \frac{15}{2}}$. Recall that in the Corollary \ref{corolarioMorrey} we have obtained that $\mathds{1}_{Q_{R_1}}div(\vw)\in \mathcal{M}^{2,\tau_1}_{t,x},$ but this is not enough to our purposes. In order to treat this condition we have:
\begin{propo}\label{divw}
Under the general hypotheses of Theorem \ref{HolderRegularity_theorem} we have $\mathds{1}_{Q_{R_2}}div(\vw)\in \M^{\frac{6}{5}, \frac{15}{2}}$.
\end{propo}
\textbf{Proof.} We first apply the divergence operator in the equation satisfied by $\vw$ (see (\ref{EquationMMP})) to obtain
$$\partial_t div(\vw)=2\Delta div(\vw)-div(\vw)-\frac{1}{2}div(((\vu+\vb)\cdot\grad)\vw).$$
Considering the localizing function $\bar\phi$ as in (\ref{ProprieteLocalisation6}) if we define $\mathcal{W}=\bar \phi div(\vw)$ we obtain the system $\partial_t \mathcal{W}=2\Delta \mathcal{W}+\mathbb{W}$ with $\mathcal{W}(0, \cdot)=0$ where $\mathbb{W}=(\partial_t \bar \phi-2\Delta \bar\phi-\bar\phi )div(\vw)-4\sum_{i=1}^3(\partial_i \bar\phi)(\partial_i div(\vw))-\frac{1}{2}\bar\phi div((\vu\cdot \grad)\vw)$, and we have
\begin{equation}\label{Ws}
\mathcal{W}(t,x)=\int_0^t e^{2(t-s)\Delta}\bigg(\underbrace{(\partial_t \bar \phi-2\Delta \bar\phi-\bar\phi )div(\vw)}_{(\mathcal{W}_1)}-4\sum_{i=1}^3\underbrace{(\partial_i \bar\phi)(\partial_i div(\vw))}_{(\mathcal{W}_2)}-\frac{1}{2}\underbrace{\bar\phi div((\vu\cdot \grad)\vw)}_{(\mathcal{W}_3)}\bigg)ds
\end{equation}
Now we will prove that each one of these term belong to $\M^{\frac{6}{5}, \frac{15}{2}}$. Indeed:
\begin{itemize}
\item For the first term $\mathcal{W}_1$ we write, following the same arguments as in (\ref{EstimationPonctuelleVV2}):
\begin{equation}\label{first}
\|\mathds{1}_{Q_{R_2}}\mathcal{W}_1\|_{\M^{\frac{6}{5}, \frac{15}{2}}}\le C\big(\|\mathds{1}_{Q_{R_2}}\mathcal{L}_1(|\mathds{1}_{Q_{\bar R}}\vw|) \|_{\M ^{\frac{6}{5}, \frac{15}{2}}}+\|\mathds{1}_{Q_{R_2}}\mathcal{L}_2(|\mathds{1}_{Q_{\bar R}}\vw|) \|_{\M^{\frac{6}{5}, \frac{15}{2}}}\big).
\end{equation}
For the first term above we set $p=\frac65$, $q=\frac{9}{2}$ and $\lambda=\frac{1}{10}$ and by Lemma \ref{Lemme_Hed} we obtain
\begin{equation*}
\|\mathds{1}_{Q_{R_2}}\mathcal{L}_1(|\mathds{1}_{Q_{\bar R}}\vw|)\|_{\M ^{\frac{6}{5},\frac{15}{2}}}\leq C\|\mathds{1}_{Q_{R_2}}\mathcal{L}_1(|\mathds{1}_{Q_{\bar R}}\vw|)\|_{\M ^{\frac{6}{\lambda 5},\frac{9}{2\lambda}}}\le C\|\mathds{1}_{Q_{R_1}}\vw \|_{\M ^{\frac{6}{5},\frac{9}{2}}}\le C \|\mathds{1}_{Q_{R_1}}\vw \|_{\M ^{3,\tau_0}}<+\infty.
\end{equation*}
For the second term of (\ref{first}), we fix $p=\frac65$, $q=\frac{12}{5}$ and $\lambda=\frac{1}{25}$, thus by Lemma \ref{Lemme_Hed}, we have 
\begin{equation*}
\|\mathds{1}_{Q_{R_2}}\mathcal{L}_2(|\mathds{1}_{Q_{\bar R}}\vw|)\|_{\M ^{\frac65,\frac{15}{2}}}\le \|\mathds{1}_{Q_{R_2}}\mathcal{L}_2(|\mathds{1}_{Q_{\bar R}}\vw|)\|_{\M ^{\frac{6}{\lambda 5},\frac{12}{\lambda 5}}}\leq C\|\mathds{1}_{Q_{R_1}}\vw \|_{\M ^{\frac{6}{5},\frac{12}{5}}}\le C \|\mathds{1}_{Q_{R_1}}\vw \|_{\M ^{3,\tau_0}}<+\infty.
 \end{equation*}
\item For the term $\mathcal{W}_2$ of (\ref{Ws}) we write for all $1\le i\le 3$:
\begin{eqnarray}
\left|\int_0^t\int_{\mathbb{R}^3}\mathfrak{g_{t-s}}(x-y)[ (\partial_i \bar \phi)\partial _i div(\vw) (s,y)]dyds \right|\leq\bigg|\mathds{1}_{Q_{R_2}}\int_0^t\int_{\mathbb{R}^3}\partial _i\mathfrak{g_{t-s}}(x-y)[(\partial_i \bar \phi ) div(\vw) (s,y)]dyds\bigg|\notag\\
+\bigg|\mathds{1}_{Q_{R_2}}\int_0^t\int_{\mathbb{R}^3}\mathfrak{g_{t-s}}(x-y)[(\partial _i^2 \bar \phi)  div(\vw) (s,y)]dyds\biggl|.\qquad \label{div}
\end{eqnarray}
The last term above can be studied just as $\mathcal{W}_2$ while for the first term of the right-hand side of (\ref{div}) we write:
\begin{equation*}
\left\|\mathds{1}_{Q_{R_2}}\int_0^t\int_{\mathbb{R}^3}\partial_i\mathfrak{g_{t-s}}(\cdot-y)[(\partial_i \bar \phi ) div(\vw) (s,y)]dyds\right\|_{\M ^{\frac65,\frac{15}{2}}}\le C\left\|\mathds{1}_{Q_{R_2}} \mathcal{L}_2(|\mathds{1}_{Q_{\bar R}}  div(\vw)|)\right\|_{\M ^{\frac65,\frac{15}{2}}}.
\end{equation*}
Taking again $p=\frac65$, $q=\frac{12}{5}$ and $\lambda=\frac{1}{25}$, applying Lemma \ref{Lemme_Hed}, we obtain
\begin{equation}\label{lastLastEstimate}
\left\|\mathds{1}_{Q_{R_2}} \mathcal{L}_2(|\mathds{1}_{Q_{\bar R}}  div(\vw)|)\right\|_{\M ^{\frac65,\frac{15}{2}}}\le C\|\mathds{1}_{Q_{R_1}}div(\vw) \|_{\M ^{\frac{6}{5},\frac{12}{5}}}\le \|\mathds{1}_{Q_{R_1}}div(\vw) \|_{\M ^{2,\tau_1}}<+\infty.
\end{equation}
\item Now, we study the last term of (\ref{Ws}). For $1\le i \le 3$ we write
\begin{align*}
\left|\mathds{1}_{Q_{R_2}}\int_0^t\int_{\mathbb{R}^3}\mathfrak{g_{t-s}}(x-y)\bar\phi\partial_i[(\vu\cdot \grad)\vw)  ]dyds\right|&\le \underbrace{\left|\mathds{1}_{Q_{R_2}}\int_0^t\int_{\mathbb{R}^3}\partial_i\mathfrak{g_{t-s}}(x-y)\bar\phi[(\vu\cdot \grad)\vw)  ]dyds\right|}_{(a)}\\
&+\left|\mathds{1}_{Q_{R_2}}\int_0^t\int_{\mathbb{R}^3}\mathfrak{g_{t-s}}(x-y)(\partial_i\bar\phi)[(\vu\cdot \grad)\vw)]dyds\right|,
\end{align*}
and we remark that in order to study the last term above it is enough to consider, for $1\le l,j,k \le 3$, the quantities
\begin{align*}
\underbrace{\left|\mathds{1}_{Q_{R_2}}\int_0^t\int_{\mathbb{R}^3}\partial_l\mathfrak{g_{t-s}}(x-y)(\partial_i\bar\phi) (u_j\omega_k)dyds\right| }_{(b)}\qquad \mbox{and}\quad\underbrace{\left|\mathds{1}_{Q_{R_2}}\int_0^t\int_{\mathbb{R}^3}\mathfrak{g_{t-s}}(x-y)(\partial_l\partial_i\bar\phi) (u_j\omega_k)dyds\right|}_{(c)}.
\end{align*}
Following the computations performed above, we have for the term $(a)$:
$$\left\|\mathds{1}_{Q_{R_2}}\int_0^t\int_{\mathbb{R}^3}\partial_i\mathfrak{g_{t-s}}(\cdot-y)\bar\phi[(\vu\cdot \grad)\vw)]dyds\right\|_{\M ^{\frac65,\frac{15}{2}}}\leq C\|\mathds{1}_{Q_{R_2}}\mathcal{L}_2(|\mathds{1}_{Q_{\bar R}}(\vu\cdot \grad)\vw|)\|_{\M ^{\frac65,\frac{15}{2}}},$$
and fixing $p=\frac65$, $q=\frac{10\tau_0}{10+3\tau_0}$ and $\lambda=\frac{10-\tau_0}{10+3\tau_0}$, by Lemma \ref{Lemme_Hed} we have
$$\|\mathds{1}_{Q_{R_2}}\mathcal{L}_2(|\mathds{1}_{Q_{\bar R}}(\vu\cdot \grad)\vw|)\|_{\M ^{\frac65,\frac{15}{2}}}\leq C\|\mathds{1}_{Q_{\bar R}}(\vu\cdot \grad)\vw\|_{\M ^{\frac65,\frac{10\tau_0}{10+3\tau_0}}}\leq C\|\mathds{1}_{Q_{\bar R}}\vu\|_{\M ^{3,\tau_0}}\|\mathds{1}_{Q_{\bar R}} \grad\otimes\vw\|_{\M ^{2,\tau_1}}<+\infty,$$
where we used the H\"older inequalities for Morrey spaces and Lemma \ref{lemma_locindi} in the last estimate.\\

For the term $(b)$, we write
$$\left\|\mathds{1}_{Q_{R_2}}\int_0^t\int_{\mathbb{R}^3}\partial_l\mathfrak{g_{t-s}}(\cdot-y)(\partial_i\bar\phi) (u_j\omega_k)dyds\right\|_{\M ^{\frac65,\frac{15}{2}}}\leq C\|\mathds{1}_{Q_{R_2}}\mathcal{L}_2(|\mathds{1}_{Q_{\bar R}}u_j\omega_k|)\|_{\M ^{\frac65,\frac{15}{2}}},$$
and applying the same arguments as in (\ref{lastLastEstimate}) we have
\begin{eqnarray*}
\|\mathds{1}_{Q_{R_2}}\mathcal{L}_2(|\mathds{1}_{Q_{\bar R}}u_j\omega_k|)\|_{\M ^{\frac65,\frac{15}{2}}}\leq C\|\mathds{1}_{Q_{\bar R}}u_j\omega_k\|_{\M ^{\frac65,\frac{12}{5}}}\leq C\|\mathds{1}_{Q_{\bar R}}u_j\omega_k\|_{\M ^{\frac32,10}}\\
\leq C\|\mathds{1}_{Q_{R_1}}\vu \|_{\M ^{3,\delta}}\|\mathds{1}_{Q_{R_1}}\vw \|_{\M ^{3,\delta}}<+\infty.
\end{eqnarray*}
For the term $(c)$, by the same ideas displayed in the study of first term of (\ref{first}) we have:
\begin{eqnarray*}
\left\|\mathds{1}_{Q_{R_2}}\int_0^t\int_{\mathbb{R}^3}\mathfrak{g_{t-s}}(\cdot-y)(\partial_l\partial_i\bar\phi) (u_j\omega_k)dyds\right\|_{\M ^{\frac65,\frac{15}{2}}}\leq C\|\mathds{1}_{Q_{R_2}}\mathcal{L}_1(|\mathds{1}_{Q_{\bar R}}u_j\omega_k|)\|_{\M ^{\frac65,\frac{15}{2}}}\leq C\|\mathds{1}_{Q_{\bar R}}u_j\omega_k\|_{\M ^{\frac65,\frac{9}{2}}}\\
\leq C\|\mathds{1}_{Q_{\bar R}}u_j\omega_k\|_{\M ^{\frac32,10}}\leq C\|\mathds{1}_{Q_{R_1}}\vu \|_{\M ^{3,\delta}}\|\mathds{1}_{Q_{R_1}}\vw \|_{\M ^{3,\delta}}<+\infty.
\end{eqnarray*}
\end{itemize}
We have thus proven that all the terms of (\ref{Ws}) belong to the Morrey space $\M^{\frac{6}{5}, \frac{15}{2}}$: the proof of the Proposition \ref{divw} is complete. \hfill$\blacksquare$\\

We have now all the hypotheses of the Proposition \ref{HolderRegularityproposition}, and thus Theorem \ref{HolderRegularity_theorem} follows.\hfill$\blacksquare$\\

\noindent{\bf Acknowledgements.} We would like to thank Jiao He for fruitful discussions. 
%%%%%%%%%%%%%%%%%%%%%%%%%%%%%%%%%%%%%%

%%%%%%%%%%%%%%%%%%%%%%%%%%%%%%%%%%%%%%
\end{document}